\documentclass[11pt]{amsart}

\usepackage[OT2, T1]{fontenc}
\usepackage{url}
\usepackage{booktabs}
\usepackage{float}
\usepackage{amsmath}
\usepackage{stackengine}
\usepackage{array}
\usepackage{graphicx}
\usepackage{amsfonts}
\usepackage{amssymb}
\usepackage{amstext}
\usepackage{amsthm}
\usepackage{enumitem}
\usepackage{bm}
\usepackage{hyperref}
\usepackage{colonequals}
\usepackage{enumitem}
\usepackage[alphabetic,lite]{amsrefs}
\usepackage{cleveref}
\usepackage[all,cmtip]{xy}
\usepackage{fullpage}
\usepackage{appendix}
\usepackage{scalerel}
\usepackage{mathrsfs}
\usepackage{comment}
\usepackage{tikz-cd}
\usepackage{tikz}
\usepackage{galois}
\usepackage{dynkin-diagrams}
\usepackage{amssymb} 
\usepackage{subcaption}
\def\acts{\curvearrowright}
\numberwithin{equation}{subsection}

\newtheorem{theorem}{Theorem}
\newtheorem{lemma}[theorem]{Lemma}
\newtheorem{proposition}[theorem]{Proposition}

\newtheorem{conj}[theorem]{Conjecture}

\theoremstyle{definition}
\newtheorem{construction}{Construction}
\newtheorem{defn}[theorem]{Definition}

\theoremstyle{remark}
\renewcommand{\comment}[1]{}
\newtheorem{remark}[theorem]{Remark}

\newtheorem*{claim}{Claim}
\newtheorem{example}[theorem]{Example}
\newtheorem{notation}[theorem]{Notation}

\numberwithin{theorem}{section}

\newcommand{\bA}{\mathbb{A}}

\newcommand{\bC}{\mathbb{C}}
\newcommand{\bD}{\mathbb{D}}

\newcommand{\bF}{\mathbb{F}}
\newcommand{\bG}{\mathbb{G}}
\newcommand{\bH}{\mathbb{H}}

\newcommand{\bK}{\mathbb{K}}

\newcommand{\bN}{\mathbb{N}}

\newcommand{\bQ}{\mathbb{Q}}

\newcommand{\bZ}{\mathbb{Z}}

\newcommand{\cA}{\mathscr{A}}

\newcommand{\cC}{\mathcal{C}}

\newcommand{\cM}{\mathcal{M}}
\newcommand{\cN}{\mathcal{N}}
\newcommand{\cO}{\mathcal{O}}

\newcommand{\Sh}{\mathcal{S}}

\newcommand{\et}{{\text{\'et}}}

\newcommand{\cris}{{\mathrm{cris}}}

\newcommand{\la}{\langle}
\newcommand{\ra}{\rangle}

\DeclareMathOperator{\FIsoc}{\textbf{$F$-Isoc}}

\newcommand{\Fpbar}{\bF}
\DeclareMathOperator{\GL}{GL}
\DeclareMathOperator{\SL}{SL}
\DeclareMathOperator{\GSp}{GSp}

\DeclareMathOperator{\Gal}{Gal}
\DeclareMathOperator{\End}{End}
\DeclareMathOperator{\Hom}{Hom}

\DeclareMathOperator{\Lie}{Lie}

\DeclareMathOperator{\Spec}{Spec}

\DeclareMathOperator{\ad}{ad}

\DeclareMathOperator{\der}{der}
\DeclareMathOperator{\Fil}{Fil}

\DeclareMathOperator{\diag}{diag}

\DeclareMathOperator{\Id}{Id}

\DeclareMathOperator{\Tr}{Tr}

\DeclareMathOperator{\gr}{gr}

\DeclareMathOperator{\Vect}{Vect}

\DeclareMathOperator{\CH}{CH}

\newcommand{\ol}[1]{\overline{#1}}

\begin{document}

\title{The Tate conjecture for abelian fourfolds}
\begin{abstract}
    We prove the Tate conjecture for non-CM and non-ordinary abelian fourfolds over positive characteristic function fields. Combining our previous works \cites{Jiang2025,Broe2026}, this establishes the Tate conjecture for abelian fourfolds over finitely generated fields. Along the way, we classify the $l$-adic monodromy of abelian fourfolds in characteristic $p$ via mixed characteristic lifts and methods from convergent and overconvergent crystalline monodromy, and show that Zarhin's minuscule weights conjecture holds for abelian varieties up to dimension 4. These enable us to prove a characteristic $p$ variant of the Mumford--Tate conjecture for most abelian fourfolds, and reduce the Tate conjecture to Markman's results on the Hodge conjecture. 
\end{abstract}
\author{Matt Broe and Ruofan Jiang}
\maketitle

\tableofcontents
\section{Introduction}
Let $L$ be a finitely generated field, and let $A/L$ be an abelian fourfold. Let $l$ be a prime different from $\mathrm{char}\,L$. The Tate conjecture predicts that for each $i\geq 0$ the Galois invariant subspace
$$[H^{2i}_{\text{ét}}(A_{\overline{L}},\bQ_l)(i)]^{\Gal(\overline{L}/L)}$$
is the $\bQ_l$-span of codimension-$i$ algebraic cycles on $A$. In this paper, we prove 
\begin{theorem}\label{thm:A}
 The Tate conjecture holds for $A$.  
\end{theorem}
The only case that requires consideration is $i=2$. In proving the conjecture for $A$, we may without loss of generality first replace $A$ by an isogenous abelian variety, or base change to a finite extension of $L$.

The Tate conjecture can be viewed as an arithmetic variant of the Hodge conjecture. Both conjectures concern whether certain spaces of cohomology classes are algebraic, and they are among the most well-known examples of \textit{cycle conjectures}. While the Hodge conjecture is stated for varieties over $\bC$, the Tate conjecture is stated for varieties over finitely generated fields (e.g., number fields or function fields over $\mathbb{F}_p$). Arithmetic properties of the base field cause subtle differences in the flavor of the Tate conjecture over different fields.  

\subsection*{Prior work on the Tate and Hodge conjectures}
    The following review of existing literature is adapted from a revised version of \cite{Broe2026}. 
    \par
    Abelian fourfolds have long been a focus of research on the Tate and Hodge conjectures. Pohlmann \cite{Pohlmann1968} presented an example, due to Mumford, of an abelian fourfold over $\bC$ whose Betti cohomology contained exceptional Hodge classes, i.e. Hodge classes which are not contained in the span of products of divisor classes. Weil \cite{Weil1979} generalized Mumford's example by defining what are now called abelian varieties of Weil type. These possess certain distinguished Hodge classes, known as Weil classes, which are exceptional for a generic abelian fourfold of Weil type. Results of Tankeev \cite{Tankeev1982} and Moonen--Zarhin (\cite{MZ1995}, \cite{MZ1999}) later reduced the Hodge conjecture for all abelian varieties of dimension $\le 5$ to the algebraicity of Weil classes on abelian fourfolds of Weil type. The latter statement was finally proved by Markman in \cite[Corollary 1.6.1]{Markman2025}.
    \par
 We now discuss the Tate conjecture.
 For a finitely generated field $L \subset \bC$ and an abelian variety $A$ over $L$, the Tate conjecture for cycles of dimension one or codimension one was famously proved by Faltings; cf. \cite[Theorem 1]{Faltings1984}. This implies the full Tate conjecture for $A$ when $\dim A \leq 3$. 
 \par
 To approach the Tate conjecture for $A$ in other codimensions, a well-known method is to use proven cases of the Hodge conjecture to construct cycles on $A_\bC$, then spread them out to $A$. However, the Hodge conjecture for $A_\bC$ does not automatically imply the Tate conjecture for $A$. Hodge classes are controlled by the Mumford--Tate group of $A_\bC$, whereas Tate classes are controlled by the $l$-adic monodromy group of $A$. In general, the Mumford--Tate conjecture is believed to mediate between these two groups. For an example of this philosophy, when $\dim A = 4$, one can use structural results on the monodromy of $A$ (\cite[Theorem 2.11 and Section 4.1]{MZ1995}, \cite[Theorems 2.19 and 6.1]{Lombardo2016}) to reduce the Tate conjecture for $A$ to the Hodge conjecture for $A_\bC$, proved in \cite[Corollary 1.6.1]{Markman2025}.%
    \par
 The situation in positive characteristic $p$ is very different. Abelian varieties in characteristic $p$ can have more cycles than in characteristic zero (as is already evident from consideration of endomorphisms of supersingular elliptic curves). In this setting, the Tate conjecture for cycles of dimension or codimension one was proved by Tate \cite{Tate1966}, Zarhin \cite{Zarhin1976}, and Mori \cite[Th\'eor\`eme XII.2.5]{MB1985}. For cycles in other codimensions, one may still hope to reduce to the Hodge conjecture, but this approach takes on added subtleties. Notably, monodromy methods along the line of the Mumford--Tate conjecture are no longer effective in general. 
 \par
 Milne showed that the Hodge conjecture for CM abelian varieties implies the Tate conjecture for abelian varieties over finite fields \cite{Milne1999b}. By reworking the method of Lefschetz groups used in \textit{loc.cit.}, and leveraging available results on the Hodge conjecture for abelian fourfolds, he later obtained unconditional cases of the Tate conjecture (\cite{Milne2001}, \cite{Milne2022}). The Tate conjecture for abelian fourfolds over finite fields was proved by the first named author in \cite{Broe2026}, using Markman's results. Other cases of the Tate conjecture for abelian varieties over finite fields were proved in \cite{Tankeev1983}, \cite{Zarhin1991}, \cite{LZ1993}, \cite{Spiess1999}, \cite{DKZ2024}, and \cite{APFV2025}. Finally, the Hodge conjecture for CM abelian varieties (and hence the Tate conjecture for all abelian varieties over finite fields) was proved in \cite{OAI2026}.
    \par
   Over characteristic $p$ function fields, the Tate conjecture for ordinary abelian fourfolds was proved by the second-named author by reducing to Markman's results and invoking a characteristic $p$ analogue of the Mumford--Tate conjecture \cite[Remark 1.18]{Jiang2025}. Given what is known, the case of Theorem \ref{thm:A} that remains to be proved is when $A$ is a non-CM and non-ordinary abelian fourfold over a finitely generated field of characteristic $p$. Theorem~\ref{thm:A} then follows from the following 
   
    \par

\begin{theorem}[Theorem~\ref{thm:Tate_simple}, Theorem~\ref{thm:Tate_nonsimple}] Let $A$ be a non-CM non-ordinary abelian fourfold over a characteristic $p$ function field $L$. Then the Tate conjecture holds for $A$.     
\end{theorem}
The overall strategy can be thought of as a combination and generalization of the finite field case in \cite{Broe2026}
and the ordinary function field case in \cite{Jiang2025}. For a large class of non-CM non-ordinary abelian fourfolds, we prove the Tate conjecture by formulating and proving analogues of the Mumford--Tate conjecture for abelian varieties in characteristic $p$, especially for non-ordinary abelian varieties\footnote{For ordinary abelian varieties, analogues of the Mumford--Tate conjecture are already treated in the previous works \cite{Jiang2023,Jiang2025}.}. This permits us to reduce to the Hodge conjecture for abelian fourfolds, proved by Markman. In the remaining cases, we explicitly study the Tate classes and deduce the Tate conjecture by lifting exceptional classes to Hodge--Weil classes in characteristic 0.

\subsection{Strong and weak Mumford--Tate conjectures in characteristic $p$} Let $A$ be an abelian variety over a finitely generated field $L\subseteq \bC$. The classical Mumford--Tate conjecture can be understood as a comparison between Betti and étale cycles: it predicts that the base change to $\bQ_l$ of the Mumford--Tate group $G_{\mathrm{B}}(A)$ coincides with the neutral component of the $l$-adic étale monodromy group $G_l(A)$. The conjecture is known to connect the Hodge and the Tate conjecture: If the Hodge conjecture and the Mumford--Tate conjecture both hold for $A$, then the Tate conjecture holds. 

\subsubsection{Analogues of the Mumford--Tate conjecture in characteristic $p$} 
Now suppose that $A$ is defined over a finitely generated field $L$ of characteristic $p$. The Mumford--Tate conjecture does not naïvely make sense here, since there is no notion of Hodge structure or Mumford--Tate group. However, in the case where $A$ is ordinary, there is a remedy. Let $\overline{L}$ be the algebraic closure of $L$, and canonically lift $A_{\overline{L}}$ to an abelian scheme $\widetilde{A}$ over the Witt ring $W(\overline{L})$. Fix an embedding $W(\overline{L})\hookrightarrow \bC$, and define the characteristic $p$ analogue of the Mumford--Tate group as $G_{\mathrm{B}}(A):=G_{\mathrm{B}}(\widetilde{A}_{\bC})$.\footnote{The original definition in \cite{Jiang2023} is slightly different from, but equivalent to, the version stated here.} The following direct analogue of the Mumford--Tate conjecture is then expected to hold for $A$.
\begin{conj}[Ordinary Mumford--Tate conjecture; cf. \cite{Jiang2023}]\label{conj:OMT}
We have $G_{\mathrm{B}}(A)\otimes\bQ_l=G_{l}(A)^\circ$.
\end{conj} 
The group $G_l(A)$ on the right-hand side is the $l$-adic étale monodromy when $l\neq p$, and is the overconvergent crystalline monodromy when $l=p$.  In \cite{Jiang2025}, it is proved that the Tate conjecture for $A$ follows from the Hodge conjecture for $\widetilde{A}_{\bC}$ and the Mumford--Tate conjecture for $A$. This is the essential mechanism that enabled the proof of the Tate conjecture when $A$ is an ordinary abelian fourfold. 
\par
When $A$ is not necessarily ordinary, the above construction fails: there is in general no canonical lift of $A$, hence no directly analogous definition of the Mumford--Tate group. In fact, simple examples already obstruct a naïve formulation as in Conjecture~\ref{conj:OMT}: if  
$A$ is supersingular, then $G_l(A)^\circ$ is a rank 1 torus, which can never be the $\bQ_l$-base change of the Mumford--Tate group of an abelian variety in characteristic 0. While it is still (conjecturally) expected that $G_l(A)^\circ$ has a $\bQ$-structure (cf. \cite{LP1995}), such a $\bQ$-structure should not be canonical, and contains little information about characteristic 0 lifts of $A$. Therefore we should not expect it to be useful in connecting the Hodge conjecture to the Tate conjecture in characteristic $p$. To maximally retain the characteristic $0$ connection, we make the following formulation of the strong/weak Mumford--Tate conjecture: 
\begin{defn}\label{def:strong MT}
Let $A$ be an abelian variety defined over a finitely generated field $L$ of characteristic $p$, and let $l\neq p$ be a prime. Let 
$W:=W(\overline{L})$, let $\overline{W}$ be the integral closure of $W$ in $\overline{W[1/p]}$, and fix an isomorphism $\overline{W[1/p]} \cong \bC$. Denote by $\widetilde{A}$ an abelian scheme over $\overline{W}$ whose special fiber is isogenous to $A_{\overline{L}}$. It is known that there is an embedding $
    G_l(A)^\circ\hookrightarrow G_{\mathrm{B}}(\widetilde{A}_{\bC})\otimes {\bQ_l}$, which is canonical up to conjugation: see Lemma \ref{lm:embeding}. We say that 
\begin{enumerate}
    \item $A$ \textbf{satisfies the Mumford--Tate conjecture} if there exists a lift $\widetilde{A}$ of $A$ such that the above embedding is an isomorphism.
    \item  $A$ \textbf{satisfies the weak Mumford--Tate conjecture} if there exists a lift $\widetilde{A}$ such that $G_l(A)^\circ$ contains the derived subgroup of $G_{\mathrm{B}}(\widetilde{A}_{\bC})\otimes {\bQ_l}$.
\end{enumerate}
\end{defn}
Note that ``satisfying the Mumford--Tate/weak Mumford--Tate conjecture'' is only a property. For example, a supersingular abelian variety does not satisfy the Mumford--Tate conjecture, but only the weak Mumford--Tate conjecture. The weak Mumford--Tate conjecture does not always hold either. The counterexamples are much more subtle; see the discussion below. 

  If $A$ satisfies the Mumford--Tate conjecture, then the Tate conjecture for $A$ follows from the Hodge conjecture for $\widetilde{A}_{\bC}$; see Proposition~\ref{prop:relationbetweencycleconjectures}. On the other hand, the weak Mumford--Tate conjecture 
does not directly link the Hodge and the Tate conjectures, but satisfying it implies satisfying Zarhin's minuscule weights conjecture:
\begin{conj}[Zarhin's minuscule weights conjecture, cf. {\cite[Conjecture 0.4.2]{Zarhin1985}}]\label{conj:MinusculeWt} 
Let $A$ be an abelian variety over a finitely generated field $L$, and let
$l \neq \mathrm{char}(L)$. Then each almost simple factor of the derived group
of $G_l(A)^\circ \otimes_{\bQ_l} \overline{\bQ}_l$ is of type $A_m$, $B_m$,
$C_m$ or $D_m$, and the monodromy representation is minuscule: for each irreducible subrepresentation $V$ of the derived group acting on the Tate module $V_l(A) \otimes_{\bQ_l} \overline{\bQ}_l$, the restriction of $V$ to
each almost simple factor of the derived group is a sum of trivial and minuscule representations.
\end{conj}
For further background on the minuscule weights conjecture, see \cite[Section 7.2]{CT2025}. In characteristic 0 the conjecture is a theorem of Pink (cf. \cite[Corollary 5.11]{Pink1998}). In positive characteristic, the conjecture was previously known in dimension at most two, and for ordinary abelian varieties of arbitrary dimension \cite[Remark 0.4.2.1 and Corollary 4.3.1]{Zarhin1985} \footnote{The assumption that $p > 2$ in \cite[Remark 0.4.2.1]{Zarhin1985} may be relaxed using the Tate conjecture for homomorphisms of abelian varieties in characteristic two (\cite{Tate1966}, \cite[Th\'eor\`eme XII.2.5]{MB1985}).}. There exist counterexamples to Conjecture \ref{conj:MinusculeWt} in general, due to B\"ultel (\cite{Bultel2005}, \cite{Bultel2013}, \cite{Bultel2022}), and other counterexamples arising from work of Katz (\cite{Katz1990}, \cite{Katz2004}, \cite[Section 7.3.2 and Appendix A]{CT2025}). In particular, B\"ultel gave counterexamples to the conjecture for abelian sixfolds \cite{Bultel2022}. Counterexamples to the minuscule weights conjecture also do not satisfy the weak Mumford--Tate conjecture.

\subsubsection{Classifying the monodromy}\label{subsub:CtM}
The failure of the minuscule weights conjecture is a warning that the classification of $l$-adic monodromy in characteristic $p$ is a subtle problem in general. However, for abelian varieties in low dimensions, there is still hope that the  Mumford--Tate conjecture holds in sufficiently many cases, so that we can  use it to deduce the Tate conjecture from the Hodge conjecture. Before stating our results on the Mumford--Tate conjecture, we define two families of abelian varieties that are ``harder to deal with than others'':\begin{defn}
Let $A$ be an abelian fourfold over $L$. We say that $A$ is of \textbf{cohomological Mumford type}, if after base-change to $\overline{\bQ}_l$, the monodromy representation of $G_l(A)^{\der}$ is the tensor product of three copies of standard representations of $\SL_2$.
 \end{defn} 
 \begin{defn}
     Let $A$ be an abelian threefold or fourfold over $L$. We say that $A$ is of \textbf{exotic type IV}, if either (1) $A$ is an abelian threefold of Newton polygon $\{0^2,(\frac{1}{2})^2,1^2\}$ and $\End^0(A_{\overline{L}})$ is an imaginary quadratic field, or (2)
     $A$ is a non-ordinary abelian fourfold and $\End^0(A_{\overline{L}})$ is a quartic CM field containing an imaginary quadratic field where $p$ splits. In case (2), it can be shown that $A$ has Newton polygon $\{0^2,(\frac{1}{2})^4,1^2\}$. 
 \end{defn}

\begin{theorem}\label{thm:B}
Let $A$ be an abelian fourfold over $L$  not of cohomological Mumford type. Then $A$ satisfies the weak Mumford--Tate conjecture. Furthermore, $A$ satisfies the Mumford--Tate conjecture in the following cases\footnote{For cohomological Mumford type, we expect that the Mumford--Tate conjecture holds, but we do not know how to prove it. For exotic type IV, the Mumford--Tate conjecture simply does not hold.}:
\begin{enumerate}
\item $A$ is ordinary (\cite{Jiang2025}),
    \item Each simple factor of $A_{\overline{L}}$ is non-CM and not of exotic type IV (Theorem~\ref{thm:MT_simple} and \ref{thm:Mumford_Tate_nonsimple}). 
\end{enumerate}
The minuscule weights conjecture~\ref{conj:MinusculeWt} holds for all abelian varieties over $L$ up to dimension 4. \footnote{For abelian fourfolds not of cohomological Mumford type, Conjecture~\ref{conj:MinusculeWt} follows from the weak Mumford--Tate; for abelian fourfolds of cohomological Mumford type, the $l$-adic monodromy is minuscule by definition. The conjecture for abelian varieties of dimension $\leq 3$ follows from  the non-simple case of dimension $4$.}
\end{theorem}

Since people have good understanding on the Mumford--Tate groups for abelian fourfolds, this theorem implicitly classifies $l$-adic monodromy of abelian fourfolds in characteristic $p$. The reader can also find explicit classifications of $l$-adic monodromy of non-CM simple abelian varieties of dimension $\leq 4$ in \S\ref{sec:simple4} (for fourfolds) and Theorem~\ref{thm:simple_leq 3} (for threefolds). These classification results are too wordy to be stated here.

The challenge of Theorem~\ref{thm:B} lies in two aspects: \begin{itemize}
    \item Finding suitable $\widetilde{A}$ for the characteristic $0$ input. This will be discussed in \S\ref{sec:LE}. 
    \item Identifying the $l$-adic monodromy group of a non-ordinary abelian variety. There are two  prominent difficulties: 
    \begin{enumerate}
        \item Unlike characteristic 0 or ordinary cases, it is not \textit{a priori} known if the monodromy representation is minuscule. Therefore, a wide range of methods based on the classification of minuscule reps are rendered invalid. This causes technical problems, for example, when one tries to classify monodromy for simple abelian varieties of type I; see Example~\ref{ex:minus}. Our strategy is to use the overconvergent crystalline monodromy, which carries out the conceptual route of ``studying the Mumford--Tate problem by unipotents'' discussed in \cite[\S 1.3.1]{Jiang2023}; see \S\ref{sub:p-adicmonodromy} for more details.
\item The exotic behavior for abelian varieties of type IV. The difficulty is of a different flavor: the rank of the center of the $l$-adic monodromy of $A$ may be small (smaller than the rank of the center of any $G_{\mathrm{B}}(\widetilde{A}_\bC)$), so $A$ may not even satisfy the Mumford--Tate conjecture.  
    \end{enumerate}We invite the reader to \S\ref{sub:examples1} for a discussion of challenges we have faced.
    \end{itemize}
\subsubsection{The Tate conjecture after Theorem~\ref{thm:B}} Since the Hodge conjecture is known for abelian fourfolds, Theorem~\ref{thm:B} establishes the Tate conjecture whenever the Mumford--Tate conjecture is satisfied. In the situation where the Mumford--Tate conjecture does not hold (or is not known) for $A$, we establish the Tate conjecture by explicitly studying the Tate classes. For example, when $A$ is of cohomological Mumford type, a similar argument as in \cite[\S 4.1 Type I(1)]{MZ1995} implies that the Tate ring of $A$ is generated by divisor classes. When $A$ has a factor of exotic type IV, the Tate ring of $A$ often contains exceptional classes. If such an $A$ is simple, then we can lift $A$ to characteristic 0 such that the exceptional classes lift to Hodge classes, which were proved to be algebraic by Markman. Finally, when $A$ is non-simple and non-CM but contains CM isogeny factors, we can often use the K\"unneth decomposition and exploit the big monodromy of the non-CM factor to show that the Tate ring of $A$ is generated by divisor classes.

\subsection{Lifting endomorphisms}\label{sec:LE} The hardest case of Theorem~\ref{thm:B} is when $A$ is simple, non-CM and non-ordinary. Before we can say anything meaningful about the monodromy, we need to construct sufficiently many ``good lifts" of $A$. The question is: what should qualify as a ``good lift''? 

The answer seems clear: a ``good lift'' should be a lift whose Mumford--Tate group is as small as possible. However, we do not know how to construct such a lift in general. What we can do is construct a lift such that a reasonably large part of the endomorphism algebra also lifts. In general, an $\widetilde{A}$ only lifting part of the endomorphism algebra is not sufficient for even the weak Mumford--Tate conjecture. However, by a low-dimensional miracle, for abelian fourfolds this is all we need.

The key technical result is the following, which is interesting in its own right: 
\begin{theorem}[$\subseteq$Theorem~\ref{thm:lift_endo_fourfold}]\label{thm:C}
    Let $A$ be a non-CM simple abelian fourfold over an characteristic $p$ algebraically closed field $K$. Then there exists a lift $\widetilde{A}$ over $\overline{W}$ whose special fiber is isogenous to $A$, such that $\End^0(\widetilde{A})=\End^0({A})$. 
\end{theorem}
An experienced reader may observe that the theorem reduces, via Serre--Tate theory, to a problem of lifting $p$-divisible groups, to which one can apply Grothendieck--Messing theory and Lubin--Tate theory. The lifting problem is particularly involved when $D$ is a quaternion algebra (over $\bQ$ or a real quadratic field). 

When $A$ has type IV, we will also need to control the signature of the lift for the application to the Mumford--Tate (for the Tate conjecture alone, there is an alternative approach in this case which does not require control over the signature: see the discussion around Lemma \ref{lemWeilAlg13}). The control over the signature is also treated in Theorem~\ref{thm:lift_endo_fourfold}. The following example explains the necessity (and subtlety) of controlling the signature. \begin{example}\label{ex:lift_signatures}
Suppose that $\End^0({A})$ is an imaginary quadratic field $E$, and let $\widetilde{A}$ be a lift as in Theorem~\ref{thm:C}. If $p$ does not split in $E$, then the $E$-action on $\Lie\widetilde{A}_{\bC}$ can have signature $(2,2)$ or $(1,3)$. Furthermore, both possibilities occur (this is not at all clear \textit{a priori}!). The Mumford--Tate group of the lift differs by a central torus between the two cases: in the $(1,3)$ case we have $G_{\mathrm{B}}(\widetilde{A}_{\bC})=\bG_m\cdot \mathrm{U}(1,3)$, while in the $(2,2)$ case we have $G_{\mathrm{B}}(\widetilde{A}_{\bC})=\bG_m\cdot \mathrm{SU}(2,2)$. On the other hand, it can be shown that the $l$-adic monodromy of $A$ is $\bG_{m,\bQ_l}\cdot \mathrm{SU}(4)_{\bQ_l}$. While either signature yields the weak Mumford--Tate conjecture for $A$, only signature $(2,2)$ yields the Mumford--Tate conjecture. 
\end{example}

\subsection{The role of $p$-adic monodromy}\label{sub:p-adicmonodromy} Having constructed $\widetilde{A}_{\bC}$ as in Theorem~\ref{thm:C}, the next goal is to show that $G_l(A)^\circ=G_{\mathrm{B}}(\widetilde{A}_{\bC})\otimes \bQ_l$ (for the Mumford--Tate conjecture) or $G_l(A)^{\circ,\der}=G_{\mathrm{B}}(\widetilde{A}_{\bC})^{\der}\otimes \bQ_l$ (for the weak Mumford--Tate conjecture). As we have discussed earlier in \S\ref{subsub:CtM}, the absence of minuscule weights, as well as the exotic behavior in type IV makes it extremely subtle.  

Instead of studying the $l$-adic étale monodromy group $G_l(A)$ alone, we also study the overconvergent crystalline monodromy group $G_p(A)$ (called the \textbf{$p$-adic monodromy} in the sequel). To make this precise, we need to spread $A$ out to an abelian scheme $\mathcal{A}$ over a smooth base $X$, and define $G_p(A)$ as the monodromy of the corresponding overconvergent $F$-isocrystal (which is a $\Breve{\bQ}_p$-group). It is known that $G_l(A)\otimes \bC_l$ and $G_p(A)\otimes \bC_p$ are isomorphic (up to neutral components) once we make the identification $\bC_l\simeq \bC_p$ (cf. \cite[Theorem 1.2.1]{D2020a}), so the comparison between $G_l(A)$ and $G_{\mathrm{B}}(\widetilde{A}_\bC)$ can be checked instead with $p$-adic monodromy. We refer the reader to \S\ref{sec:mono p as} for more details.  

The advantage of $p$-adic monodromy is that it makes explicit the information coming from Frobenius. We will explain this from two perspectives. Spread out $A$ to an abelian scheme $\mathcal{A}$ over a smooth base $X$ with constant Newton polygon, and let $x\in X$ be an $\Fpbar$-point.\begin{enumerate}
    \item (Parabolicity) The parabolic $P_{x}\subseteq G_p(A)$ associated to the Newton quasi-cocharacter $\nu_x$ has a specific geometric meaning: it is the monodromy group of the convergent $F$-isocrystal associated to  $\mathcal{A}$. This is a special case of Crew's parabolicity conjecture, proved in \cite{D2020b}. One can then study the structure of $P_{x}$ via the convergent $F$-isocrystal in question, which is closely linked to the deformation theory of the fiber $\mathcal{A}_x$. This is the key idea in a series of works (\cite{Jiang2023,Jiang2025}, as well as \cite{DAddezio2022}). One application of the parabolicity conjecture in our paper is to show that when $A$ is non-CM then $G_p(A)$ has nontrivial unipotent with respect to $\nu_x$; see Theorem~\ref{thm:nonvan}. The non-vanishing result will be combined with Frobenius structure discussed below to show that the unipotent satisfies rigid structural properties. 
    \item (Frobenius structure) Let $\sigma\in \mathrm{Gal}(\Breve{\bQ}_p/\bQ_p)$ be the absolute Frobenius. The monodromy group $G_p(A)$ is equipped with a Frobenius structure $\Phi:G_p(A)^{(\sigma)} \rightarrow G_p(A)$ that comes, via the Tannakian formalism, from the Frobenius structure on the $\sigma$-isocrystal associated to the geometric fiber $\mathcal{A}_x$; see \S\ref{subsub:Frobstructure} for more details. To explain this conceptually, the category of $F$-isocrystals on $X$ is $\bQ_p$-linear, while $G_p(A)$ is only a $\Breve{\bQ}_p$-group. The Frobenius structure $\Phi$ recovers the information that $G_p(A)$ alone fails to capture. When $A$ is non-ordinary, the Frobenius structure $\Phi$ often carries non-trivial information about the structure of the group $G_p(A)$. For example, when $A$ has Newton polygon $\{(\frac{1}{4})^4,(\frac{3}{4})^4\}$, we can show that the Lie algebra of the unipotent radical of the parabolic $P_x$ of $G_p(A)$ is a $\sigma$-isocrystal of slope $1/2$. This forces the unipotent radical to be even-dimensional. See \S\ref{example:1/4-3/4}.
\end{enumerate} 
The above two structural properties of the $p$-adic monodromy allow us to understand $G_p(A)$ in a systematic way. We expect that the method is valuable for studying higher dimensional (especially non-ordinary) abelian schemes in characteristic $p$ as well. 

\subsection{Examples}\label{sub:examples1} We collect some purely characteristic $p$ phenomena that were once challenges when we were attacking the conjecture.

\begin{example}[When minuscule weight is missing]\label{ex:minus} As we have discussed, the $l$-adic monodromy representation of a non-ordinary abelian variety is not known to be minuscule. In the following we get a sense of how things get  complicated when minuscule-ness is absent.

Consider a non-ordinary abelian variety $A$ with endomorphism algebra $\bQ$. To determine its $l$-adic monodromy, we wish to answer the following question: let $\mathfrak{sp}_{2g}\acts V_l$ be the std rep, what are possible reductive subalgebras $\mathfrak{n}\subseteq \mathfrak{sp}_{2g}$ acting irreducibly on $V_l$?  

Let $g=4$. If $\mathfrak{n}\acts V_l$ is minuscule, then it is classically known to be either the std rep of $\mathfrak{sp}_8$ or the $\text{std}^{\otimes 3}$-rep of  $\mathfrak{sl}_2^3$. However, if we don't require it to be minuscule, then it can furthermore be $\mathrm{Sym}^7$ of $\mathfrak{sl}_2$; see \S\ref{subsub:TypeI(1)} for more details. Let $g=3$. If $\mathfrak{n}\acts V_l$ is minuscule, then it is classically known to be geometrically the std rep of $\mathfrak{sp}_6$. If it is not minuscule then there are two more possibilities: $\mathfrak{sl}_2$ acting by $\mathrm{Sym}^5$, or $\mathfrak{sl}_2\times \mathfrak{sl}_2$ acting by $\mathrm{std}\otimes \mathrm{Sym}^2$; see the proof of Theorem~\ref{thm:simple_leq 3}.  These extra possibilities are ruled out by the use of $p$-adic monodromy. 
\end{example}

\begin{example}[Exotic type IV fourfolds]\label{ex:exoticIV}
Let $A$ be an abelian fourfold over number field with quartic CM endomorphism algebra. It is known that the Tate ring of any power of $A$ is generated by divisor classes, so that the Mumford--Tate and the Tate conjecture hold for any power of $A$, cf. \cite[\S 7.5 Type IV(2,1)]{MZ1995}. This fails in characteristic $p$, with the obstruction coming exactly from exotic type IV abelian varieties. For example, if $E$ is a quartic CM field containing an imaginary quadratic subfield in which $p$ splits, $S$ is a special subvariety of $\mathcal{A}_4$ parametrizing abelian fourfolds with $E$-action, and $\mathcal{C}$ is a central leaf in the $\{0^2,(\frac{1}{2})^4,1^2\}$ Newton stratum of the mod $p$ reduction of $S$, then the generic abelian variety over $\mathcal{C}$ is exotic type IV.

For an abelian fourfold over number field with quartic CM endomorphism algebra, or a 
non-exotic type IV abelian fourfold over a characteristic $p$ function field with quartic CM endomorphism algebra, the 
Lie algebra of the $l$-adic monodromy group has a rank 3 toric center, with semi-simple part geometrically isomorphic to $\mathfrak{sl}_2\times \mathfrak{sl}_2$.  For an exotic type IV fourfold, the 
Lie algebra still has semi-simple part geometrically isomorphic to $\mathfrak{sl}_2\times \mathfrak{sl}_2$, but only has a rank 2 toric center. See \S\ref{subsub:Type IV(2,1)}. In the latter case, the absence of a rank 1 central torus causes $A$ to have exceptional Tate classes in $H^4_{\text{ét}}$, and to fail to satisfy the Mumford--Tate conjecture.
\end{example}
\begin{example}[Exotic special correspondence between mod $p$ special subvarieties]\label{Exoexample} This example extends the discussion in Example~\ref{ex:lift_signatures}. It is well-known that two non-isogenous CM abelian varieties in characteristic 0 can become isogenous modulo $p$. The exotic correspondence between the mod $p$ reductions of $\mathrm{U}(1,3)$ and $\mathrm{U}(2,2)$ special subvarieties is in some sense an upgraded version of this phenomenon. 
\par
Suppose that $E$ is an imaginary quadratic field where $p$ does not split. Take a central leaf $\mathcal{C}$ in the $\{0^2,(\frac{1}{2})^4,1^2\}$ Newton polygon of a $\mathrm{U}(1,3)$ special subvariety with $E$-multiplication. Then $\mathcal{C}$ is simultaneously a central leaf in the $\{0^2,(\frac{1}{2})^4,1^2\}$ Newton polygon of a $\mathrm{U}(2,2)$ special subvariety with $E$-multiplication. Therefore, while $\mathrm{U}(1,3)$ and $\mathrm{U}(2,2)$ special subvarieties are very different in characteristic 0, their reductions can share  central leaves. This is related to (but slightly different from) the notion of \textbf{exotic Hecke correspondence} of Xiao and Zhu; cf. \cite[Remark 1.3.5]{xiao_zhu}. \end{example}
\subsection*{Organization of the paper} In Section 2 we recall basic notions related to Tate classes, and discuss the moduli space of abelian varieties, and foliations on the Newton stratification of its mod $p$ fiber. In Sections 3 and 4, we recall basic notions of $p$-adic monodromy. In Section 5 we classify endomorphisms of simple abelian fourfolds based on their Newton polygons. In Section 6 we prove some preliminary results on the mod $p$ Mumford--Tate conjecture. Section 7 is devoted to mixed characteristic lifting, which is used in Section 8 to prove the Tate conjecture for simple non-CM abelian fourfolds.  Section 9 treats non-simple abelian fourfolds. 
\subsection*{AI usage} In writing up the mixed characteristic lifting theorem in \S\ref{sec:lsaf}, R. J. had multiple rounds of conversation with ChatGPT 6 Astra, for the purposes of gaining familiarity with the lifting problem, computing quaternion algebras and Dieudonné crystals, and getting inspiration. The two lifting lemmas~\ref{lm:QCLL} and \ref{lm:graph_lemma} were assisted by ChatGPT 6 Astra. In particular, the symmetric graph construction in Lemma \ref{lm:graph_lemma} is based on ideas originating in interactions with ChatGPT 6 Astra. ChatGPT was also used for literature search, and ChatGPT and Claude Fable 5.2 were used for proofreading. In particular, Claude identified and helped to fix several errors in a near-final draft. The text of the paper was written by the authors.

\subsection*{Acknowledgments}
The authors thank Ananth Shankar, Yunqing Tang and Ziquan Yang for valuable discussions. R. J. also thanks  the organizers (Si Ying Lee, Keerthi Madapusi, George Pappas and Alex Youcis) of the AIM workshop \textit{New $p$-adic perspectives on canonical integral models for Shimura varieties}, where he learned of the notion of exotic Hecke correspondence from Xinwen Zhu. The shadow of this notion appears in Example~\ref{Exoexample}, and helped him to conceptualize the behavior of the otherwise confusing $(1,3)$$\leftrightarrow$$(2,2)$ signature shift. Part of this work was done while M. B. was visiting ICERM in fall 2026, and he thanks the institute for its support.

\section{Preliminaries}
In this section we recall the notions of Tate and Weil classes, review previously known cases of the Tate conjecture for abelian fourfolds, and discuss Newton strata and foliations in the moduli space of abelian varieties in characteristic $p$.
\subsection{The Tate ring and exceptional classes} Let $L$ be a finitely generated field, and let $A$ be an abelian variety over $L$ of dimension $g$. Let $G_L$ be the absolute Galois group of $L$, and let $l\neq \mathrm{char}\, L$. For $r \ge 0$, denote the space of \textbf{Tate classes} of degree $r$ by
    \begin{gather*}
        \mathscr{T}^r_l(A) := H^{2r}_{\et}(A_{\overline{L}}, \bQ_l(r))^{G_L}.
    \end{gather*}
    The Tate conjecture predicts that $\mathscr{T}^r_l(A) $ is the $\bQ_l$-span of algebraic cycles in $\CH^r(A)_\bQ$.

\begin{defn}\label{def:TateRing}
    For $r \ge 0$, let 
    \begin{gather*}
        \mathscr{D}^r(A) \subseteq \CH^r(A)_\bQ
    \end{gather*}
    denote the $\bQ$-subspace spanned by products of divisors. Elements of $\mathscr{D}^r(A)$ are called \textbf{Lefschetz classes}. We let
    \begin{gather*}
        \mathscr{D}^r_l(A) \subseteq \mathscr{T}^r_l(A)
    \end{gather*}
    denote the $\bQ_l$-linear span of $\mathscr{D}^r(A)$ under the cycle class map. We also refer to elements of $\mathscr{D}^r_l(A)$ as \textbf{Lefschetz classes} when no confusion may arise. We let
    \begin{align*}
        \mathscr{T}^\bullet_l(A) &:= \bigoplus\limits_{r \ge 0} \mathscr{T}^r_l(A),\\
        \mathscr{D}^\bullet(A) &:= \bigoplus\limits_{r \ge 0} \mathscr{D}^r(A),\\
        \mathscr{D}^\bullet_l(A) &:= \bigoplus\limits_{r \ge 0} \mathscr{D}^r_l(A).
    \end{align*}
    Then $\mathscr{T}^\bullet_l(A)$ and $\mathscr{D}^\bullet_l(A)$ are graded $\bQ_l$-algebras under cup product, and $\mathscr{D}^\bullet(A)$ is a graded $\bQ$-algebra under intersection product. We refer to $\mathscr{T}^\bullet_l(A)$ as the \textbf{Tate ring} of $A$. 
    \par
    We say that a Tate class $x \in \mathscr{T}^r_l(A)$ is \textbf{exceptional} if $x \notin \mathscr{D}^r_l(A)$.
\end{defn}
\begin{remark}\label{rmk:TateIsogInv}
Let $T(A,L,l)$ be the statement that ``the $l$-adic Tate conjecture holds for the abelian variety $A$ over $L$''. It is known that the truth value of $T(A,L,l)$ remains unchanged if we replace $A$ by an isogeny, $L$ by a finite extension, and $l$ by a different prime . So in practice, we can often assume that $A$ is principally polarized and $L$ is sufficiently large.
\end{remark}
\subsubsection{} The known Tate conjecture for divisors on $A$ implies that to prove the full Tate conjecture for $A$, it suffices to show that all exceptional Tate classes on $A$ are algebraic. In proving the Tate conjecture for abelian fourfolds, a central role is played by Weil classes, as defined below.
\begin{defn}[Weil classes]\label{def:WeilClasses}
    Let $H$ be a Weil cohomology theory on smooth projective varieties over $L$, with coefficients in a characteristic zero field $M$. For example, if we fix an embedding $L \subseteq \bC$, we may take $H$ to be Betti cohomology with coefficients in $\bQ$; otherwise we may take $H$ to be $l$-adic cohomology.
    \par
    Let $k$ be a CM number field, let $d = [k:\bQ]$, and fix an embedding $\iota: k \subseteq \End^0(A)$. Let $N = M \otimes_\bQ k$. Then  $H^1(A)$ has the natural structure of a free module of rank $2g / d$ over $N$. We let $\mathscr{W}_{k, H}(A)$ denote the top exterior power
    \begin{gather*}
        \mathscr{W}_{k, H}(A) := \bigwedge\limits_N^{2g/d} H^1(A),
    \end{gather*}
    and refer to it as the space of \textbf{Weil classes} associated with $(A, \iota, H)$. We naturally regard $\mathscr{W}_{k, H}$ as an $M$-subspace of 
    \begin{gather*}
        \bigwedge\limits_M^{2g/d} H^1(A) \cong H^{2g/d}(A).
    \end{gather*}
    When $H$ is Betti cohomology with coefficients in $\bQ$, we let $\mathscr{W}_{k, B}(A) := \mathscr{W}_{k, H}(A)$; when $H$ is $l$-adic cohomology, we let $\mathscr{W}_{k, l}(A) := \mathscr{W}_{k, H}(A)$. When the choice of $H$ is clear from context, we may simply write $\mathscr{W}_k(A) := \mathscr{W}_{k, H}(A)$.
    \par
    Assume that $2g / d$ is even. When $H$ is Betti cohomology and $\mathscr{W}_{k, B}(A)$ consists of Hodge classes, we refer to $\mathscr{W}_{k, B}(A)$ as the space of \textbf{Hodge--Weil} classes associated with $(A, \iota)$.
\end{defn}
\subsection{Convention on simple and CM abelian varieties} Let $A$ be a dimension $g$ abelian variety over an arbitrary field $K$. For our purposes, the distinction between $A$ and the base change of $A$ to a finite extension of $K$ is largely unimportant (Remark~\ref{rmk:TateIsogInv}). Thus, when we say that $A$ is \textbf{simple}, we mean that $A_{\overline{K}}$ is simple (whereas in the literature this is usually termed as \textbf{geometrically simple}). Similarly, a \textbf{simple factor} of $A$ is a simple factor of a base change of $A$.

We say that $A$ is CM (for \textbf{complex multiplication}) if $\End^0(A_{\overline{K}})$ contains a commutative $\bQ$-algebra of dimension $2g$ (whereas in the literature this is usually called \textbf{potential CM}). A CM abelian variety over $\overline{K}$ can always be descended to a finitely generated subfield. In characteristic $p$, a CM abelian variety over $\overline{K}$ can always be descended \textit{up to isogeny} to a finite field. Here it is essential to allow for isogenies: there exist non-isotrivial 1-dimensional families of CM abelian varieties over finite fields (e.g. a non-isotrivial family of supersingular abelian surfaces over $\mathbb{P}^1_{\mathbb{F}}$ is given by the universal abelian scheme over a component of the supersingular locus of the mod $p$ Siegel space $\mathscr{A}_{2, \mathbb{F}}$). 
\subsection{The Tate conjecture for abelian fourfolds}
The Tate conjecture for dimension and codimension one cycles was proved for abelian varieties over all finitely generated fields in \cite{Tate1966}, \cite{Zarhin1976}, \cite{Faltings1984}, and \cite{MB1985}. Markman's work \cite{Markman2025} establishing the Hodge conjecture for abelian fourfolds opens an approach to the remaining part of the Tate conjecture for abelian fourfolds. 
\subsubsection{Characteristic $0$}
The full Tate conjecture for abelian fourfolds in characteristic 0 follows from Markman's results and earlier works of Moonen--Zarhin \cite{MZ1995,MZ1999}. The latter authors proved the Mumford--Tate conjecture for abelian fourfolds, except for the case where  $\End(A_{\overline{L}})=\bZ$. When $\End(A_{\overline{L}})\neq\bZ$, one can use the well-known implication $$\text{Hodge} + \text{Mumford--Tate} \Rightarrow \text{Tate}.$$
When $\End(A_{\overline{L}})=\bZ$, one can check the Tate conjecture by hand; \cite[\S 4.1]{MZ1995}. This establishes the Tate conjecture for abelian fourfolds in characteristic 0. 
\subsubsection{Finite fields}
The full Tate conjecture for abelian fourfolds over finite fields was first proved in \cite{Broe2026}. This also proves the Tate conjecture for all CM abelian fourfolds in characteristic $p$. 

We briefly review the method of \cite{Broe2026}. For this, suppose that $A$ is an abelian fourfold over a finite field $\mathbb{F}_q$. The first step of the proof shows that (possibly after passing to a finite extension of $\bF_q$) the Tate classes $\mathscr{T}^2_l(A)$ are spanned by Lefschetz classes and Weil classes \cite[Corollary 5.5]{Broe2026}. In other words, we have
\begin{gather*}
    \mathscr{T}^2_l(A) = \mathscr{D}^2_l(A) + \sum\limits_{\iota: k \hookrightarrow \End^0(A)} \mathscr{W}^2_{k, l}(A)(2).
\end{gather*}
Here the above sum runs over embeddings of imaginary quadratic fields $k$ into $\End^0(A)$, such that all nonzero elements of the associated space of Weil classes $\mathscr{W}_{k, l}(A)(2)$ are exceptional Tate classes. This step of the proof adapts $p$-adic Hodge theoretic methods of Ancona \cite{Ancona2021}. 
\par
Fixing an embedding $\iota: k \hookrightarrow \End^0(A)$ with the above properties, the proof now reduces to showing that $\mathscr{W}^2_{k, l}(A)(2)$ consists of algebraic classes. For this, we use Honda--Tate theory to lift $A$ (up to isogeny) to a CM abelian fourfold $\mathcal{A}$ in characteristic zero, equipped with an embedding $k \hookrightarrow \End^0(\mathcal{A})$ which induces $\iota$ upon specialization. The argument then breaks into cases depending on the Hodge type of $\mathscr{W}_{k, B}(\mathcal{A}) \otimes_\bQ \bC$. This Hodge type is of the form $(a, 4-a) + (4-a, a)$ for some $a \in \{2, 3, 4\}$.
\par
If $a = 2$, then $\mathscr{W}_{k, B}(\mathcal{A})(2)$ consists of Hodge classes, which are algebraic by \cite[Corollary 1.6.1]{Markman2025}. By specialization, $\mathscr{W}_{k, l}(A)(2)$ also consists of algebraic classes. 
\par
If $a = 4$, then $p$-adic Hodge theoretic arguments of Ancona show that $A$ is supersingular \cite[Remark A.8(2)]{Ancona2021}, in which case it is well-known that the Tate ring of $A$ is generated by divisor classes. This contradicts our assumption that $A$ possesses exceptional Tate classes.
\par
Now assume that $a = 3$. This case is more complicated: Ancona constructed examples of pairs $(A, \iota)$ such that $\mathscr{W}_{k, l}(A)(2)$ consists of exceptional Tate classes, and $\mathscr{W}_{k, B}^2(\mathcal{A}) \otimes_\bQ \bC$ has Hodge type $(3, 1) + (1, 3)$ for \textit{every} choice of CM lift $\mathcal{A}$ \cite[Proposition A.1]{Ancona2021}. In particular, these examples fail the Mumford--Tate conjecture. The proof circumvents this failure by using results of Markman to construct algebraic cycles on an auxiliary abelian \textit{sixfold} in characteristic zero, specializing them to positive characteristic, and transferring them to $A$ via an algebraic correspondence. This method has further applications to the Tate conjecture over infinite finitely generated fields of positive characteristic: see Lemma \ref{lemWeilAlg13} for details.
\subsubsection{Ordinary abelian fourfolds in characteristc $p$}
The full Tate conjecture for ordinary abelian fourfolds over characteristic $p$ function fields is a consequence of the characteristic $p$ Mumford--Tate machinery developed in \cite{Jiang2023,Jiang2025}. Roughly speaking, there is a well-defined notion of the Mumford--Tate group for ordinary abelian varieties that behaves similarly to its characteristic $0$ counterpart. One can also formulate a version of the  Mumford--Tate group for ordinary abelian varieties; see Conjecture~\ref{conj:OMT}.
In \cite{Jiang2025} it is proved  that, similar to the characteristic $0$ situation, there is again an implication
$$\text{Hodge} + \text{ordinary Mumford--Tate} \Rightarrow \text{ordinary Tate}.$$
Here ``ordinary Mumford--Tate'' and ``ordinary Tate'' respectively refer to the Mumford--Tate and Tate conjectures restricted to ordinary abelian varieties. By adapting classical methods, \cite{Jiang2025} proves Conjecture~\ref{conj:OMT} for all ordinary abelian fourfolds which are not of cohomological Mumford type. The Tate conjecture then follows for these fourfolds. When $A$ is of cohomological Mumford type, the Tate conjecture can be checked by hand, similarly to \cite[\S 4.1]{MZ1995}. 

For reader's convenience, we briefly explain the conceptual reason why we can adapt classical methods to the ordinary characteristic $p$ setting: in ordinary case, the Hodge cocharacter of the canonical lift provides a cocharacter of $G_p(A)$ with weights $\{0,1\}$, hence the representation forms a weak Mumford--Tate pair in the sense of \cite{Pink1998} (and is in particular minuscule). So the classical group theoretic classifications apply.
\subsubsection{The remaining cases} It only remains to prove the Tate conjecture all non-CM, non-ordinary abelian fourfolds $A$ in characteristic $p$.

\subsection{Siegel modular varieties}\label{sub:Siegel moduli}
Let $H$ be a self-dual symplectic $\bZ_{(p)}$-lattice. One can form a Siegel Shimura datum $(\GSp(H_\bQ),\mathfrak{H}^{\pm})$ with reflex field $\bQ$. Let $\bK=\bK_p{\bK}^p$ be a hyperspecial level structure such that $\bK_p=\GSp(H)(\bZ_p)$, ${\bK}^p\subseteq \GSp(H_{\bQ})(\bA_f^p)$ compact open, and $\bK$ leaves $H_{\widehat{\bZ}}$ stable. For $\bK'$ sufficiently small, we get a smooth variety $\Sh_{\bK}(\GSp(H_\bQ),\mathfrak{H}^{\pm})$ over $\bQ$, called the Siegel modular variety, which admits a smooth integral canonical model $\mathscr{A}_{g,\bK}$ over $\bZ_{(p)}$ (\cite[Theorem 2.10]{Mil92}). When $\bK$ is sufficiently small, $H$ gives rise to a universal principal polarized abelian scheme over $\cA_{g,\bK}$ (canonical up to prime-to-$p$ isogeny), whose first $\bZ_{(p)}$-Betti cohomology is the local system induced by $H$. The abelian scheme can be extended to $\mathscr{A}_{g,\bK}$, and we denote by $A^{\mathrm{univ}}\rightarrow\mathscr{A}_{g,\bK}$ the universal abelian scheme over  $\mathscr{A}_{g,\bK}$. Let $\bH_\text{B},\bH_{\text{dR}}, \bH_{l,\text{ét}}$ be the rational Betti, de Rham, $l$-adic étale ($l\neq p$) relative first cohomology of $A^{\mathrm{univ}}\rightarrow \mathscr{A}_{g,\bK}$, and let $\bH_{\text{cris}}$ be the first rational crystalline cohomology of $A^{\mathrm{univ}}_{\mathbb{F}_p}\rightarrow \mathscr{A}_{g,\bK,\Fpbar_p}$. 

We often suppress the level structure and write the moduli space simply as $\mathscr{A}_{g}$.   
Now suppose that we have an abelian variety $A$ over a finitely generated field $L$ and we want to prove the Tate conjecture for $A$. 
By Remark~\ref{rmk:TateIsogInv}, we have the freedom to base change $A$ to a finite extension, and also assume that $A$ is principal polarized. Then $A$ gives rise to a map
$\Spec L\rightarrow \mathscr{A}_{g}$
under which $A^{\mathrm{univ}}$ pulls back to $A$. 
\subsubsection{Spreading out}\label{subsub:spreadout} Let's keep the assumption that $A$ is principally polarized. We will often need to spread out $A$ in order to apply geometric treatment. For example, the crystalline monodromy method in \S\ref{sec:mono p as} relies on spreading out. Here we briefly explain what spreading out means: (up to a finite extension) there is a smooth geometrically connected variety $X_0$ over the base prime field whose fraction field is $L$, such that there is a principally polarized abelian scheme $A_{X_0}$ over $X_0$ whose generic fiber is $A$. Spreading out is always possible. 

We briefly settle down the notation for working with abelian schemes in characteristic $p$.  
Let $X_0$ denote a geometrically connected smooth variety over $\mathbb{F}_q$, let $X=(X_{0})_{\Fpbar}$, and let $x$ be a fixed $\Fpbar$-point of $X_0$. We fix a map
\begin{equation}\label{eq:setupmapsfinitemodel}
    f_0:X_0\rightarrow \mathscr{A}_{g,\mathbb{F}_q},
\end{equation}
and denote by $A_{X_0}$ the abelian scheme over $X_0$ pulled back via the map $f_0$. We will assume that the image of $f_0$ lies in a single Newton stratum. We let
\begin{equation}\label{eq:setupmaps}
    f:X \rightarrow \mathscr{A}_{g,\mathbb{F}_q},
\end{equation}
be the base change of $f_0$ to the algebraic closure $\Fpbar:=\overline{\Fpbar}_p$.

\subsection{Newton strata on $\mathscr{A}_{g,\Fpbar_p}$ and leaves} The mod $p$ fiber of the Siegel modular variety admits a \textbf{Newton stratification} indexed by \textbf{Newton polygons}. We will usually use $\xi$ to denote a Newton polygon and use $W_\xi$ to denote the corresponding Newton stratum. Oort proved that each Newton stratum is foliated by \textbf{central leaves} and \textbf{isogeny leaves}, and is an almost product of a central leaf with an isogeny leaf; cf. \cite{Oort2004}.  Let's briefly review what these objects are: Let $x\in W_\xi$ be an $\Fpbar$-point, corresponding to an abelian variety $A_x$. Then the central leaf $\mathcal{C}_{\xi,x}$ passing through $x$ consists of all points in the Newton stratum whose $p$-divisible group is geometrically isomorphic to $A_x[p^\infty]$. It can be shown that $\mathcal{C}_{\xi,x}$ can be endowed with the structure of a reduced subscheme. On the other hand, the isogeny leaf $\mathcal{I}_{\xi,x}$ passing through $x$ is a maximal reduced subscheme parametrizing ``families of $p$-power isogenies of $A_x$''. It can be shown that $\mathcal{I}_{\xi,x}$ is proper. 
We will often omit $x$ from the subscript and call a central leaf (resp. isogeny leaf) simply $\mathcal{C}_\xi$ (resp. $\mathcal{I}_\xi$). The most well-known example is the supersingular Newton stratum, in which central leaves are 0 dimensional, and isogeny leaves are union of irreducible components of the Newton stratum.

The following is a complete chart (taken from \cite[(8.1)]{Oort2004}) recording Newton strata in the moduli space of abelian fourfolds $\mathscr{A}_{4, \bF_p}$, and the dimension of their central and isogeny leaves. 
\begin{table}[H]
\centering
\renewcommand{\arraystretch}{1.35}
\begin{tabular}{cccc}
\toprule
$\xi$ & 
$\dim W_\xi$ &  $\dim\mathcal C_\xi$ & $\dim\mathcal I_\xi$\\
\midrule
 $\{0^4,1^4\}$ 
   & $10$  & $10$ & $0$\\
 $\{0^3,(\tfrac12)^2,1^3\}$
  & $9$ & $9$ & $0$\\
  $\{0^2,(\tfrac12)^4,1^2\}$
 & $8$ & $7$ & $1$\\
 $\{0,(\tfrac13)^3,(\tfrac23)^3,1\}$
 & $7$  & $6$ & $1$\\
 $\{0,(\tfrac12)^6,1\}$
& $6$ & $4$ & $2$\\
 $\{(\tfrac14)^4,(\tfrac34)^4\}$
 & $6$  & $5$ & $1$\\
$\{(\tfrac13)^3,(\tfrac12)^2,(\tfrac23)^3\}$
 & $5$ & $3$ & $2$\\
 $\{(\tfrac12)^8\}$ 
 & $4$  & $0$ & $4$\\
 \hline
\end{tabular}
\end{table}

\section{Monodromy of mod $p$ abelian schemes}\label{sec:mono p as} 
Let $X_0$ be a smooth and geometrically connected variety over $\mathbb{F}_q$ and let $x\in X_0(\Fpbar)$. Let $W_0=W(\Fpbar_q)$, $W=W(\Fpbar)$ and $K=W[\frac{1}{p}]$, and let $\Breve{\bQ}_p$ denote the completion of the maximal unramified extension of $\bQ_p$. Let $F$ be the absolute Frobenius on $X_0$ and let $\sigma$ be the absolute Frobenius of $\Breve{\bQ}_p$. For a number field $L$, we let $\mathrm{fpl}(L)$ denote the set of finite places of $L$.

\subsection{Monodromy of $F$-isocrystals}\label{subsub:mciiso} The category of convergent $F$-isocrystals over $X_0$, denoted by $\FIsoc(X_0)$, is Tannakian, and we can attach a monodromy group to it, following \cite{Crew1992}. 
\par
Consider an object $\mathcal{M}\in\FIsoc(X_0)$. We denote by $\la\mathcal{M}\ra^{\otimes}$ the tensor abelian subcategory generated by $\mathcal{M}$, and let $\la\mathcal{M}\ra^{\otimes}_K$ be its scalar extension to a field extension $K \supseteq\bQ_p$. We define a fiber functor: 
\begin{align*}
  \omega_x: \langle \mathcal{M}\rangle^\otimes_{\Breve{\mathbb{Q}}_{p}} &\rightarrow \Vect_{\Breve{\mathbb{Q}}_{p}}\\
(\mathcal{N},F)&\rightarrow \cN_x.
\end{align*}

This realizes $\langle \mathcal{M}\rangle^\otimes_{\Breve{\mathbb{Q}}_{p}}$ as a $\Breve{\mathbb{Q}}_{p}$-linear neutral Tannakian category. The monodromy group of $\mathcal{M}$ is defined to be the Tannakian fundamental group $\text{Aut}^\otimes(\omega_x)\subseteq \GL(\omega_x)$. We denote it by $G(\cM,x)$. 
\subsubsection{Frobenius structure}\label{subsub:Frobstructure} The representation category of $G(\cM,x)$ is isomorphic to the scalar extension $\langle \mathcal{M}\rangle^\otimes_{\Breve{\mathbb{Q}}_{p}}$. This is good enough, but in many applications we need the more refined information of $\langle \mathcal{M}\rangle^\otimes$. This is treated in \cite[\S 2.2]{Crew1992}, where the author imposes a \textbf{Frobenius structure} $\Phi$ on $G(\cM,x)$, and shows that $\langle \mathcal{M}\rangle^\otimes$ is equivalent to the representation category of $(G(\cM,x),\Phi)$.

We briefly recall what a Frobenius structure is. Let $G$ be an affine algebraic group over $\Breve{\mathbb{Q}}_{p}$. Then a Frobenius structure on $G$ is a group isomorphism $$\Phi:G^{(\sigma)}\xrightarrow{\sim} G$$  
where $\sigma\in\Gal(\Breve{\mathbb{Q}}_{p}/\bQ_p)$ is the absolute Frobenius. A morphism of algebraic groups with Frobenius structure is a morphism of algebraic groups intertwining the Frobenius structure. 

Let $V$ be a $\Breve{\mathbb{Q}}_{p}$-vector space with a $\sigma$-isocrystal structure $\Psi:V^{(\sigma)}\rightarrow V$. This gives rise to a Frobenius structure $\mathrm{Ad}(\Psi):\GL(V)^{(\sigma)}\rightarrow \GL(V)$ defined by  $\mathrm{Ad}(\Psi)(g): g\rightarrow \Psi g \Psi^{-1}$. A \textbf{representation} of $(G,\Phi)$ on $(V,\Psi)$ is a morphism of algebraic groups with Frobenius structure $(G,\Phi)\rightarrow  (\GL(V),\mathrm{Ad}(\Psi))$. Representations of $(G,\Phi)$ form a rigid $\bQ_p$-linear tensor abelian category, which we denote by $F\text{Rep}_{\Breve{\mathbb{Q}}_{p}}(G,\Phi)$.  

We now return to the discussion of $\langle \mathcal{M}\rangle^\otimes$. The fiber space $\omega_x(\mathcal{M})$ carries a $\sigma$-isocrystal structure $\Psi_x$ that arises from $\mathcal{M}_x$. Then $G(\cM,x)$ is equipped with a Frobenius structure $\Phi_x:=\mathrm{Ad}(\Psi_x)$. There is a tautological representation of $(G(\cM,x),\Phi_x)$ on $(\omega_x(\mathcal{M}),\Psi_x)$, denoted by $\mathrm{std}_x$. Let $\la \mathrm{std}_x\ra^{\otimes}\subseteq F\text{Rep}_{\Breve{\mathbb{Q}}_{p}}(G(\cM,x),\Phi_x)$ be the tensor abelian subcategory generated by $\mathrm{std}_x$. Then 
\begin{theorem}[Crew]
    There is an equivalence of tensor abelian categories $\langle \mathcal{M}\rangle^\otimes\simeq \la \mathrm{std}_x\ra^{\otimes}$.   
\end{theorem}
\begin{proof}
    This follows from \cite[Proposition 2.4]{Crew1992}.
\end{proof}
\subsubsection{Monodromy of overconvergent $F$-isocrystals}\label{subsub:moiso} Let $\FIsoc^\dagger(X_0)$ be the $\bQ_p$-linear tensor abelian category of overconvergent $F$-isocrystals over $X_0$. Let $\mathrm{Fgt}:\FIsoc^\dagger(X_0)\rightarrow \FIsoc(X_0)$ be the forgetful map, and fix an object $\mathcal{M}^\dagger$ in $\FIsoc^\dagger(X_0)$. The monodromy group $G(\mathcal{M}^\dagger, x)$ of $\mathcal{M}^\dagger$ is defined as the Tannakian fundamental group of $\la\mathcal{M}^\dagger\ra^{\otimes}_{\Breve{\bQ}_p}$ with fiber functor $\omega_x \comp \mathrm{Fgt}$. 

As in the convergent case, $G(\cM^\dagger,x)$ is equipped with a Frobenius structure $\Phi_x= \mathrm{Ad}(\Psi_x)$, and there is a standard representation $(G(\cM^\dagger,x),\Phi_x)\acts (\omega_x \comp \mathrm{Fgt}(\mathcal{M}^\dagger),\Psi_x)$. Crew has shown in \cite[Proposition 2.4]{Crew1992} that $\la\mathcal{M}^\dagger\ra^{\otimes}$ is equivalent to the tensor abelian category generated by the standard representation of $(G(\cM^\dagger,x),\Phi_x)$.

\subsection{$l$-adic and $p$-adic monodromy of mod $p$ abelian schemes}\label{sub:becm}

Consider a map as in  (\ref{eq:setupmapsfinitemodel}).
\begin{defn}[Étale and overconvergent crystalline monodromy]\label{def:mon}  For a prime $l \neq p$, write $\bH_{l,f_0}$ or $\bH_{l,X_0}$ for the pullback étale sheaf on $X_0$ via the map $f_0$. We also write $\bH_{\cris,f_0}$ or $\bH_{\cris,X_0}$ for the pullback overconvergent $F$-isocrystal on $X_0$ via the map $f_0$. To ease notation, we will often write them as $\bH_{p,X_0}$ or $\bH_{p,f_0}$. Let $u\in \mathrm{fpl}(\bQ)$, an element $\bH_{u,X_0}$ is called a \textbf{$u$-adic coefficient object}. 

The collection $\{\bH_{l,X_0}\}_{l\neq p}$ is a \textbf{$\bQ$-compatible system of coefficient objects} in the sense of \cite{D2020b}. The collection $\{\bH_{u,X_0}\}_{u\in \mathrm{fpl}(\bQ)}$ is a \textbf{weakly $\bQ$-compatible system of coefficient objects} in the sense of \cite[\S 3.4]{Jiang2023}. It can be upgraded to a $\bQ$-compatible system if we replace the overconvergent $F$-isocrystal $\bH_{u,X_0}$ by an overconvergent $F^a$-isocystal for a suitable $a\in \bN$. 

Let $(G(\bH_{u,X_0},x),\rho_u)$ be the monodromy group of $\bH_{u,X_0}$, together with its tautological representation. When $u\neq p$, this is just the étale monodromy group. When $u=p$, this is the monodromy group of the overconvergent $F$-isocrystal $\bH_{p,X_0}$ with respect to the fiber functor $\omega_x\circ \mathrm{Fgt}$; cf. \S\ref{subsub:moiso}. It is known that $G(\bH_{u,X_0},x)$ is a reductive $\bQ_u$-group. When $u=l\neq p$ (resp. $u=p$), we will call it the \textbf{$l$-adic étale monodromy} (resp. \textbf{overconvergent crystalline monodromy}) of $A_{X_0}$. For convenience, we will usually just call it the \textbf{$u$-adic monodromy}. To ease notation, we will often denote $G(\bH_{u,X_0},x)$ simply by $G_u(f_0)$. It is known that the neutral component $G_u(f_0)^\circ$ together with $\rho_u$ is 
``independent of $u$'' in the sense that it admits a rational model over a sufficiently large number field in the sense of \cite[Theorem 1.2.4]{D2020a}. We let $H_u(A) = (G_u(A) \cap \SL(V_u(A)))^\circ$. 
\end{defn}
\begin{notation}[Independence from the choice of model]\label{not:hidingfinitefield} In practice, we often work with coefficient objects and their monodromy without referring to a specific finite field model. It is convenient to set up a model-free system of notation. Let $f_0$ and $f$ be as in \ref{eq:setupmapsfinitemodel}. For each $u$, the neutral component $G(\bH_{u,X_0},x)^\circ$ depends only on $f$, and not on $f_0$. There exist models $f_0$ of $f$ defined over certain sufficiently large finite fields $\bF_q$ such that $G(\bH_{u,X_0},x)$ is the smallest among all possible $G(\bH_{u,X_0},x)$ over all 
possible finite field models. We will write $\bH_{u,X}$ or $\bH_{u,f}$ for the coefficient objects over those finite field models, and write $G(\bH_{u,f},x)$ for the monodromy group. To ease notation, we will often denote $G(\bH_{u,f},x)$ simply by $G_u(f)$. Similarly, we will write $\bH_{p,X}^-$ for the underlying convergent $F$-isocrystal of $\bH_{p,X}$.%
\end{notation}

The following is a mild structure result of the monodromy groups, which is an analogue of Bogomolov's result in characteristic 0:
\begin{lemma}\label{lemHomotheties}
    For $u\in \mathrm{fpl}(\bQ)$, $G(\bH_{u,X_0},x)$ contains the rank 1 diagonal torus of homothety. %
\end{lemma}
\begin{proof}
It suffices to prove it for $u\neq p$. Let $x_0$ be a finite field model of $x$ so that $G(\bH_{u,x_0},x)$ is connected. Then $G(\bH_{u,x_0},x)\hookrightarrow G(\bH_{u,X_0},x)$. It suffices to show that $G(\bH_{u,x_0},x)$ contains the torus of homothety. Let $\alpha_1,\dots,\alpha_{2g}$ be the eigenvalues of the geometric Frobenius. Picking a suitable basis on $\omega_x(\bH_{u,x_0})\otimes \overline{\bQ}_u$, we have $\diag(\alpha_1,\dots,\alpha_{2g})\in G(\bH_{u,x_0},x)(\overline{\bQ}_u)$. Let $\diag(c,c,.\dots, c)$ be a $\overline{\bQ}_u$-point of the homothety torus. Then it suffices to show that for any $e_i\in \bZ$ with $\prod_{i=1}^{2g} \alpha_i^{e_i}=1$, we also have $\prod_{i=1}^{2g} c^{e_i}=1$, for any $c\in \overline{\bQ}_u^*$. This follows from the Weil conjectures: choosing an embedding into $\bC$, we have $\prod_{i=1}^{2g} |\alpha_i|^{e_i}=1$, but $|\alpha_i|=\sqrt{q}$ for all $i$, so we get $\sum_{i=1}^{2g}e_i=0$. This implies the desired result.
\end{proof}
\begin{remark}\label{remBogomolov}
There is a subtle difference between the above lemma and Bogomolov's theorem. Let $u\neq p$. If we only take the Lie algebra of the image of the étale fundamental group without taking Zariski closure (as what Bogomolov does in characteristic 0), it may not contain the torus of homothety; cf. \cite{Zarhin2007}. 
\end{remark}
\subsection{Parabolicity and non-vanishing unipotent}
Let $f_0:X_0\rightarrow \mathscr{A}_{g,\Fpbar_q}$ be as in (\ref{eq:setupmapsfinitemodel}) and let $x\in X_0(\Fpbar)$. Recall that $\bH_{p,X_0}$ is the overconvergent Dieudonné isocrystal of $A_{X_0}$, with its underlying convergent isocrystal $\bH_{p,X_0}^-$.  
Let $\nu_x$ be the Newton cocharacter associated to $x$. Since $\bH_{p,X_0}$ has constant Newton polygon, Crew's parabolicity conjecture (now a theorem of D'Addezio; cf. \cite{D2020b}) implies the following
\begin{theorem}[D'Addezio]
    $G(\bH_{p,X_0}^-,x)$ is the parabolic subgroup of $G(\bH_{p,X_0},x)$ fixing $\nu_x$. 
\end{theorem}
Let $\Fil_\bullet \bH_{p,X_0}$ be the slope filtration on $\bH_{p,X_0}$, and let $\gr_\mu\bH_{p,X_0}$ be the slope $\mu$ isoclinic graded object. Then the Levi of $G(\bH_{p,X_0},x)$ fixing $\nu_x$ is the monodromy group of $G(\bigoplus_\mu\gr_\mu\bH_{p,X_0},x)$. It is also known that the unipotent of  $G(\bH_{p,X_0}^-,x)$ fixing $\nu_x$ is the monodromy group of the formal $F$-isocrystal $\bH_{p,X_0}|_{X_\Fpbar^{/x}}$, which is proved by Chai in the ordinary case (cf. \cite{Chai2003}), and by D'Addezio and van Hoften in general (cf. \cite[Theorem 3.4.4]{DAddezio2022}).

In the following, let $U_{\nu_x}$ be the unipotent of $G(\bH_{p,X_0},x)$ associated to $\nu_x$. The unipotent is intimately related to the formal geometry of the completed map $f_{\Fpbar}^{/x}:X_{0,\Fpbar}^{/x}\rightarrow \mathscr{A}_{g,\Fpbar}^{/f_{\Fpbar}(x)}$. In fact, when $f_0$ maps to the ordinary stratum, $U_{\nu_x}$ admits a canonical $\bQ_p$-structure $U_{\nu_x}^{\Phi}$ which can be canonically identified with the rational cocharacter group of the smallest formal subtorus of the Serre--Tate formal torus $\mathscr{A}_{g,\Fpbar}^{/f_{\Fpbar}(x)}$ that contains the image of $f_{\Fpbar}^{/x}$; cf. \cite{Jiang2023}. 

For higher Newton strata, a similar scenario is conjectured to hold: suppose that the image of $f_0$ lies in a central leaf $\cC_\xi$, then $U_{\nu_x}$ can be canonically identified with a certain \textbf{Dieudonné--Lie structure} on the smallest \textbf{strongly Tate-linear subspace} of $\mathscr{A}_{g,\Fpbar}^{/f_{\Fpbar}(x)}$ that contains the image of $f_{\Fpbar}^{/x}$. We won't explain the terminologies here, since they won't be used. We refer the readers to \cite[Conjecture 5.5.7 and \S6]{DAddezio2022} for more details. 

For our application, we only need the following: 
\begin{theorem}[Non-vanishing unipotent]\label{thm:nonvan}
Notation as above. Suppose that $A_{X_0}$ is non-CM, and let $x\in X_0(\Fpbar)$. Then $U_{\nu_x}$ is positive dimensional. 
\end{theorem}
\begin{proof}
We can reduce to the case where $\dim X_0=1$ and (up to a finite cover) $X_0$ projects non-constantly to a central leaf. In fact, for general $X_0$ and $x\in X_0(\mathbb{F})$, we just find a curve $C_0\rightarrow X_0$ whose image passes by $x$ and whose projection to a central leaf has positive dimension, and note that the monodromy group of $C_0$ embeds into that of $X_0$. If the monodromy of $C_0$ has nontrivial unipotnent, so is $X_0$.

Recall that $\bH_{p,X_0}$ admits a slope filtration $\Fil_\bullet \bH_{p,X_0}$. Let $\gr_{\mu}\bH_{p,X_0}$ be the $\mu$-isoclinic graded piece. In general, each step in the slope filtration, as well as each $\gr_{\mu}\bH_{p,X_0}$, lies only in $\FIsoc(X_0)$ but not $\FIsoc^\dagger(X_0)$. In the following, it is more convenient to work in larger categories $F^n\textbf{-Isoc}(X_0)$ and $F^n\textbf{-Isoc}^\dagger(X_0)$ in the sense of \cite{D2020a}, where $n$ is a common denominator of the slopes (alternatively, one can just work with $F^\infty$-isocrystals, which are also defined in \textit{loc.cit.}). We can twist each $\gr_{\mu}\bH_{p,X_0}$ by a constant object in $F^n\textbf{-Isoc}({X}_0)$ to get a unit-root $F^n$-isocrystal $\mathscr{E}_\mu$. By \cite[Theorem 3.7]{Kedlaya2017}, each $\mathscr{E}_\mu$ corresponds to a $\bQ_{p^n}$-coefficient $p$-adic lisse sheave $\mathscr{E}^{\et}_\mu$ over $X_0$. 

Let $\overline{X}_0$ be the smooth projective curve containing $X_0$. We have a map $\overline{f}_0:\overline{X}_0\rightarrow \mathscr{A}_{g,\Fpbar_p}^{\Sigma}$ extending $f_0$, where the later being a smooth toroidal compactification. Suppose for the sake of contradiction that $\nu_x$ admits a trivial unipotent. Then the parabolic fixing $\nu_x$ is $G(\bH_{p,X_0},x)$ itself. By Crew's parabolicity conjecture (\cite{D2020b}), each $\gr_{\mu}\bH_{p,X_0}$ need to be an overconvergent $F$-isocrystal as well, and $\bH_{p,X_0}=\bigoplus_{\mu} \gr_{\mu}\bH_{p,X_0}$ in the category $\FIsoc^\dagger(X_0)$. Tsuzuki's theorem (\cite[Theorem 3.9]{Kedlaya2017}) implies that for all $\mu$, $\mathscr{E}_\mu^{\et}$ is potentially unramified. Possibly replacing $X_0$ by an étale cover, we see that $\bH_{p,X_0}$ extends to an $F^n$-isocrystal $\bH_{p,\overline{X}_0}$ with constant slopes over the compactification $\overline{X}_0$.

By companion conjecture for curves (see \cite[Theorem 1.2]{Ked2022}), for $l\neq p$ the $l$-adic étale lisse sheaf $\bH_{l,X_0}$, after some scalar extension, also extends to $\overline{X}_0$. Néron–Ogg–Shafarevich then says that the abelian scheme $A_{X_0}$ also extends to $\overline{X}_0$, so $\overline{X}_0$ maps entirely into $\mathscr{A}_{g,\Fpbar_p}$. 
Now by full-faithfulness \cite[Theorem 5.3]{Kedlaya2017} the Dieudonné crystal of $A_{\overline{X}_0}$ equals $\bH_{p,\overline{X}_0}$. So $A_{\overline{X}_0}$ has constant slopes  over $\overline{X}_0$, meaning that $\overline{X}_0$ maps completely inside a Newton stratum, hence by our assumption has non-constant projection to a central leaf. Since central leaves are quasi-affine (cf. \cite{Oort1999}) and $\overline{X}_0$ is proper, this is a contradiction.

\end{proof}

\subsection{More on Frobenius structures}\label{subsec:MFS}
According to \S\ref{subsub:Frobstructure}, $G(\bH_{p,X_0},x)$ is equipped with a Frobenius structure $\Phi_x=\mathrm{Ad}(\Psi_x)$ induced from the $F$-isocrystal structure $\Psi_x$ on $\omega_x(\bH_{p,X_0})$. 
\begin{lemma}\label{lm:Frob--Sigma}
  $\Phi_x$ restricts to Frobenius structures on the parabolic, Levi, and unipotent of $G(\bH_{p,X_0},x)$ associated to $\nu_x$. Moreover, $\Phi_x$ endows $\Lie G(\bH_{p,X_0},x)$ with the structure of a $\sigma$-isocrystal. Similar for the Lie algebras of  the parabolic, Levi, and unipotent.
\end{lemma}
\begin{proof}
    Note that $\Phi_x$ restricts to Frobenius structures on parabolic and Levi since they are the monodromy groups of $F$-isocrystals $\bH_{p,X_0}^-$ and $\bigoplus_\mu\gr_\mu\bH_{p,X_0}$, respectively. One can similarly prove the unipotent case by realizing it as the monodromy of  $\bH_{p,X_0}|_{X_\Fpbar^{/x}}$. But we can do a more direct proof: $\Phi_x$ sends $G(\bH_{p,X_0}^-,x)^{(\sigma)}$ to $G(\bH_{p,X_0}^-,x)$, and it must maps the unipotent radical of $G(\bH_{p,X_0}^-,x)^{(\sigma)}$ to the unipotent radical of $G(\bH_{p,X_0}^-,x)$. Therefore $\Phi_x$ restricts to a Frobenius structure on the unipotent radical, as desired.  Finally, the $\sigma$-isocrystal structure on  $\Lie G(\bH_{p,X_0},x)$ follows from the fact that $\Phi_x:\Lie G(\bH_{p,X_0},x)^{(\sigma)}\rightarrow \Lie G(\bH_{p,X_0},x)$ is an isomorphism. 
\end{proof}
\section{Examples of $p$-adic monodromy for abelian fourfolds}
In this section, we will consider $p$-adic monodromy of abelian fourfolds with Newton polygon $\xi= \{(\frac{1}{4})^4,(\frac{3}{4})^4\}$ or $\{0^2,(\frac{1}{2})^4,1^2\}$. We will run concrete computation for the $p$-adic monodromy. The examples serve two purposes: first, they familiarize the reader with the computation; second, they will be used later as technical inputs for the Mumford--Tate conjecture. The notation for this chapter will be consistent with that of the last chapter.

\subsection{Newton polygon $\{(\frac{1}{4})^4,(\frac{3}{4})^4\}$}\label{example:1/4-3/4}
Let $x\in W_{\xi}(\Fpbar)$. We can choose a $\sigma$-invariant basis on $\omega_x(\bH_{p})$ so that the Frobenius structure $\Psi_x$ and the symplectic pairing $\psi_x$ are given by  $$\Psi_x=\begin{pmatrix}
\begin{array}{c|c}
     \begin{array}{cccc}
  &&&p \\
   
  1&&&\\
  &1&& \\
  &&1&  \\
 \end{array}&  \\\hline
     &  \begin{array}{cccc}
  &&&{1} \\
  
  p&&&\\
  &p&& \\
  &&p&  \\
 \end{array}
\end{array}
\end{pmatrix},\qquad\psi_x=\begin{bmatrix}
    &I_4\\
    -I_4&
\end{bmatrix}.$$
The the whole group $\GSp(\omega_x(\bH_{p}),\psi_x)$ is stable under the Frobenius structure $\Phi_x=\ad \Psi_x$ (in fact, $\GSp(\omega_x(\bH_{p}),\psi_x)$ equals the $p$-adic monodromy group $G(\bH_{p,\mathcal{C}_{\xi}},x)^\circ$, where $\mathcal{C}_{\xi}$ is a central leaf). The Newton quasi-cocharacter is given by $\nu_x(t)=\diag[t,t,t,t,t^3,t^3,t^3,t^3]$. Let $U_{\nu_x}^\xi$ be the unipotent of $\GSp(\omega_x(\bH_{p}),\psi_x)$ with respect to $\nu_x$. In the following, we will compute the Frobenius structure on $U_{\nu_x}^\xi$. In fact, we have $$\Lie U_{\nu_x}^\xi=\left\{\begin{pmatrix}
     \begin{array}{c|c}
 \mathbf{0}_4 &  \mathbf{0}_4\\
 \hline
 D & \mathbf{0}_4
 \end{array}
\end{pmatrix}: D=D^{\mathrm{t}}\in \mathrm{M}_4(\Breve{\bQ}_p).\right\}
$$
Interpret $\Phi_x$ as a $\sigma$-linear map. Its action on $\Lie U_{\nu_x}^\xi$ is given by
\begin{equation}\label{eq:matrix of Phi 1}
    D\xrightarrow{\Phi_x}   \begin{pmatrix}
  &&&1 \\
   
  p&&&\\
  &p&& \\
  &&p&  \\
 \end{pmatrix}\sigma D \begin{pmatrix}
  &1&& \\
   
  &&1&\\
  &&&1 \\
 p^{-1} &&&  \\
 \end{pmatrix}.
\end{equation}
Recall from Lemma~\ref{lm:Frob--Sigma} that $\Lie U_{\nu_x}^\xi$ is equipped with a $\sigma$-isocrystal structure arising from $\Phi_x$. To see it clearly, let's write 
$$D=\begin{pmatrix}
   b_{11}&b_{12}&b_{13}&b_{14}\\
 b_{12}&b_{22}&b_{23}&b_{24}\\
 b_{13}&b_{23}&b_{33}&b_{34}\\
b_{14}&b_{24}&b_{34}&b_{44}
\end{pmatrix}, \quad b_{ij}\in \Breve{\bQ}_p.$$
Then $$\Phi_x(D) = \sigma \begin{pmatrix}
   \frac{1}{p}b_{44}&b_{14}&b_{24}&b_{34}\\
 b_{14}&pb_{11}&pb_{12}&pb_{13}\\
 b_{24}&pb_{12}&pb_{22}&pb_{23}\\
b_{34}&pb_{13}&pb_{23}&pb_{33}
\end{pmatrix}.$$
It is easy to see (exercise for the reader) from the above explicit computation that $\Lie U_{\nu_x}$ decomposes into five $\Phi_x$-invariant subspaces of dimension 2, each of which is isomorphic to the unique irreducible $\sigma$-isocrystal of slope $1/2$. Just to name one: the two dimensional subspace corresponding to $(b_{13},b_{24})$, where the Frobenius structure acts as
$(b_{13},b_{24})\xrightarrow{\Phi_x}  \sigma(b_{24},pb_{13})$, is a simple $\sigma$-isocrystal of slope $1/2$.  We remark that the geometric reason behind the fact that  $\Lie U_{\nu_x}^\xi$ decomposes into 5 pieces is that $\dim \cC_\xi= 5$. In fact, it is known from the work of Chai that $\cC_\xi^{/x}$ is a direct sum of five Lubin--Tate formal groups of slope $1/2$; cf. \cite{Chai2005}. 
\begin{lemma}\label{lm:evendim two slope}
    Let $f_0:X_0\rightarrow \mathscr{A}_{4,\Fpbar_q}$ be a morphism whose image lies generically in $\mathcal{C}_\xi$. Let $x\in X_0(\Fpbar)$ mapping to $\mathcal{C}_\xi$. Then the unipotent group $U_{\nu_x} \subseteq G(\bH_{p,X_0},x)$ has even dimension. If $A_{X_0}$ is non-CM, then $\dim U_{\nu_x}$ is furthermore positive. 
\end{lemma}
\begin{proof}
    Note that $U_{\nu_x}\subseteq U_{\nu_x}^\xi$ interwining the Frobenius structure. Now $\Phi_x$ makes $\Lie U_{\nu_x}^\xi$ a $\sigma$-isocrystal of slope $1/2$. Therefore $\Lie U_{\nu_x}$ is also a $\sigma$-isocrystal of slope $1/2$, hence having even dimension. When $A_{X_0}$ is non-CM, the unipotent is nontrivial by Theorem~\ref{thm:nonvan}. 
\end{proof}
\subsubsection{Variants} The following lemma is a variant of Lemma~\ref{lm:evendim two slope} for abelian threefolds which will be used later: 
\begin{lemma}\label{lm:3slope} Let $f_0:X_0\rightarrow \mathscr{A}_{3,\Fpbar_q}$ be a morphism whose image lies generically in $\mathcal{C}_\xi$, where $\xi=\{(\frac{1}{3})^3,(\frac{2}{3})^3\}$.  Let $x\in X_0(\Fpbar)$ mapping to $\mathcal{C}_\xi$. Then the unipotent group $U_{\nu_x} \subseteq G(\bH_{p,X_0},x)$ has  dimension divisible by $3$. If $A_{X_0}$ is non-CM, then $\dim U_{\nu_x}$ is either 3 or 6. 
\end{lemma}
\begin{proof}
   A similar computation shows that $\Lie U_{\nu_x}$ is a $\sigma$-isocrystal of slope $1/3$, hence having dimension divisible by 3. When $A_{X_0}$ is non-CM, the unipotent is nontrivial by Theorem~\ref{thm:nonvan}. 
\end{proof}
\begin{remark}
    Let $g \ge 3$ be odd. The above lemma extends in an obvious way to abelian $g$-folds with Newton polygon $\{(\frac{1}{g})^g,(\frac{g-1}{g})^g\}$: the dimension of $U_{\nu_x}$ is then divisible by $g$, and is positive if
    $A_{X_0}$ is non-CM.
\end{remark}

\subsection{Newton polygon $\{0^2,(\frac{1}{2})^4,1^2\}$}\label{example:0-1/2-1} Let $x\in W_{\xi}(\Fpbar)$. We can choose a $\sigma$-invariant basis on $\omega_x(\bH_{p})$ so that  
$$\Psi_x=\begin{pmatrix}
1&&&&&&&&\\
&1&&&&&&&\\
&&&&1&&&&\\
&&&-p&&&&&\\
&&&&&&1&&\\
&&&&&-p&&&\\
&&&&&&&p&\\
&&&&&&&&p
\end{pmatrix},\qquad\psi_x=\begin{bmatrix}
    &&&I_2\\
    &H&&\\
   &&H& \\ 
-I_2&&&
\end{bmatrix},\qquad H=\begin{bmatrix}
    &1\\
    -1&
\end{bmatrix}.$$
Again, the whole group $\GSp(\omega_x(\bH_{p}),\psi_x)$ is stable under the Frobenius structure $\Phi_x=\ad \Psi_x$ (in fact, $\GSp(\omega_x(\bH_{p}),\psi_x)$ equals the $p$-adic monodromy group $G(\bH_{p,\mathcal{C}_{\xi}},x)^\circ$). The Newton cocharacter is $\nu_x(t)=\diag[1,1,t,t,t,t,t^2,t^2]$. In this case, it is easy to see that 
$$\Lie U_{\nu_x}^\xi=\left\{\begin{pmatrix}
\begin{array}{cc|cccc|cc}
&&&&&&&\\
&&&&&&&\\\hline
a_1&b_1&&&&&&\\
a_2&b_2&&&&&&\\
a_3&b_3&&&&&&\\
a_4&b_4&&&&&&\\\hline
c_1&c_3&a_2&-a_1&a_4&-a_3&&\\
c_3&c_2&b_2&-b_1&b_4&-b_3&\;&\;
\end{array}
\end{pmatrix}: a_i,b_j,c_k\in \Breve{\bQ}_p\right\}.
$$
The action of $\Phi_x$ on $\Lie U_{\nu_x}^\xi$ is
\begin{equation}\label{eq:matrix of Phi 2}
   \begin{pmatrix}
\begin{array}{cc|cccc|cc}
&&&&&&&\\
&&&&&&&\\\hline
a_1&b_1&&&&&&\\
a_2&b_2&&&&&&\\
a_3&b_3&&&&&&\\
a_4&b_4&&&&&&\\\hline
c_1&c_3&a_2&-a_1&a_4&-a_3&&\\
c_3&c_2&b_2&-b_1&b_4&-b_3&\;&\;
\end{array}
\end{pmatrix}\xrightarrow{\Phi_x}    \sigma\begin{pmatrix}
\begin{array}{cc|cccc|cc}
&&&&&&&\\
&&&&&&&\\\hline
a_2&b_2&&&&&&\\
-pa_1&-pb_1&&&&&&\\
a_4&b_4&&&&&&\\
-pa_3&-pb_3&&&&&&\\\hline
pc_1&pc_3&-pa_1&-a_2&-pa_3&-a_4&&\\
pc_3&pc_2&-pb_1&-b_2&-pb_3&-b_4&\;&\;
\end{array}
\end{pmatrix}.
\end{equation}
Note that $\dim \cC_\xi= 7$, and $\Lie U_{\nu_x}^\xi$ decomposes into 7 $\Phi_x$-invariant subspaces, 
3 of dimension 1 and 4 of dimension 2. The dimension 2 ones are isomorphic to the simple $\sigma$-isocrystal of slope $1/2$. The dimension 1 ones are isomorphic to the simple $\sigma$-isocrystal of slope $1$. The deformation space $\cC_\xi^{/x}$ is not simply a formal group, but a Tate linear formal variety in the sense of \cite{Chai2024}.

\begin{lemma}\label{lm:keylemma}
 Let $f_0:X_0\rightarrow \mathscr{A}_{4,\Fpbar_q}$ be a morphism whose image lies generically in $\mathcal{C}_\xi$. Let $x\in X_0(\Fpbar)$ mapping to $\mathcal{C}_\xi$. Suppose that $\End^0(A_{X_0})$ admits a complex imaginary subfield where $p$ splits and the Rosati involution acts as conjugation. If $\End^0(A_{X_0})$ splits $\bH_{p,X_0}$ into two pieces which are of Newton polygons $\{0^2,(\frac{1}{2})^2\}$ and $\{1^2,(\frac{1}{2})^2\}$ respectively, then the unipotent group $U_{\nu_x} \subseteq G(\bH_{p,X_0},x)$ has dimension 0, 2 or 4.  
If $A_{X_0}$ is non-CM, then $\dim U_{\nu_x}$ is furthermore positive. 
\end{lemma}
\begin{proof}
Let $p$ split as $v\overline{v}$ in the complex field. Note that $\bH_{p,X_0}$  splits into two mutually dual $F$-isocrystals: $$\bH_{p,X_0}\simeq \bH_{p,X_0}\otimes_{E_{\bQ_p}} E_{\bQ_p}\simeq(\bH_{p,X_0}\otimes{E_v})\oplus(\bH_{p,X_0}\otimes{E_{\overline{v}}}).$$
We denote the two pieces by  $\bH_{v,X_0}$ and $\bH_{\overline{v},X_0}$, respectively. We can assume that $\bH_{\overline{v},X_0}$ has generic Newton polygon $\{0^2,(\frac{1}{2})^2\}$. 

In this case, $G(\bH_{v,X_0},x)$ is the projection of $G(\bH_{\overline{p},X_0},x)$ to $\GL(\omega_x(\bH_{v,X_0}))$ (similar for $G(\bH_{\overline{v},X_0},x)$), and  $G(\bH_{\overline{p},X_0},x)^{\der}\subseteq G(\bH_{v,X_0},x)^{\der}\times G(\bH_{\overline{v},X_0},x)^{\der}$ is the graph of an isomorphism. 

Let $\nu_v$ be composition of $\nu_x$ with the projection to  $\GL(\omega_x(\bH_{v,X_0}))$. Then $\nu_v$ has weight $\{0^2,1^2\}$. Let $U_{\nu_v}$ be the unipotent of $G(\bH_{v,X_0},x)$ with respect to $\nu_v$. Since projection induces an isomorphism on derived groups, we have $U_{\nu_x}\simeq U_{\nu_v}$.  

Now we do a computation that is similar to the above example but much easier. We can choose a $\sigma$-invariant basis on $\omega_x(\bH_{v,X_0})$ so that the Frobenius is the following matrix
$$\Psi_x=\begin{pmatrix}
1&&&&&\\
&1&&&&\\
&&&&1&\\
&&&-p&&\\
\end{pmatrix}.$$
The cocharacter $\nu_v$ at $x$ is $\diag[1,1,t,t]$. In this case, it is easy to see that 
$$\Lie U_{\nu_v}\subseteq L:=\left\{\begin{pmatrix}
\begin{array}{cc|cc}
&&&\,\,\,\\
&&&\,\,\,\\
\hline
a_1&b_1&&\\
a_2&b_2&&
\end{array}
\end{pmatrix}: a_i,b_j\in \Breve{\bQ}_p\right\}.
$$
The action of $\Phi_x$ on $L$ is
\begin{equation*}
  \begin{pmatrix}
\begin{array}{cc|cc}
&&&\,\,\,\\
&&&\,\,\,\\
\hline
a_1&b_1&&\\
a_2&b_2&&
\end{array}
\end{pmatrix}\xrightarrow{\Phi_x} \sigma\begin{pmatrix}
\begin{array}{cc|cc}
&&&\,\,\,\\
&&&\,\,\,\\
\hline
a_2&b_2&&\\
-pa_1&-pb_1&&
\end{array}
\end{pmatrix}.
\end{equation*}
Note that $L$ decomposes into a direct sum of two simple $\sigma$-isocrystals of slope $1/2$. Since $\Lie U_{\nu_v}$ is a $\Phi_x$-invariant subspace, it is a sub-$\sigma$-isocrystal of $L$, hence must also be of slope $1/2$. Therefore  $\dim\Lie U_{\nu_v}$ is even. The positivity is again a consequence of Theorem~\ref{thm:nonvan}. The lemma follows. 
\end{proof}

\section{Endomorphisms of simple non-CM abelian fourfolds}
\subsection{The Albert classification} Let $A$ be a simple abelian fourfold which is non-CM over a characteristic $p$ algebraically closed field. Then the Albert classification yields the following possibilities for the its endomorphism algebra $D$:

\begin{table}[H]
\centering
\renewcommand{\arraystretch}{1.35}
\begin{tabular}{|c|c|c|}
\toprule
\multicolumn{2}{|c|}{Type}   &$D$ \\
\midrule
 & I(1) & $\bQ$ \\\cline{2-2}
 I & I(2) & Real quadratic field \\\cline{2-2}
  & I(4) & Totally real quartic field \\
  \hline
 II  & II(1) & Indefinite quaternion over $\bQ$ \\\cline{2-2}
  & II(2) & Indefinite quaternion over a real quadratic field \\\hline
 III  & III(1) & Definite quaternion over $\bQ$ \\\cline{2-2}
  & III(2)* & Definite quaternion over a real quadratic field \\\hline
 & IV(1,1) & Imaginary quadratic field \\\cline{2-2}
  IV  & IV(2,1) & Quartic CM field \\\cline{2-2}
  & IV(2)*  & Quaternion over an imaginary quadratic field\\  \hline 
\end{tabular}
\caption{Possible endomorphism algebras}\label{tb:Endo4}
\end{table}

In the table, the types with a ``\,*\,'' does not happen: having type III(2) endomorphism algebra means that the $l$-adic monodromy group $G_l(A)$ is contained in the centralizer of $D$ in $\GSp(H_l^1(A))$, which is a torus.  We will see that IV(2) cannot happen either, see Theorem~\ref{eq:Newton--endo}.

\subsection{Classification with prescribed Newton polygons}\label{sub:final_classification} Let $A$ be a simple abelian variety of dimension $g$ over a characteristic $p$ algebraically closed field $L$ of Newton polygon $\xi$. Let $K=W(L)$. Let $D=\End^0(A)$ and let $F$ be a field lying in the center of $D$ and let $\Sigma_p$ be the set of places of $F$ lying above $p$. We have $F\otimes_{\bQ} \bQ_p\simeq \prod_{v\in \Sigma_p} F_v e_v$, where $e_v$'s are idempotents that satisfy $\sum_{v\in \Sigma_p} e_v=1$ and $e_ve_{w}=\delta_{vw}$.  This induces a splitting  $$D\otimes_{\bQ} \bQ_p\simeq D\otimes_{F} (F\otimes_{\bQ}\bQ_p)= \prod_{v\in \Sigma_p} D_ve_v,\qquad D_v=D\otimes_F F_v.$$
The splitting also induces a splitting of the rational Dieudonné module $$\bH_{p,A}:=\bD(A)[1/p]=\bigoplus_{v\in \Sigma_p} \bH_{v,A},\qquad\bH_{v,A}:=e_v \bH_{p,A}.$$ 
The action of $D\otimes_{\bQ} \bQ_p$ on $\bH_{p,A}$ also splits, yielding  \begin{equation}\label{eq:embdv}\iota_{vw}:D_v\rightarrow\End(\bH_{w,A}).
\end{equation}
The following is clear: 
\begin{lemma}\label{lm:split}
 $\bH_{p,A}$ is a free module over $F\otimes_{\bQ}K$ and  each $\bH_{v,A}$ is a free module over ${F_v}\otimes_{\bQ_p} K$ of the same rank. In particular $\dim_{K}\bH_{v,A}=2g{[F_v:\mathbb{Q}_p]}/{[F:\bQ]}$. 
\end{lemma} 
Now we consider the distribution of slopes among $\mathbb{H}_{v,A}$. 
\begin{lemma}\label{lm:RosatiD}
Let $\dagger$ be the Rosati involution on $D$.  The following are true:
\begin{enumerate}
    \item If $\dagger$ is identity over $F$, then each  $\mathbb{H}_{v,A}$ is self-dual. 
    \item If $F$ is CM and  $\dagger$ is a complex conjugation map $x\rightarrow \overline{x}$, then $\mathbb{H}_{\overline{v},A}$ is dual to $\mathbb{H}_{v,A}$.
\end{enumerate}
\end{lemma}
\begin{proof}
    Write $\psi$ for the polarization. For any $x\in \mathbb{H}_{v,A}$ and $y\in \mathbb{H}_{w,A}$, we have  $$\psi(x, y)=\psi(e_{v}x, y)=\psi(x,e^\dagger_vy).$$
  If $\dagger$ is identity over $F$, then $e^\dagger_vy=e_vy$, which is 0 whenever $w\neq v$. Therefore $\psi$ induces a self-duality over $\mathbb{H}_{v,A}$. If $\dagger$ is conjugation over $F$, then $e^\dagger_vy=e_{\overline{v}}y$, which is 0 whenever $w\neq \overline{v}$. Therefore  $\psi$ induces a duality between $\mathbb{H}_{v,A}$ and $\mathbb{H}_{\overline{v},A}$. 
\end{proof}

In the following we suppose furthermore that $F$ is the center of $D$. 

Let $d=\sqrt{\dim_F D}$ be the degree of $D$. Then $D_v$ is also a central simple algebra over $F_v$ of degree $d$. 
Note that $\End(\bH_{v,A})$ splits as a product of the endomorphism algebras of its isoclinic parts, and in (\ref{eq:embdv}) $D_v$ needs to embed into the endomorphism algebra of each isoclinic part. Let $\bH_{v,A,\lambda}$ be the isoclinic part of slope $\lambda$, then $D_v\hookrightarrow \End(\bH_{v,A,\lambda})$. %
\begin{lemma}
Notation as above, we have \begin{equation}\label{eq:d[]|h}
    d\cdot[F_v:\bQ_p]\, |\, \dim_{K}\bH_{v,A,\lambda}.
\end{equation}
\end{lemma}
\begin{proof}
   Let $[F_v:\bQ_p]=n$ and $\, \dim_{K}\bH_{v,A,\lambda}=m$. Then $\End(\bH_{v,A,\lambda})$ is a central simple algebra over $\bQ_p$ of degree $m$, and $D_v\subseteq \End(\bH_{v,A,\lambda})$ is a central simple algebra over $F_v$ of degree $d$. The lemma follows from the double centralizer theorem (cf. \cite[tag 074T]{stacks-project}): Let $B:=C_{\End(\bH_{v,A,\lambda})}(F_v)$ be the centralizer of $F_v$ in $\End(\bH_{v,A,\lambda})$, then $B$ is a central simple algebra over $F_v$ of degree $m/n$. Note that $D_v\subseteq B$, so $C_{B}(D_v)$, the centralizer of $D_v$ in $B$, is a central simple algebra over $F_v$. As a result $d$ needs to divide $m/n$. So $dn|m$. 
\end{proof}

\begin{lemma}\label{lm:field_lemma}
 Suppose that $g=4$. If $\xi$ has any slope of odd height, then $D$ is a field.
\end{lemma}
\begin{proof}
 If $D$ is not a field, then by Albert classification $D$ is a quaternion over its center, so $d=2$. Let $\lambda$ be the slope that has odd height. Then there is at least one $v$ such that $\dim_{K}\mathbb{H}_{v,A,\lambda}$ is odd. This contradicts   (\ref{eq:d[]|h}). 
\end{proof}

\begin{theorem}\label{eq:Newton--endo}
    Let $A$ be a non-ordinary simple non-CM abelian fourfold over a characteristic $p$ algebraically closed field. Let $D=\End^0(A)$. The following table lists possible types of $D$ associated to a given Newton polygon $\xi$,

\begin{table}[H]
\centering
\renewcommand{\arraystretch}{1.35}
\begin{tabular}{c|c}
\toprule
$\xi$& Type of $D$  \\
\midrule
$\{0^3,(\tfrac12)^2,1^3\}$ & I(1), I(2), I(4)\\
$\{0,(\tfrac13)^3,(\tfrac23)^3,1\}$ & I(1), I(4), IV(1,1)\\  
$\{0,(\tfrac12)^6,1\}$& I(1), I(2), I(4)\\  
$\{(\tfrac13)^3,(\tfrac12)^2,(\tfrac23)^3\}$& I(1), I(4)\\ \hline 
$\{0^2,(\tfrac12)^4,1^2\}$ & I(1), I(2), I(4), II(1),\\
&II(2), III(1), IV(1,1), IV(2,1)\\\hline
$\{(\tfrac14)^4,(\tfrac34)^4\}$ & I(1), I(2), I(4), II(2),\\
&IV(1,1), IV(2,1)\\\hline 
\end{tabular}
\caption{Type of $D$ vs. Newton polygon}\label{tb:Endo4Newton}
\end{table}

\end{theorem}
\begin{proof}
\begin{enumerate}
    \item $\{0^3,(\tfrac12)^2,1^3\}$.
 It follows from Lemma~\ref{lm:field_lemma} that $D$ is a field.  We only need to exclude the possibility that $D$ is a CM field. Suppose that $D$ is an imaginary quadratic extension of $\bQ$. If $p$ does not split in $D$, let $v$ be the unique prime lying above $p$. Then by (\ref{eq:d[]|h}) we have $2=[D_{v}:\bQ_p]$ dividing the height of the slope 0 part of 
    $\mathbb{H}_{p,A}$, which is 3. Contradiction. If $p$ splits, then Lemma~\ref{lm:RosatiD} tells that $\mathbb{H}_{p,A}$ splits into two mutually dual $F$-subcrystals, which cannot happen. Now suppose that $D$ is a quartic CM field with the RM subfield $D^\dagger$. If there is a unique prime $\mathfrak{p}$ of $D^\dagger$ lying above $p$, then we use a similar argument   as above (by discussing whether $\mathfrak{p}$ splits in $D$) to arrive at a contradiction. If $p$ splits as $\mathfrak{p}_1\mathfrak{p}_2$ in $D^\dagger$, then Lemma~\ref{lm:split} and \ref{lm:RosatiD} implies that $\mathbb{H}_{p,A}$ splits into two rank 4 self-dual pieces $\mathbb{H}_{\mathfrak{p}_i,A}$, $i=1,2$. There is exactly one summand, say $\mathbb{H}_{\mathfrak{p}_1,A}$, whose Newton polygon is $\{0,(\frac{1}{2})^2,1\}$. Since $\mathbb{H}_{\mathfrak{p}_1,A}$ does not split into two mutually dual summands, there is a unique prime $v$ of $D$ that lies above  $\mathfrak{p}_1$. Then by (\ref{eq:d[]|h}), $2=[D_v:\bQ_p]$ divides $\dim_{K}\mathbb{H}_{v,A,0}=1$. Contradiction. This rules out the possibility that $D$ is a quadric CM field.

\item $\{0,(\tfrac13)^3,(\tfrac23)^3,1\}$.      It follows from Lemma~\ref{lm:field_lemma} that $D$ is a field. We first rule out the case where $D$ is a quartic CM field. Suppose $D$ is quartic CM, and let $D^\dagger$ be the RM subfield. Since $\mathbb{H}_{p,A}$ cannot split into two rank 4 summands, there is a unique prime $\mathfrak{p}$ of $D^\dagger$ lying above $p$ by Lemma~\ref{lm:split}. Let $v$ be any prime of $D$ lying above $\mathfrak{p}$, we have $[D_v:\mathbb{Q}_p]=2$ or $4$. This soon runs into trouble with (\ref{eq:d[]|h}) since $\mathbb{H}_{p,A}$ contains height 1 slope 0 (or 1) part. Using a similar argument we can show that $D$ is not a quadratic real extension. We leave this to the reader. 

\item $\{0,(\tfrac12)^6,1\}$. It follows from Lemma~\ref{lm:field_lemma} that $D$ is a field. Since $\mathbb{H}_{p,A}$ cannot split into two mutally dual summands, and contains odd height slope, $D$ cannot be a quadratic imaginary extension. If $D$ is quartic CM, then $p$ need to split in $D^\dagger$. Say $p=\mathfrak{p}_1\mathfrak{p}_2$. Then this corresponds to the splitting of $\mathbb{H}_{p,A}$ into two self-dual pieces $\mathbb{H}_{\mathfrak{p}_1,A}\oplus \mathbb{H}_{\mathfrak{p}_2,A}$. One of the pieces, say $\mathbb{H}_{\mathfrak{p}_1,A}$, has Newton polygon $\{0,(\frac{1}{2})^2,1\}$. The rest of the proof is identical to what we have in (1).

\item 
$\{(\tfrac13)^3,(\tfrac12)^2,(\tfrac23)^3\}$. It follows from Lemma~\ref{lm:field_lemma} that $D$ is a field. Note that the awkward shape of $\xi$ prevents its splitting into two rank 4 mutually dual summands nor two rank 4 self-dual summands. It is then easy to see that $D$ cannot be a quadratic extension. If $D$ is quartic CM, then $p$ cannot split in $D^\dagger$, and this will again violate (\ref{eq:d[]|h}).  
\item $\{0^2,(\tfrac12)^4,1^2\}$. We explain why IV(2) cannot happen. Let $F$ be the center of $D$, which is an imaginary quadratic extension. If $p$ does not split in $F$, let $v$ be the unique prime in $F$ over $p$. Then $d\cdot [F_v:\bQ_p]$ in (\ref{eq:d[]|h}) is $4$, while $\dim_{K}\bH_{v,A,0}=2$. Since $4\nmid 2$, this case cannot happen. If $p$ splits in $F$, say $p=v\overline{v}$. Then up to swapping the conjugates the Newton polygon of $\mathbb{H}_{v,A}$ can only be $\{0^2,(\frac{1}{2})^2\}$ or $\{0,(\frac{1}{2})^2,1\}$. The second is obvious impossible by (\ref{eq:d[]|h}). For the first one, $D_v$ needs to embed in both $\End(\bH_{v,A,1/2})$ (which is a non-split quaternion) and $\End(\bH_{v,A,0})$ (which is $\mathrm{M}_2(\bQ_p)$). This is again impossible.
\item $\{(\tfrac14)^4,(\tfrac34)^4\}$. Let $D_4$ be the endomorphism algebra of the simple slope-$\frac{1}{4}$ $F$-isocrystal $\bH_{p,A,\frac{1}{4}}$. It is a degree 4 central simple division $\bQ_p$-algebra that has order 4 in $\mathrm{Br}(\bQ_p)$. We claim that $D_4$ cannot contain a central simple $\bQ_p$-algebra of degree 2. In fact, if $B\subseteq D_4$ is a degree 2 central simple $\bQ_p$-algebra, then by double centralizer theorem \cite[tag 074U]{stacks-project}, we have $D_4\simeq B\otimes C$ (where $C$ is another central simple $\bQ_p$-algebra of degree 2). Since $B\otimes C$ has order at most 2 in $\mathrm{Br}(\bQ_p)$, we arrive at a contradiction. 

This implies that $D$ cannot be of type II(1), III(1), nor IV(2) when $p$ splits in the imaginary quadratic field. The case left is type IV(2) when $p$ does not split in the imaginary quadratic field. Then the whole $D_{\bQ_p}$ needs to embed into $D_4$. We will show that $D$ contains a quaternion  $D_{0}$ over $\bQ$, then the non-existence of an embedding $D_{0,\bQ_p}\hookrightarrow D_4$ shown in the last paragraph will bring us the desired contradiction. To show that such $D_{0}$ exists, it suffices to look at the usual quaternion involution $\iota$ on $D$ (which acts trivially on the center). Let $\dagger'$ be the composition of $\iota$ with the Rosati involution $\dagger$. Define $D_0$ to be the invariant subalgebra of $\dagger'$. One can check that $D_0$ is a quaternion over $\bQ$ by looking at the base change $D\otimes_{\bQ}\overline{\bQ}$ and explicitly working with matrix algebras. 
\end{enumerate}  
 
\end{proof}
\subsection{Signature of quadratic CM}\label{rmk:signatureissue}
Let $D$ be an imaginary quadratic field. In characteristic 0, the signature of a $D$-multiplication on an abelian variety is an (unordered) pair of integers recording the dimensions of the splitting of its Lie algebra under the $D$-action. Roughly speaking, abelian varieties  with signature $(m,n)$ are classified by $\mathrm{U}(m,n)$ Shimura varieties. For abelian fourfolds in characteristic 0, the signature can either be $(1,3)$ or $(2,2)$; cf. \cite{MZ1995}. 

The discussion in this section generalizes to the case where $D$ is a CM field. Since a general theory won't be used, we choose not to complicate  things. 
\subsubsection{Signature over $p$-adic DVR}\label{subsub:signature_DVR} Now let $R$ be a $p$-adic discrete valuation ring with fraction field $K$. Let $\widetilde{\mathscr{G}}$ be a $p$-divisible group over $R$ with $D$-multiplication. Then we can again talk about the signature of $\widetilde{\mathscr{G}}$: it is an unordered pair of integers determined by the dimensions of the two pieces of $t_{\widetilde{\mathscr{G}}}(\mathbb{C}_p)$ splitted via the $D$-multiplication. When $p$ splits in $D$,  $\widetilde{\mathscr{G}}$ correspondingly splits as $\widetilde{\mathscr{G}}_{1}\oplus \widetilde{\mathscr{G}}_{2}$ up to isogeny. The signature of $\widetilde{\mathscr{G}}$ in this case is simply $(\dim \widetilde{\mathscr{G}}_{1},\dim\widetilde{\mathscr{G}}_{2})$. 

Let $\widetilde{A}$ be an abelian scheme over $R$ with $D$-multiplication. The signature of $\widetilde{A}_{K}$ equals the signature of its $p$-divisible group $\widetilde{\mathscr{G}}=A[p^{\infty}]$. This is an easy consequence of the Hodge--Tate decomposition: 
$(\Lie \widetilde{\mathscr{G}})\otimes \mathbb{C}_p\simeq (\Lie \widetilde{\mathscr{G}}^\vee)\otimes \mathbb{C}_p$ identifies with 
$H^1(\widetilde{A}_{\bC_p},\cO_{\widetilde{A}_{\bC_p}})=\Lie \widetilde{A}_{\bC_p}$.

\begin{example}\label{ex:first_ex_of_diagonal_ss}
    This example will be relevant to us. Let $k$ be an algebraically closed field of characteristic $p$, and let $R=W(k)$. Suppose that $p$ is inert in $D$. Let $\widetilde{\mathscr{H}}_{1/2}$ be an one dimensional $p$-divisible group over $R$ with $D$-multiplication and supersingular reduction. Such $\widetilde{\mathscr{H}}_{1/2}$ exists, say, by Lubin--Tate theory. Let $\widetilde{\mathscr{G}}=\widetilde{\mathscr{H}}_{1/2}^2$. There are two ways that $D$ can embed into $\End^0(\widetilde{\mathscr{G}})=\mathrm{M}_2(D_{\bQ_p})$. The first  is the diagonal embedding $d\rightarrow \diag(d,d)$, and the second is $d\rightarrow \diag(d,\overline{d})$. They are two different $D$-multiplications on  $\widetilde{\mathscr{G}}$: the first has signature $(2,0)$, while the second has signature $(1,1)$. 
\end{example}

\subsubsection{Signature in characteristic $p$}\label{subsub:sig_charp} We will assume that $p$ splits in $D$. Let $k$ be a field of characteristic $p$ and let $\mathscr{G}$ be a $p$-divisible group over $k$ with $D$-multiplication. Then $\mathscr{G}$ splits up as $\mathscr{G}_1\oplus \mathscr{G}_2$ up to isogeny. The signature is again defined to be $(\dim\mathscr{G}_1,\dim\mathscr{G}_2)$. If $A$ is an abelian variety over $k$, then its signature is defined to be the signature of its $p$-divisible group. The following lemma is immediate from the definition, so we omit the proof.
\begin{lemma}\label{lm:signature_when_split}
Let $R$ be a $p$-adic DVR with residue field $k$, and suppose that $p$ splits in $D$. If there is an abelian scheme $\widetilde{A}$ with $D$-multiplication whose special fiber is $A$, then the signature of $\widetilde{A}$ equals the signature of $A$.  
\end{lemma}

\begin{example}
Let $k=\overline{k}$ and $A/k$ be an abelian fourfold. Suppose $p$ splits as $ \mathfrak{p}\overline{\mathfrak{p}}$ and $\mathfrak{p}\ne \overline{\mathfrak{p}}$ in $D$. Then as we have noted earlier, $\bH_{p,A}$
 splits as mutually dual pieces $\bH_{\mathfrak{p},A}\oplus \bH_{\overline{\mathfrak{p}},A}$. The signature of $A$ is then $(\dim \bH_{\mathfrak{p},A},\dim \bH_{{\overline{\mathfrak{p}},A}})$ (here the dimension is not $K$-dimension, but the dimension of the corresponding $p$-divisible group). Note that $\dim \bH_{\mathfrak{p},A}$ is just the sum of all slopes in its slope multiset. For example, when $\xi=\{0,(\tfrac13)^3,(\tfrac23)^3,1\}$, if the slopes of $\bH_{\mathfrak{p},A}$ are $\{0,(\tfrac13)^3\}$ (resp. $\{0,(\tfrac23)^3\}$), then the signature is $(1,3)$ (resp. $(2,2)$).
\end{example}

\section{Mod $p$ Mumford--Tate conjectures} %
Let $A$ be a principally polarized abelian variety defined over a finitely generated field $L$ of characteristic $p$, and let $l\neq p$ be a prime. Let 
$W:=W(\overline{L})$, with an embedding $W\hookrightarrow \bC$. Let $\mathrm{S}(A)$ be the set of principally polarized abelian schemes over $\overline{W}$ whose special fiber is isogenous to $A_{\overline{L}}$. 
We usually denote by $\widetilde{A}$ an element in $\mathrm{S}(A)$. For $\widetilde{A}\in \mathrm{S}(A)$, 
we also define $\mathcal{S}_{\widetilde{A}}$ to be the smallest special subvariety of $\mathcal{A}_{g}$ containing the moduli point corresponding to $\widetilde{A}_\bC$. The generic Mumford--Tate group of $\mathcal{S}_{\widetilde{A}}$ is isomorphic to $G_{\mathrm{B}}(\widetilde{A}_\bC)$.
\begin{lemma}\label{lm:embeding}
For any $l\neq p$, there is an embedding $G_l(A)^\circ\hookrightarrow G_{\mathrm{B}}(\widetilde{A}_\bC)\otimes \bQ_l$ intertwining their canonical representations, canonical up to conjugation. Similarly, there is an embedding $G_p(A)^\circ\hookrightarrow G_{\mathrm{B}}(\widetilde{A}_\bC)\otimes \Breve{\bQ}_p$ intertwining their canonical representations, canonical up to conjugation.
\end{lemma}
\begin{proof}
Possibly base changing to a finite extension of $L$, we spread $A$ out to some smooth geometrically connected
$X_0$ whose generic point we denote by $\eta$. Let $\mathscr{S}$ be the Zariski closure of $\mathcal{S}_{\widetilde{A}}$ in $\mathscr{A}_{g, \overline{\bZ}}$, and let $\widetilde{\mathscr{S}}$ be its normalization. We will regard  $\mathscr{S}$ and $\widetilde{\mathscr{S}}$ as defined over a suitable number field. Let $\epsilon$ be the generic point of $\mathcal{S}_{\widetilde{A}}$ with a geometric point $\overline{\epsilon}$. Possibly replacing $X_0$ by a finite cover, we can assume that the map $\eta\rightarrow \mathscr{A}_{g,\Fpbar}$ factors through $\widetilde{\mathscr{S}}$. We have $G_l(A)^\circ=G(\mathbb{H}_{l,\eta},\overline{\eta})^\circ$. On the other hand, we also have the following relations between several étale fundamental groups: 
$$ {\pi_1^{\et}(\eta,\overline{\eta})}\rightarrow \pi_1^{\et}(\widetilde{\mathscr{S}},\overline{\eta}) \xleftarrow{\sim}\pi_1^{\et}(\widetilde{\mathscr{S}},\overline{\epsilon})\xleftarrow{g} \pi_1^{\et}(\epsilon,\overline{\epsilon}).$$
The map $g$ is surjective since $\widetilde{\mathscr{S}} $ is normal (\cite[0BQM]{stacks-project}). Recall that the monodromy groups are the Zariski closures of the images of étale fundamental groups. Choosing appropriate identifications of the fibers, we see that $G(\bH_{l,\eta},\overline{\eta})^\circ\subseteq G(\mathbb{H}_{l,\epsilon},\overline{\epsilon})^\circ$. The later group is known to be contained in $G_{\mathrm{B}}(\widetilde{A}_{\bC})\otimes {\bQ_l}$; cf. \cite[I]{DP1982}. 

Now we treat the crystalline case. Let $\Fpbar_q$ be a sufficiently large finite field. By abuse of notation, we will also write $\mathscr{S}_{\Fpbar_q}$ for its reduced induced structure. By \cite[Lemma 4.8]{Jiang2023}, possibly replacing $X_0$ by a smooth variety dominating it, we can assume that $X_0$ factors through a component $Y_0$ of a smooth alteration of  $\mathscr{S}_{\mathbb{F}_q}$. It suffices to show that $G(\bH_{p,Y_0})^\circ\subseteq G_\mathrm{B}(\widetilde{A}_\bC)\otimes \Breve{\bQ}_p$. 
By \cite[Lemma 4.8]{Jiang2023} again, we can shrink $Y_0$ without change the monodromy. Therefore possibly shrinking $Y_0$ we can assume that $Y_0$ maps to a smooth open subscheme $U_0\subseteq {\mathscr{S}}_{\mathbb{F}_q}$. Let $x\in U_0(\Fpbar)$, it suffices to show that $G(\bH_{p,U_0},x)^\circ\subseteq G_\mathrm{B}(\widetilde{A}_\bC)\otimes \Breve{\bQ}_p$.  
Now we can use the characterization of $F$-isocrystals due to Berthelot and Ogus; cf. \cite{BO1983}. Let $]U_0[$ be the rigid tubular neighborhood of $U_0$ in $(\mathcal{S}_{\widetilde{A}})^{\mathrm{an}}_{\bQ_q}$. The  
$F$-isocrystal $\bH_{p,U_0}$ is a Frobenius bundle with connection $\mathcal{H}$ over $ ]U_0[$. The Hodge cycles on 
$(\mathcal{S}_{\widetilde{A}})^{\mathrm{an}}_{\bQ_q}$ defining $G_\mathrm{B}(\widetilde{A}_\bC)$ induce 
(étale local) Frobenius invariant horizontal sections $\{s_\alpha\}\subseteq \mathcal{H}^\otimes$ cutting out 
$G_\mathrm{B}(\widetilde{A}_\bC)\otimes \Breve{\bQ}_p$ (cf. \cite[Theorem 5.3]{Bl1994}) and are overconvergent. Since $]U_0[\subseteq (\mathcal{S}_{\widetilde{A}})^{\mathrm{an}}_{\bQ_q}$, the monodromy group $G(\bH_{p,U_0},x)^\circ$ must preserve $\{s_\alpha\}$. Therefore $G(\bH_{p,U_0},x)^\circ\subseteq  G_\mathrm{B}(\widetilde{A}_\bC)\otimes \Breve{\bQ}_p$, as desired. 
\end{proof}

The notion of (weak) Mumford--Tate conjecture is already introduced in Definition~\ref{def:strong MT}.

\begin{proposition}\label{prop:relationbetweencycleconjectures} Assume that $A$ satisfies the Mumford--Tate conjecture and the Hodge conjecture holds for $\widetilde{A}_\bC$. Then the Tate conjecture holds for $A$.   
\end{proposition}
\begin{proof} There is an  $\widetilde{A}\in \mathrm{S}(A)$ such that $G_l(A)^\circ\hookrightarrow G_{\mathrm{B}}(\widetilde{A}_{\bC})\otimes \bQ_l$ is an isomorphism. Let $\mathscr{H}^i(\widetilde{A}_{\bC})\subseteq H^{2i}_\mathrm{B}(\widetilde{A}_{\bC})(i)$ be the subspace fixed by the Mumford--Tate group $G_{\mathrm{B}}(\widetilde{A}_{\bC})$. The Hodge conjecture implies that $\mathscr{H}^i(\widetilde{A}_{\bC})$ is spanned by algebraic cycles. These cycles are defined over $\overline{K}=\overline{W}[1/p]$, so we can spread them out to algebraic cycles over $\overline{W}$ and take reduction modulo $p$ to get algebraic cycles on $A_{\overline{L}}$. To prove that these algebraic cycles span the Tate classes $\mathscr{T}_l^i(A_{\overline{L}}):=\lim_{L\subseteq L'\subseteq \overline{L}} \mathscr{T}_l^i(A_{L'})$, we simply note that $\mathscr{T}_l^i(A_{\overline{L}})$ is the subspace of $H^{2i}_{\text{ét}}(A_{\overline{L}},\bQ_l)(i)$ fixed by $G_l(A)^\circ\simeq G_{\mathrm{B}}(\widetilde{A}_{\bC})\otimes \bQ_l$. 
\end{proof}

\subsection{Producing elements in $\mathrm{S}(A)$}\label{sub:PSE}  We present a method of producing reasonably small elements in $\mathrm{S}(A)$. If we are lucky enough, we would find among them an element that verifies the (weak) Mumford--Tate conjecture.

\begin{construction}\label{const:1}%
We discuss how to find a small $\widetilde{A}\in \mathrm{S}(A)$ in practice. From Serre--Tate theory, it suffices to consider the much more tractable question of lifting the polarized $p$-divisible group $(\mathscr{G}:=A_{\overline{L}}[p^\infty],\psi)$ up to isogeny. For a slope $0\leq \mu\leq 1$, let $\mathscr{H}_{\mu}$ be the 
simple slope-$\mu$ $p$-divisible group over $\overline{L}$. Dieudonné--Manin classification  tells us that $\mathscr{G}$ is isogenous to $\bigoplus_{\mu} \mathscr{H}_\mu^{m_\mu}$, which is equipped with a quasi-polarization $\psi$ induced from that on $\mathscr{G}$. Note that $\psi$ decomposes into quasi-polarizations $\psi_{\mu}$ on $\mathscr{H}_\mu^{m_\mu}\oplus \mathscr{H}_{1-\mu}^{m_{1-\mu}}$ (resp. $\mathscr{H}_\mu^{m_\mu}$) when $\mu<\frac{1}{2}$ (resp. $\mu=\frac{1}{2}$).

To lift $(\mathscr{G},\psi)$ up to isogeny, it suffices to lift each $\mathscr{H}_\mu^{m_\mu}$ separately to $\overline{W}$ up to isogeny and taking direct sum (more care is needed for the quasi-polarization to lift). Lifting connected $p$-divisible groups of pure slope can be done by Lubin--Tate theory or Grothendieck--Messing theory, and is relatively well understood. This often helps us to construct a lift of $(\mathscr{G},\psi)$ up to isogeny with controlled endomorphism algebra. Finally, use a further isogeny to descend the polarization to a  principal one. By Serre--Tate theory, we get a principal polarized abelian scheme $\widetilde{A}$ over $\overline{W}$ whose special fiber is isogenous to $A_{\overline{L}}$. 

\end{construction}
Let $D=\End^0(A_{\overline{L}})$. In practice, one can often lift $\mathscr{G}$ up to isogeny so that a subalgebra $D'\subseteq D$ lifts. %
Towards this end, the following is particularly useful: 
\begin{proposition}[Chai--Conrad--Oort]\label{prop:lifting_field}
Suppose that $D$ is a field. Up to an isogeny, there is a lift $\widetilde{A}$ such that  $\End^0(\widetilde{A})=D$. 
\end{proposition}
\begin{proof} 
This is well-known. We briefly sketch the proof. When $D$ is of maximal degree, one can approach this via Serre--Tate theory and CM lifting for $p$-divisible groups; cf. \cite[Proposition 2.4]{CCO2013}. More precisely, one uses \cite[Proposition 2.4(4)]{CCO2013} to construct a lift of the $p$-divisible group up to isogeny such that $D_{\bQ_p}$ also lifts, then uses \cite[Proposition 2.4(1)]{CCO2013} to construct a rational lift of polarization. When $D$ has a smaller degree, we can find a maximal $\dagger$-invariant étale $\bQ_p$-subalgebra $\mathcal{K}\subseteq \End^0(A_{\overline{L}}[p^\infty])$  containing $D_{\bQ_p}$. We then reduce to the case of CM lifting for $p$-divisible groups. A similar argument using \cite[Proposition 2.4]{CCO2013} shows that up to isogeny $A_{\overline{L}}[p^\infty]$ lifts to a polarized $p$-divisible group such that $\mathcal{K}$ also lifts. In particular, $D$ lifts as well. Apply Serre--Tate theory to conclude. 
\end{proof}

\section{Lifting simple abelian fourfolds}\label{sec:lsaf}
This section logically follows \S\ref{sub:PSE}. We will construct good lifts of simple non-CM abelian fourfolds, producing small elements in $\mathrm{S}(A)$. The small elements will then be used in subsequent chapters to prove that a large family of such abelian fourfolds satisfy the Mumford--Tate conjecture. To notation will be consistent with the last section: $A$ is used to denote a simple non-CM abelian fourfold over a characteristic $p$ algebraically closed field $L=\overline{L}$ and $\mathscr{G}=A[p^\infty]$. We write $D$ for the endomorphism algebra $\End^0(A)$. Let $W=W(\overline{L})$ and $K=W[1/p]$, and let $\overline{W}$ be the integral closure of $W$ in $\overline{K}$. The main result of this section is 
\begin{theorem}
\label{thm:lift_endo_fourfold}
  Notation as above. Then\begin{enumerate}[label=(\alph*)]
      \item  there is an element $\widetilde{A}\in \mathrm{S}(A)$ such that  $\End^0(\widetilde{A})=D$. 
      \item if $D$ is an imaginary quadratic field where $p$ does not split, then we can require that the $D$-multiplication on $\widetilde{A}$ in (a) to have signature $(2,2)$.
  \end{enumerate} 
\end{theorem} 
For (a), the liftability follows from Proposition~\ref{prop:lifting_field} when $D$ is a field.  When $D$ is not a field, then by Theorem~\ref{eq:Newton--endo}, we are in two cases:\begin{enumerate}
   
    \item  $\xi=\{0^2,(\frac{1}{2})^4,1^2\}$ and $D$ is of type II(1), II(2), or III(1), 
    \item  $\xi=\{(\frac{1}{4})^4,(\frac{3}{4})^4\}$ and $D$ is of type II(2). 
\end{enumerate} 
So this is a problem of \textit{quaternionic lifting}. 

For (b),  by Theorem~\ref{eq:Newton--endo} the only Newton polygons that support IV(1,1) type endomorphisms are $$\{0^2,(\frac{1}{2})^4,1^2\}, \{(\frac{1}{4})^4,(\frac{3}{4})^4\} \text{ or } \{0,(\frac{1}{3})^3,(\frac{2}{3})^3,1\}.$$ But for the last Newton polygon, $p$ need to split in $D$. So again, we are reduced to consider the first two Newton polygons.

We will prove the theorem in next two subsections. 

\begin{notation}
    Notation as in Construction~\ref{const:1}. In the sequel, we will make heavy use of the Dieudonné module attached to $\mathscr{G}$. We briefly set up the notation here. We will use $\mathbb{D}_{\mathrm{cov}}(-)$ for the \textbf{covariant} Dieudonné functor and $\mathbb{D}^0_{\mathrm{cov}}(-)$ for its rationalization. We have the following splitting of polarized rational Dieudonné modules: 
\begin{equation}\label{eq:standard_set_up}
\mathbb{D}_{\mathrm{cov}}^0(\mathscr{G})=  \bigoplus_{\mu} \mathbb{D}_{\mathrm{cov}}^0(\mathscr{H}_\mu)^{\oplus m_\mu}.\end{equation} 
We will again write $\dagger$ for the Rosati involution induced by $\psi$. 

\end{notation}

\subsection{Newton polygon $\{0^2,(\frac{1}{2})^4,1^2\}$}\label{sub:NP-0-1/2-4: lifting} Decompose 
\begin{equation}\label{eq:a_split}
(\mathbb{D}_{\mathrm{cov}}^0(\mathscr{G}),\psi)=\left(\mathbb{D}_{\mathrm{cov}}^0(\mathscr{H}_0^2\oplus \mathscr{H}_1^2),\psi\right)\oplus\left(\mathbb{D}_{\mathrm{cov}}^0(\mathscr{H}_{1/2}^2),\psi\right) .
    \end{equation} Needless to say, the key of the whole lifting problem lies in the supersingular part $(\mathscr{H}_{1/2}^{2},\psi)$. Indeed, to treat the ordinary part $\mathscr{H}_0^2\oplus \mathscr{H}_1^2$, we can just take the canonical lift, and the lifting of polarization and endomorphisms is unobstructed. 

To ease notation, we will just write $ H$ for $\mathbb{D}_{\mathrm{cov}}^0(\mathscr{H}_{1/2})$. We will think of $H$ as a dimension 2 $K$-space equipped with a $\sigma$-linear Frobenius action $F$ and a symplectic pairing $\psi$ that is compatible with the Frobenius. Let $\Delta$ be the unique quaternion division algebra over $\bQ_p$. We have the following identification 
$$\End^0(\mathscr{H}_{1/2}^2)\simeq\End_F(H^2)\simeq \mathrm{M}_2(\bQ_p)\otimes \Delta\simeq \mathrm{M}_{2}(\Delta).$$
The Rosati involution on $\End_F(H^2)$ induced by $\psi$ is still denoted by $\dagger$. 

On the other hand, we can think of $H$ as the $F$-isocrystal of a supersingular elliptic curve, carrying a principal polarization $\psi_{\mathrm{s}}$ whose induced Rosati involution on $\Delta$ is the classical quaternion conjugation $c\rightarrow \overline{c}$. This gives us a \textbf{framing pairing} $\psi_0=\psi_{\mathrm{s}}\oplus\psi_{\mathrm{s}}$, which induces a \textbf{framing involution} on $\mathrm{M}_2(\Delta)$ by $a\rightarrow \overline{a}^t$. 
The following lemma explains the name ``framing'': 
\begin{lemma}\label{lm:refpairing1}
    There is a bijection between the following sets \begin{align*}
   \left\{\begin{aligned}
       & \text{Invertible }Q\in \mathrm{M}_2(\Delta)\\ &\text{ with $Q=\overline{Q}^t$} 
    \end{aligned}\right\}&\xrightarrow{\sim}  \left\{\begin{aligned}
        &\text{Symplectic pairings $\psi$ on $H^2$}  \\
   &\text{with $\psi(F-,F-)=\sigma p\psi(-,-)$}
    \end{aligned}\right\} \\
   Q&\rightarrow  \psi_0(-,Q-).
    \end{align*}
    Under the above correspondence, the involution $\tau$ on $H^2$ induced by $\psi$ can be expressed by conjugating the reference involution, i.e., $\tau(a)=Q^{-1}\overline{a}^t Q$.
\end{lemma}

\begin{proof}
    Easy to check that $\psi_0(-,Q-)$ is a symplectic pairing: the condition $Q=\overline{Q}^t$ is used to guarantee that $\psi_0(-,Q-)$ is alternating. Now we construct an map in the inverse direction: given a $\psi$ from the RHS set, we see that $\psi$ and $\psi_0$ both induce isomorphism of $F$-isocrystals $H^2\rightarrow (H^2)^{\vee}(-1)$. So by composing one isomorphism with the inverse of the other, we find a $Q\in \mathrm{M}_2(\Delta)$ with $\psi(-,-)=\psi_0(-,Q-)$. Since $\psi_0(-,-)$ and  $\psi_0(-,Q-)$ are both alternating, we must have $Q=\overline{Q}^t$.  The assertion about involution is clear.
\end{proof}
The following lemma says that we have finer description of $Q$ base on the type of $D$: 
\begin{lemma}\label{lm:Refpairinglemma2} Let $\psi$ be the polarization pairing on $H^2$, which by Lemma~\ref{lm:refpairing1} corresponds to a $Q$ with $\psi(-,-)=\psi_0(-,Q-)$. Let $I=\begin{pmatrix}
            1&\\
            &1
        \end{pmatrix}$ and $J=\begin{pmatrix}
            &1\\
            -1&
        \end{pmatrix}$, both regarded as matrices in $\mathrm{M}_2(\Delta)$. \begin{enumerate}
    \item If $D_{\bQ_p}=\mathrm{M}_2(\bQ_p)\subseteq \mathrm{M}_2(\Delta)$ and $\dagger$ is an orthogonal involution, then $Q$ is a symmetric matrix in $\mathrm{M}_2(\bQ_p)$.
    \item  If $D_{\bQ_p}=\mathrm{M}_2(\bQ_p)\subseteq \mathrm{M}_2(\Delta)$ and $\dagger$ is a symplectic involution, then there exists a $c\in \Delta$ with $\overline{c}=-c$, such that $Q=cJ$.
    \item If $D_{\bQ_p}=\Delta \cdot I$, and $\dagger$ is an orthogonal involution ($I$ is the identity matrix of $\mathrm{M}_2(\Delta)$). Then there a $c\in \Delta$ with $\overline{c}=-c$, such that 
    $Q=cJ$. 
     \item Suppose that $\iota:L\subseteq \Delta$ is a quadratic field extension of $\bQ_p$. If $D_{\bQ_p}=L\cdot I{\subseteq}\Delta \cdot I$ with $\dagger$ being complex conjugation, then $Q\in \mathrm{M}_2(L){\subseteq}\mathrm{M}_2(\Delta)$ is an Hermitian matrix. 

      \item Suppose that $D_{\bQ_p}\simeq \bQ_p\times \bQ_p$ embedded diagonally in $\mathrm{M}_2(\Delta)$, and $\dagger$ acts by swapping the two summands. Then there is an element $c\in \Delta$ such that $Q=\begin{pmatrix}
          & c\\
          \overline{c} & 
      \end{pmatrix}$. 
\end{enumerate} 
 \end{lemma}
\begin{proof}
For all $a\in D_{\bQ_p}$ we have  $\psi(ax,y)=\psi(x,a^\dagger y)$, and this implies that 
$\psi_0(x,\overline{a}^tQy)=\psi_0(ax,Qy) = \psi_0(x,Qa^\dagger y)$. Therefore \begin{equation}\label{eq:aQa_commute}
\overline{a}^tQ=Qa^\dagger
\end{equation}
Now we begin case by case study: 
    \begin{enumerate}
        \item We can write $a^\dagger=S^{-1}a^tS$ for symmetric $S\in \GL_2(\bQ_p)$. Plug into (\ref{eq:aQa_commute}) to get $${a}^tQ=\overline{a}^tQ=Qa^\dagger= QS^{-1}a^tS.$$
        So $QS^{-1}$ commutes with $a^t$ for all $a\in \mathrm{M}_2(\bQ_p)$. As a result $QS^{-1}$ lies in the centralizer of $\mathrm{M}_2(\bQ_p)$, which is $\Delta \cdot I$. Let $QS^{-1}=c\cdot I$ for some $c\in \Delta$. Then $Q=cS$. On the other hand, we have $Q=\overline{Q}^t$. So $cS=(\overline{cS})^t=\overline{c}S$, hence $c=\overline{c}\in \bQ_p$. Therefore $Q=cS$ is a symmetric matrix in $\GL_2(\bQ_p)$. 
        \item We can write $a^\dagger=J^{-1}a^tJ$. Plug into (\ref{eq:aQa_commute}) to get $${a}^tQ=\overline{a}^tQ=Qa^\dagger= QJ^{-1}a^tJ.$$
         So similar to the previous case, we have $QJ^{-1}=cI$ for some $c\in \Delta$. It follows that $\overline{c}=-c$ and $Q=cJ$. 
         \item There exists some $c\in \Delta$ with $\overline{c}=-c$, such that $a^\dagger=c^{-1}\overline{a}c$. Plug into (\ref{eq:aQa_commute}) to get $$ (\overline{a}I)Q=(\overline{a}^tI)Q=Q(a^\dagger I)= Q(c^{-1}\overline{a}cI).$$ So $Qc^{-1}I$ commutes with $aI$. As a result $Qc^{-1}I\in \mathrm{M}_2(\bQ_p)$. So we can write $Q=cR$ for some $R\in \mathrm{M}_2(\bQ_p)$. Now $Q=\overline{Q}^t$ implies that $R^t=-R$. So $R=rJ$, hence  $Q=rcJ$. Absorb $r$ into $c$ to conclude.   \item  Since the involution induced by $\psi$ and $\psi_0$ both induce conjugation on $L\cdot I$, $Q$ needs to commute with $L\cdot I$. Therefore $Q\in \mathrm{M}_2(L)$. Hermitian follows from the fact that $Q=\overline{Q}^t$. 
         \item This is the easiest one, and is left as an exercise.  
    \end{enumerate}
\end{proof}

 \begin{lemma}\label{lm:QCLL}
     Let $B\subseteq \End_F(H^2)$ be a subalgebra and let $C\subseteq \End_F(H^2)$ be the centralizer of $B$. Suppose that there exists a quadratic subfield $L\subseteq C$ with the property that $a^\dagger=\overline{a}$ for $\forall a\in L$. Then up to isogeny, $(\mathscr{H}^2_{1/2},\psi)$, together with the $B$-multiplication, lifts to $W\cO_L$ (this ring equals $W$ when $p$ is unramified in $L$, and is $W[\sqrt{p}]$ when $p$ is ramified). 
 \end{lemma}
 \begin{proof}
    Fix an embedding $\iota: L\hookrightarrow \Delta$. By Skolem--Noether, the embedding $L\subseteq \End_F(H^2)=\mathrm{M}_2(\Delta)$ is conjugate to the diagonal one $\diag(\iota,\iota)$. Up to a self-isogeny of $\mathscr{H}^2_{1/2}$, we can assume that $L\subseteq \End_F(H^2)$ is exactly the diagonal embedding $\diag(\iota,\iota)$.

Note that the centralizer of $\diag(\iota,\iota)(L)$ is exactly $\mathrm{M}_2(L)$, viewed as a subalgebra of $\mathrm{M}_2(\Delta)$ via the embedding $\mathrm{M}_2(\iota)$. It follows that $B\in\mathrm{M}_2(L)$. Let $Q$ be as in Lemma~\ref{lm:refpairing1} so that $\psi(-,-)=\psi_0(-,Q-)$. 
Since both $\psi$ and $\psi_0$ induce conjugation on $\diag(\iota,\iota)(L)$, we see that $Q$ commutes with $\diag(\iota,\iota)(L)$. It follows that $Q\in \mathrm{M}_2(L)$ as well. 

By Lubin--Tate theory, we can find a $p$-divisible group with quasi-polarization $(\widetilde{\mathscr{H}}_{1/2},\widetilde{\psi}_{\mathrm{s}})$ over $W\cO_L$ lifting $(\mathscr{H}_{1/2},\psi_{\mathrm{s}})$ so that the $L$-multiplication also lifts. Since $B$ and $Q$ both lie in $\mathrm{M}_2(L)$,
they lift to $\widetilde{\mathscr{H}}_{1/2}^2$. It follows that the pairing $\widetilde{\psi}:=\widetilde{\psi}_{0}(-,Q-)$ lives over  $\widetilde{\mathscr{H}}_{1/2}^2$ and lifts $\psi$. As a result, $(\widetilde{\mathscr{H}}_{1/2}^2,\widetilde{\psi})$ is a lift of $({\mathscr{H}}_{1/2}^2,{\psi})$ such that the $B$-multiplication lifts as well. 
 \end{proof}

\begin{proof}[Proof of Theorem~\ref{thm:lift_endo_fourfold}(a) for Newton polygon $\{0^2,(\frac{1}{2})^4,1^2\}$]As we have seen earlier, to solve the lifting problem it suffices to lift the supersingular part $\mathscr{H}_{1/2}^2$ with prescribed polarization and $D_{\bQ_p}$-multiplication.  
\begin{enumerate} 
    \item Type II(1). Note that $D_{\bQ_p}$ embeds into $\End^0(\mathscr{H}_{0}^{2})$, hence $D_{\bQ_p}\simeq \mathrm{M}_2(\bQ_p)$. On the other hand $D_{\bQ_p}\hookrightarrow\End_F(H^2)\simeq \mathrm{M}_2(\Delta)$, and Skolem--Noether implies that the image is conjugate to $\mathrm{M}_2(\bQ_p)$. Up to a self-isogeny of $\mathscr{H}_{1/2}^2$, we can assume that $D_{\bQ_p}$ embeds as exactly $\mathrm{M}_2(\bQ_p)$. Recall that the Rosati involution $\dagger$ on $D_{\bQ_p}$ is orthogonal. Apply Lemma~\ref{lm:Refpairinglemma2}(1) to conclude that $\psi(-,-)=\psi_0(-,Q-)$ for some symmetric $Q\in \mathrm{GL}_2(\bQ_p)$.  

   Take $B=D_{\bQ_p}$ in Lemma~\ref{lm:QCLL}. Then $C=\Delta\cdot I$. We can take $L$ to be any quadratic subfield of $C$. Then $\dagger$ induces conjugation on $C$. Apply Lemma~\ref{lm:QCLL} to conclude.
   \item Type III(1). Similar to the type II(1) case, $D_{\bQ_p}\simeq \mathrm{M}_2(\bQ_p)$ and up to a self-isogeny of $\mathscr{H}_{1/2}^{2}$ we can assume $D_{\bQ_p}=\mathrm{M}_2(\bQ_p)\subseteq \mathrm{M}_2(\Delta)$. The Rosati involution is symplectic. Apply  Lemma~\ref{lm:Refpairinglemma2}(2) to conclude that $\psi(-,-)=\psi_0(-,cJ-)$ for some $c\in \Delta$ with $\overline{c}=-c$.   

    Take $B=D_{\bQ_p}$ in Lemma~\ref{lm:QCLL}. Then $C=\Delta\cdot I$. We can take $L=\bQ_p(c)$. Then $\dagger$ induces conjugation on $C$. Apply Lemma~\ref{lm:QCLL} to conclude. 
    \item Type II(2).  By (\ref{eq:d[]|h}), $p$ must split as two distinct primes $\mathfrak{p}_1\mathfrak{p}_2$ in the center of $D$, leaving us $D_{\bQ_p}=D_{\mathfrak{p}_1}\times D_{\mathfrak{p}_2}$ embedded block diagonally into $\End^0(\mathscr{H}_{0}^{2})\times \End^0(\mathscr{H}_{1/2}^{ 2})$. It is immediate that $D_{\mathfrak{p}_1}$ is a matrix algebra. But $D_{\mathfrak{p}_2}$ can be either matrix or division algebra. The crux now lies in lifting $\mathscr{H}_{1/2}^2$ with the polarization and the $D_{\mathfrak{p}_2}$-multiplication. 

If $D_{\mathfrak{p}_2}$ is a matrix algebra, then the argument in type II(1) case applies. 

If $D_{\mathfrak{p}_2}$ is a division algebra, then it is isomorphic to $\Delta$. By Noether--Skolem, and up to a self-isogeny of $\mathscr{H}_{1/2}^2$, we can assume that $D_{\mathfrak{p}_2}$ embeds diagonally into $\mathrm{M}_2(\Delta)$. Apply  Lemma~\ref{lm:Refpairinglemma2}(3) to conclude that $\psi(-,-)=\psi_0(-,cJ-)$ for some $c\in \Delta$ with $\overline{c}=-c$.  Take $B=D_{\bQ_p}$ in Lemma~\ref{lm:QCLL} and $L=\bQ_p(c)$ to conclude. 
\end{enumerate}
\end{proof}

\begin{proof}[Proof of Theorem~\ref{thm:lift_endo_fourfold}(b) for Newton polygon $\{0^2,(\frac{1}{2})^4,1^2\}$] The goal is to construct a lift $\widetilde{\mathscr{G}}$ together with (rational) lifts of the polarization and the $D$-multiplication, such that $D$ acts via signature $(2,2)$. Again, the ordinary part exhibits no obstruction, and the canonical lift contributes to signature $(1,1)$.  Our goal is to lift the supersingular part $\mathscr{H}_{1/2}^2$ so that it contributes to another $(1,1)$. In the following, we will first explain how to lift $\mathscr{H}_{1/2}^2$ so that $D$ acts with signature $(2,0)$, then ``twist'' one factor by conjugation to obtain another lift with signature $(1,1)$. We encourage the reader to first read Example~\ref{ex:first_ex_of_diagonal_ss} as a warm up (note that the example is much simpler since we didn't consider the polarization). 

Let $L=D_{\bQ_p}$ and fix $\iota: L\subseteq \Delta$. Up to a self-isogeny of $\mathscr{H}_{1/2}^2$, $D_{\bQ_p}$ embeds in $\mathrm{M}_2(\Delta)$ as $\diag(\iota,\iota)$. Apply Lemma~\ref{lm:refpairing1}(4) to show that $Q\subseteq \mathrm{M}_2(L)$ is an Hermitian matrix. Diagonalized $Q$ by a unitary matrix $X\in \GL_2(D_{\bQ_p})$ into $\diag(a,b)$ for some $a,b\in \bQ_p^*$. Therefore, up to a self-isogeny of $\mathscr{H}^{2}_{1/2}$ induced by a scalar multiple of $X$, $$\psi(-,-)=\psi_0(-,\diag(a,b)-)=a\psi_{\mathrm{s}}\oplus b\psi_{\mathrm{s}},$$
while the embedding $D_{\bQ_p}\hookrightarrow \mathrm{M}_2(\Delta)$ stays $\diag(\iota,\iota)$.

If we want to get a lift of ${\mathscr{H}}_{1/2}^2$ with signature $(2,0)$, we simply apply the Lubin--Tate theory to lift $\mathscr{H}_{1/2}$, $\psi_{\mathrm{s}}$ and $\iota (D_{\bQ_p})$ together, then take sum of two copies of such lift.  

On the other hand, a lift of ${\mathscr{H}}_{1/2}^2$ with signature $(1,1)$ requires twisting one factor. 
Let $\mathscr{C}_1:=\mathscr{H}^{2}_{1/2}$ be as above. Let $\overline{\iota}$ be the conjugation of $\iota$. Then Skolem--Noether implies that there exists a $u\in \Delta^*$ such that $\overline{\iota}=u^{-1}{\iota}u$. Up to a scaling we can assume that $u^{-1}\in \End(\mathscr{H}_{1/2})$. Now consider the self-isogeny $$\mathscr{C}_1=\mathscr{H}_{1/2}\times \mathscr{H}_{1/2}\xrightarrow{(\mathrm{id},u^{-1})} \mathscr{H}_{1/2}\times \mathscr{H}_{1/2}.$$
To distinguish, let's denote the second copy of $\mathscr{H}_{1/2}^2$ by $\mathscr{C}_2$. Then the induced embedding $D_{\bQ_p}\hookrightarrow \End^0(\mathscr{C}_{2})\simeq\mathrm{M}_2(\Delta)$ identifies with $\diag(\iota,\overline{\iota})$. The induced polarization $\psi'$ on $\mathbb{D}_{\mathrm{cov}}^0(\mathscr{C}_2)$ is given by the formula $$\psi':=a\psi_{\mathrm{{s}}}\oplus b\psi_{\mathrm{{s}}}(u-,u-)=a\psi_{\mathrm{{s}}}\oplus b\psi_{\mathrm{{s}}}( -,\overline{u}u-) = a\psi_{\mathrm{{s}}}\oplus b\,\mathrm{Nm}(u)\psi_{\mathrm{{s}}}.$$  
By Lubin--Tate theory, we can lift $\mathscr{H}_{1/2}$ to $\overline{W}$ so that $\psi_{\mathrm{s}}$ and $\iota(D_{\bQ_p})$ both lift. Take a sum of two copies of such lift to get a lift $\widetilde{\mathscr{C}}_2$ of $\mathscr{C}_2$ such that $\psi'$ and $\diag(\iota,\overline{\iota})(D_{\bQ_p})$ both lift. Then the $D_{\bQ_p}$-multiplication on $\widetilde{\mathscr{C}}_2$ has signature $(1,1)$ as desired.
\end{proof}

\subsection{Newton polygon $\{(\frac{1}{4})^4,(\frac{3}{4})^4\}$} To ease notation, let $P=\mathbb{D}^0_{\mathrm{cov}}(\mathscr{H}_{1/4})$, $Q=\mathbb{D}^0_{\mathrm{cov}}(\mathscr{H}_{3/4})$. %

\begin{lemma}\label{lm:graph_lemma}
Suppose that there is a non-degenerate symmetric $K$-bilinear form $s:P\times P\rightarrow K$ such that $s(ax,y)=s(x,a^{\dagger}y)$ for all $a\in D_{\bQ_p}$. Then up to isogeny, $(\mathscr{H}_{1/4}\oplus \mathscr{H}_{3/4},\psi)$, together with the $D_{\bQ_p}$-multiplication, lifts to $W$. %
\end{lemma}
\begin{proof} Let $N:=P\oplus Q$. We will apply Grothendieck--Messing theory. In order to do this, we need a two step filtration $$\Fil^\bullet N:\Fil^0N=N\supseteq \Fil^1N$$ such that $\Fil^1N$ is $D_{\bQ_p}$-stable and maximally isotropic with respect to $\psi$. We will construct this filtration from the symmetric form $s$. 

First, since both $\psi$ and $s$ are non-degenerate, there is a unique $K$-linear isomorphism $T:P\rightarrow Q$ with the property that $\psi(x,Ty)=s(x,y)$. For any $a\in D_{\bQ_p}$, we have 
$$\psi(x,Tay) = s(x,ay)= s(a^\dagger x,y) = \psi(a^{\dagger}x,Ty)=\psi(x,aTy). $$
So $Ta=aT$, meaning that $T$ is also $D_{\bQ_p}$-linear. Now define $\Fil^1N$ to be the graph of $T$   \begin{equation}\label{eq:the graph}
\Fil^1N:=\{x+Tx:x\in P\}\subseteq P\oplus Q.
\end{equation}
We claim that $\Fil^1N$ is $D_{\bQ_p}$-stable and  maximally isotropic with respect to $\psi$. Let $a\in D_{\bQ_p}$ and $x\in P$, then $$a(x+Tx)= ax+aTx=(ax)+T(ax)\in \Fil^1N.$$  So $\Fil^1N$ is $D_{\bQ_p}$-stable. On the other hand, since $P$ and $Q$ are all isotropic, we have 
        \begin{align*}
            \psi(x_1+Tx_1, x_2+Tx_2)&= \psi(x_1,Tx_2)+\psi(Tx_1,x_2)\\&=\psi(x_1,Tx_2)-\psi(x_2,Tx_1)\\
            &=s(x_1,x_2)-s(x_2,x_1)=0.
        \end{align*}
Since $\dim \Fil^1N=4$, it is maximally isotropic. 

It is easy to check that $\Fil^\bullet N$ is weakly admissible: the only nontrivial sub-objects are $P$ and $Q$, and $\mathcal{F}_s \cap P=\mathcal{F}_s\cap Q= {0}$. By \cite[Theorem 2.1]{Lf1980}, there exists a strongly divisible $W$-lattice $M\subseteq N$ with $$\Fil^1M=M\cap \Fil^1N,$$ such that 
$F(M)\subseteq M$, $F(\Fil^1M)\subseteq pM$, and $M=F(M)+p^{-1}F(\Fil^1M)$. Therefore, $M$ is a Dieudonné lattice with an one step filtration $M\supseteq \Fil^1M$ such that $$\Fil^1 (M/pM)= VM/pM=\Fil^1M \bmod p.$$ 
By Grothendieck--Messing theory, we get a $p$-divisible group $\widetilde{\mathscr{G}}$ over $W$ from $\Fil^1M$ that lifts $\mathscr{G}$ up to isogeny\footnote{Some care is needed when $p=2$: the divided power of $(p)$ in $W$ is not nilpotent. However, a generalized version of Grothendieck--Messing theory is still available since our $p$-divisible group is connected; cf. \cite[Remark 2.6]{Lau2013}. }. Since $\Fil^1N=\Fil^1M \otimes \bQ_p$ is $D_{\bQ_p}$-stable and is maximally isotropic with respect to $\psi$, we see that 
$\psi$ and the $D_{\bQ_p}$-multiplication both lift rationally to $\widetilde{\mathscr{G}}$.
\end{proof}%

\begin{proof}[Proof of Theorem~\ref{thm:lift_endo_fourfold}(a) for Newton polygon $\{(\frac{1}{4})^4,(\frac{3}{4})^4\}$] The only case that we need to consider is type II(2), and $p$ does not split in the center $E\subseteq D$. Let $\mathscr{K}:=E\otimes K$. Then $P$ is a free $\mathscr{K}$-module of rank 2. We have $$D_{\bQ_p}\otimes K\simeq  \End_{\mathscr{K}}(P)\simeq \mathrm{M}_2(\mathscr{K}).$$ 
The Rosati involution $\dagger$ extends to $\End_{\mathscr{K}}(P)$ and is orthogonal, i.e., there exists an invertible symmetric matrix $S\in  \End_{\mathscr{K}}(P)$ such that $a^\dagger= S^{-1}a^tS$. Then $S$ induces a nondegenerating symmetric bilinear pairing 
$$S:P\times P\rightarrow \mathscr{K}, \,\,(x,y)\rightarrow x^tSy.$$
Note that for $a\in \End_{\mathscr{K}}(P)$, we have $$S(ax,y)=x^ta^{t}Sy = x^t Sa^{\dagger}y=S(x,a^\dagger y).$$
Now define $s:P\times P\rightarrow K$ by $$s=\Tr_{\mathscr{K}/K} S.$$
Then the above computation shows that $s$ is symmetric, nondegenerate and satisfies $s(ax,y)=s(x,a^\dagger y)$. Apply Lemma~\ref{lm:graph_lemma} to win.
\end{proof}
\begin{proof}[Proof of Theorem~\ref{thm:lift_endo_fourfold}(b) for Newton polygon $\{(\frac{1}{4})^4,(\frac{3}{4})^4\}$] 
Here $D=\bQ(\sqrt{d})$ is an imaginary quadratic field where $p$ does not split. Then $\mathscr{K}:=D\otimes K$ is $K(\sqrt{p})$ when $p$ is ramified, or $K\times K$ when $p$ is inert. Picking a basis of $P$ as a free $\mathscr{K}$-module, we can write elements in $P$ as a pair $(x,y)\in \mathscr{K}^2$.  Define a symmetric pairing $s:P\times P\rightarrow K$ by 
$$s((x_1,y_1),(x_2,y_2))= x_1\overline{x}_2+x_2\overline{x}_1+y_1\overline{y}_2+y_2\overline{y}_1.$$
Easy to check that $s$ is nondegenerate. Moreover, $\alpha:=\sqrt{d}\otimes 1\in \mathscr{K}$ acts as
\begin{align*}
    s(\alpha(x_1,y_1),(x_2,y_2))&=\alpha x_1\overline{x}_2-\alpha x_2 \overline{x}_1+\alpha y_1\overline{y}_2-\alpha y_2\overline{y}_1\\
    &=  s((x_1,y_1),\overline{\alpha}(x_2,y_2)).
\end{align*}
So $s(a -, -)=s(-,a^\dagger-)$ for all $a\in \mathscr{K}$. 

Feeding $s$ to Lemma~\ref{lm:graph_lemma}, we get a lift $\widetilde{\mathscr{G}}$ over $W$ up to isogeny such that $\psi$ and $D_{\bQ_p}$-multiplication all lift. We are left to show that the $D_{\bQ_p}$ action on $\widetilde{\mathscr{G}}$ has signature $(2,2)$. For this we need to go back to the construction in the proof: In fact, $\widetilde{\mathscr{G}}$ is constructed using $\Fil^1N$, and Grothendick--Messing theory claims that $$(N/\Fil^1N)_{\bC_p}=(\Lie \widetilde{\mathscr{G}})_{\bC_p}.$$
So to get the signature, it suffices to compute the characteristic polynomial of $\alpha=\sqrt{d}\otimes 1$ on $N/\Fil^1N$. Note that the projection $\Fil^1N\rightarrow P$ is a $D_{\bQ_p}$-linear isomorphism. The characteristic polynomial of $\alpha$ is $(t^2-d)^2$ on both $P$ and $Q$, hence is  $(t^2-d)^4$ on $N$. Therefore the characteristic polynomial of $\alpha$ on $N/\Fil^1N$ is $(t^2-d)^2=(t-\sqrt{d})^2(t+\sqrt{d})^2$, having signature $(2,2)$ as desired. 
\end{proof}
\subsection{More on type IV(2,1)} 
In this section we prove some results about type IV(2,1) which will be used later. The reason that we proof these results in this section is because they are proved using lifting techniques. The following proposition is not restricted to exotic type IV. 
\begin{proposition}\label{lifting mom} Suppose that $A$ has endomorphism algebra $E$ of type IV(2,1) and Newton polygon $\xi=\{0^2,(\frac{1}{2})^4,1^2\}$. Then $p$ is not split in $E_0=E^\dagger$. 
\end{proposition}
We will adopt the notation in \cite[\S 7.5. Type IV(2,1)]{MZ1995}. Let $\rho_1,\rho_2$ be two real embeddings of $E_0$. Let $\sigma_i,\tau_i$ be the two complex embeddings of $E$ extending $\rho_i$. In characteristic 0, abelian fourfolds of type IV(2,1) satisfy a strict signature condition: Let $(n_{\sigma_i},n_{\tau_i})$ be the signature of the $E$-action on such an abelian fourfold along the real embedding $\rho_i$. Then by \cite[(7.5.1)]{MZ1995}, (up to a permutation) we must have: \begin{equation}\label{eq:sigsig}
(n_{\sigma_1},n_{\tau_1})=(2,0),\qquad (n_{\sigma_2},n_{\tau_2})=(1,1).
\end{equation}
\begin{proof} Suppose for a contradiction that $p$ splits in $E_0$. Then $p$ splits as $\mathfrak{p}_1\mathfrak{p}_2$ in $E_0$. Write the decomposition $\bH_{p,A}=\bH_{\mathfrak{p}_1,A}\oplus \bH_{\mathfrak{p}_2,A}$ induced by the splitting of $p$ in $E_0$. One of the summands in this decomposition (say, $\bH_{\mathfrak{p}_2,A}$) must be supersingular: otherwise, since the summands are self-dual of rank four, they must each have Newton polygon $\{0, (1/2)^2, 1\}$. This contradicts (\ref{eq:d[]|h}), by a case-by-case analysis of the possible splitting behaviors of $\mathfrak{p}_2$ in $E$. 
\begin{enumerate}
    \item 
Suppose that $\mathfrak{p}_2$ splits as $\mathfrak{P}_{2}\overline{\mathfrak{P}}_2$ in $E$. 
In the following, we construct a lift $\widetilde{A}$ that violates the signature condition (\ref{eq:sigsig}). We will use the notation in \S\ref{sub:NP-0-1/2-4: lifting}.

First, we canonical lift the ordinary part $\mathscr{H}_0^2\oplus \mathscr{H}_1^2$. We claim that the $E\otimes_{E_0}E_{0,\mathfrak{p}_1}$-action has signature $(1,1)$. When $\mathfrak{p}_1$ does not split, this is clear. When $\mathfrak{p}_1$ splits as $\mathfrak{P}_1\overline{\mathfrak{P}}_1$, then $\bH_{\mathfrak{p}_1,A}$ splits into mutually dual pieces $\bH_{\mathfrak{P}_{1},A}\oplus \bH_{\overline{\mathfrak{P}}_1,A}$. If $\bH_{\mathfrak{P}_{1},A}$ was isoclinic, then the Newton cocharacter factors through the center of $G_p(A)$, contradicting Theorem~\ref{thm:nonvan}. Therefore $\bH_{\mathfrak{P}_{1},A}$ has slope $\{0,1\}$. It then follows that the action also has signature $(1,1)$.

On the other hand, up to a self-isogeny of $\mathscr{H}_{1/2}^2$, the action by $E\otimes_{E_0}E_{0,\mathfrak{p}_2}=E_{\mathfrak{P}_2}\times E_{\overline{\mathfrak{P}}_2}\simeq\bQ_p\times 
   \bQ_p$ on $H^2=\mathbb{D}^0_{\mathrm{cov}}(\mathscr{H}_{1/2}^2)$ is diagonal. The Rosati involution sends $(a,b)$ to $(b,a)$. By Lemma~\ref{lm:Refpairinglemma2}(5), there is a $c\in \Delta$ such that $\psi(-,-)=\psi_0(-,\begin{pmatrix}
       &c\\
       \overline{c}&
   \end{pmatrix}-)$. It then suffices to lift $\mathscr{H}_{1/2}$ so that the action by the field $\bQ_p(c)$ also lifts, then take direct sum. The $E\otimes_{E_0}E_{0,\mathfrak{p}_2}$-action on the lift has signature $(1,1)$. By Serre--Tate theory, we get a lift $\widetilde{A}$ with signature $(1,1)$ along both real embeddings of $E_0$, contradicting (\ref{eq:sigsig}).%

   \item If $\mathfrak{p}_2$ does not split, then there is a unique prime $\mathfrak{P}_{2}$ of $E$ lying above $\mathfrak{p}_2$. By a similar argument as above, the $E\otimes_{E_0} E_{0,\mathfrak{p}_1}$-multiplication lifts to the canonical lift of $\mathscr{H}_0^2\oplus \mathscr{H}_1^2$ and contributes to signature $(1,1)$. On the other hand, the $E_{\mathfrak{P}_{2}}$-multiplication on  $\mathscr{H}_{1/2}^2$ can be lifted with signature $(1,1)$, using exactly the same argument in the proof of Theorem~\ref{thm:lift_endo_fourfold}(b). We again get a lifting $\widetilde{A}$ that violates the signature restriction.
\end{enumerate}
 We have proven that $p$ does not split in $E_0$.
\end{proof}

\begin{lemma}\label{lm:exo_criterion}
Suppose that $A$ is a simple abelian fourfold of type IV. Then \begin{enumerate}
    \item If $A$ is of exotic type IV, then it has Newton polygon $\{0^2,(\frac{1}{2})^4,1^2\}$. 
    \item If $A$ has Newton polygon $\{0^2,(\frac{1}{2})^4,1^2\}$, then $E$ has Galois group $\bZ/2\times \bZ/2$ if and only if $A$ is of exotic type IV. 
    \item  If $A$ is of exotic type IV, there are exactly two imaginary quadratic subfields $k,k'\subseteq E$ such that $p$ splits in $k$ but does not split in $k'$. In addition, there is a lift $\widetilde{A}\in \mathrm{S}(A)$ such that $k'$ lifts and acts with signature $(2,2)$.
\end{enumerate}  
\end{lemma}
\begin{proof}
The only Newton polygons that support IV(2,1) type endomorphisms are $$\{0^2,(\frac{1}{2})^4,1^2\}\text{ or } \{(\frac{1}{4})^4,(\frac{3}{4})^4\}.$$ However, $\{(\frac{1}{4})^4,(\frac{3}{4})^4\}$ does not support an exotic type IV: if there is exotic type IV then the Dieudonné crystal decomposes into two mutually dual parts and the Newton cocharacter factors centrally through $G_p(A)$, contradicting Theorem~\ref{thm:nonvan} (more precisely, to make sense of $G_p(A)$ we need to first descend $A$ to a geometric family: this can be done by taking the Zariski closure of the image of the corresponding moduli map to $\mathscr{A}_{g,\Fpbar_q}$). This proves (1). Part (2) and the $k,k'$ assertion in part (3) are  direct consequences of Proposition~\ref{lifting mom}. 
The existence of the lifting can 
again be obtained via lifting the $p$-divisible group: we canonical lift the ordinary part, and lift the supersinuglar part as in  Theorem~\ref{thm:lift_endo_fourfold}(b) to get such an $\widetilde{A}$ with $k'$-multiplication of signature $(2,2)$.

\end{proof}

\section{Monodromy of simple abelian fourfolds and the Tate conjecture}\label{sec:simple4}
Let $L$ be a finitely generated field over $\bF_q$ and let $A/L$ be a (principally polarized) simple non-CM abelian fourfold. In this section we treat the Mumford--Tate property and the Tate conjecture for $A$. The main theorem is the following
\begin{theorem}\label{thm:MT_simple}
    Notation as in \ref{not:compareMZ}. Suppose $A$ is not of cohomological  Mumford type. Then  
    \begin{enumerate}
          \item If $A$ is of exotic type IV, then $A$ satisfies the weak Mumford--Tate conjecture, but not the Mumford--Tate conjecture. 
            \item In other cases, $A$ satisfies the Mumford--Tate conjecture.
    \end{enumerate}
       
\end{theorem}

\begin{theorem}\label{thm:Tate_simple}
 The Tate conjecture holds for $A$.    
\end{theorem}
The theorems will be proved case by case in a style similar to that of Moonen--Zarhin \cite{MZ1995}. We make the notation slightly different from the previous chapters but closer to \cite{MZ1995} to make it easier for the reader to compare our work to \cite{MZ1995}.

\begin{notation}\label{not:compareMZ}
  \begin{enumerate}
     \item On the $l$-adic side, let $l\neq p$ be a prime. Let $V_l:=\bH_{l,A}$ be the first rational $l$-adic étale cohomology of $A$, equipped with polarization $\varphi_l$. Let $G_l(A)\subseteq \GSp(V_l,\varphi_l)$ be the $l$-adic monodromy and let $\mathfrak{g}_l$ be its Lie algebra. Let $\mathfrak{h}_l\subseteq \mathfrak{g}_l$ be the subalgebra of trace 0, so that $\mathfrak{g}_l=\mathfrak{h}_l\oplus \mathbb{Q}_l\cdot\Id$. Let $\mathbb{C}_l=\widehat{\overline{\bQ}}_l$.
     \item On the crystalline side, 
     we spread out $A$ to an abelian scheme $A_{X_0}$ of constant slope; cf. \S\ref{subsub:spreadout}. We fix a point $x\in X_0(\Fpbar)$.
     Let $V_p$ be the first rational crystalline cohomology of $A_x$, equipped with polarization $\varphi_p$. Let $G_p(A)\subseteq \GSp(V_p,\varphi_p)$ be the $p$-adic monodromy with base point $x$ and let $\mathfrak{g}_p$ be its Lie algebra. Let $\mathfrak{h}_p\subseteq \mathfrak{g}_p$ be the subalgebra of trace 0, so that $\mathfrak{g}_p=\mathfrak{h}_p\oplus \Breve{\mathbb{Q}}_p\cdot\Id$. Let $\mathbb{C}_p=\widehat{\overline{\bQ}}_p$.
     \item On the Betti side, fix an $\widetilde{A}\in \mathrm{S}(A)$ be constructed in Theorem~\ref{thm:lift_endo_fourfold}. More precisely, we require that $\widetilde{A}$ satisfying $\End^0(\widetilde{A})=\End^0(A)$ and, when $A$ is of type IV(1,1) and $p$ does not split in $\End^0(A)$, we further require that $D$ acts on $\End^0(\widetilde{A})$ with signature $(2,2)$.

    Let $V:=\bH_{\mathrm{B},\widetilde{A}_{\bC}}$ be the first Betti cohomology of $\widetilde{A}_{\bC}$, equipped with polarization $\varphi$. Let $G_{\mathrm{B}}(A)\subseteq \GSp(V,\varphi)$ be the generic Mumford--Tate group of $\widetilde{A}_{\bC}$, and let $\mathfrak{g}$ be its Lie algebra. Let $\mathfrak{hg}\subseteq \mathfrak{g}$ be the subalgebra of trace 0, so that $\mathfrak{g}=\mathfrak{hg}\oplus \bQ\cdot \Id$.  We know that $G_l(A)^\circ\subseteq G_\mathrm{B}(A)\otimes \bQ_l$ (resp. $G_p(A)^\circ\subseteq G_\mathrm{B}(A)\otimes \Breve{\bQ}_p$), canonical up to conjugation; cf. Lemma~\ref{lm:embeding}.
     \item We also adopt the convention of \cite[\S 2.6]{MZ1995}. Suppose that $E\subseteq \End^0(A)$ is totally real. Then there exists a $\psi:V\times V\rightarrow E$ whose trace is $\varphi$. We write $\mathfrak{sp}_E(V,\psi)$ as the scalar restriction of the $E$-Lie algebra $\mathfrak{sp}(V,\psi)$ to $\bQ$. Similarly for $\mathfrak{sp}_{E_l}(V_l,\psi_l)$ and $\mathfrak{sp}_{E_p}(V_p,\psi_p)$ in the $l$-adic étale and crystalline cases (but note $E_l$ or $E_p$ may not be a field anymore).
     \par
     Similarly, suppose that $E\subseteq \End^0(A)$ is CM, and the polarization is chosen so that the Rosati involution induces complex conjugation on $E$ (such a choice is always possible according to \textit{loc.cit.}).  Let $i\in E$ be such that $\overline{i}=-i$. Then there is a unique $E$-Hermitian form $V\times V\rightarrow E$ such that $\varphi$ is the trace of $i\psi$. We write $\mathfrak{u}_E(V,\psi)$ as the scalar restriction of the $E$-Lie algebra $\mathfrak{u}(V,\psi)$ to $\bQ$, and similarly for $\mathfrak{u}_{E_l}(V_l,\psi_l)$ and $\mathfrak{u}_{E_p}(V_p,\psi_p)$ in the $l$-adic étale and crystalline cases. We let $U_E$ be the torus over $\bQ$ whose functor of points is
     \begin{gather*}
         U_E(R) = \{x \in (E \otimes_\bQ R)^\times: x \overline{x} = 1\}.
     \end{gather*}
 \end{enumerate}
\end{notation} 
Before we begin a case-by-case analysis in the style of \cite{MZ1995}, we present without a proof a well-known group theoretic lemma that is implicitly used in most of the case by case analysis of \cite{MZ1995}, which adapts to our situation as well. See \cite[\S 2.2]{MZ1995} for more details.
\begin{lemma}[Hazama--Murty]
    Suppose that $A$ is not of type III, and that $\mathfrak{h}_l=\mathfrak{sp}_{\End(A)\otimes \bQ_l}(V_l,\varphi_l) $ (resp. $\mathfrak{h}_p=\mathfrak{sp}_{\End(A)\otimes \Breve{\bQ}_p}(V_p,\varphi_p)$), then the Tate ring of $A^n$ is generated by divisor classes.  
\end{lemma}
We will also use this lemma implicitly. Needless to say, if the Tate ring of $A$ is generated by divisor classes, then the Tate conjecture holds for $A$. 
 
\subsection{Type I}
\subsubsection{Type I(1)}\label{subsub:TypeI(1)} \textit{The geometric endomorphism algebra of $A$ is $\bQ$. Then the Lie algebra $\mathfrak{hg}$ together
with its representation on $V$ (resp. the Lie algebra $\mathfrak{h}_l$ together with its representation on $V_l$) is over $\bC$ (resp. $\bC_l$) is isomorphic to one of the following:\begin{enumerate}
    \item $\mathfrak{sp}_8$ with the standard representation. 
    \item $\mathfrak{sl}_2^3$ with the representation $\mathrm{std}^{\otimes 3}$. 
\end{enumerate}
The Tate ring of $A$ is generated by divisor classes.}
\begin{proof}
  For $\mathfrak{hg}$ the proof of \cite[\S 4.1.Type I(1)]{MZ1995} adapts verbatim. The case where $\mathfrak{h}_l$ is not simple also follows from the discussion of \textit{loc.cit.} So the only case we need to treat is when $\mathfrak{h}_l$ is simple (note that the minuscule weight lemma \cite[Lemma 4.4]{Chi1992} cited in \textit{loc.cit.} does not apply in our case).  To only case that we need to rule out is when geometrically  $\mathfrak{h}_l=\mathfrak{sl}_2$ and the representation $\mathfrak{h}_l\acts V_l$ is isomorphic to $\mathrm{Sym}^7$ (to see why this is the only case: let $ \mathfrak{m}\supseteq \mathfrak{h}_l$ be a proper maximal subalgebra that acts irreducibly, then by complexifying \cite[Table 27]{Graf2020}, one finds that $\mathfrak{m}$ is either $\mathfrak{sl}_2$ acting by $\mathrm{Sym}^7$ or $\mathfrak{sl}_2^{3}$ acting by $\text{std}^{\otimes 3}$, it then follows by our simplicity assumption that $\mathfrak{m}= \mathfrak{h}_l=\mathfrak{sl}_2$ acting by $\mathrm{Sym}^7$). 

  In the following it is easier to work with the $p$-adic monodromy: geometrically again we have $\mathfrak{h}_p=\mathfrak{sl}_2\acts V_p$ as $\mathrm{Sym}^7$. The Newton cocharacter need to factor nontrivially through $\mathfrak{h}_p$. However, $\mathrm{Sym}^7$ has 8 different weights, meaning $A$ has 8 different Newton slopes, which does not happen. 
\end{proof}
\subsubsection{Type I(2)}\textit{The geometric endomorphism algebra of $A$ is a real quadratic field $E$. Then 
$\mathfrak{hg}\simeq \mathfrak{sp}_{E}(V,\varphi)$ (resp.
$\mathfrak{h}_l\simeq \mathfrak{sp}_{E}(V_l,\varphi_l)$), and is geometrically a product of two copies of $\mathfrak{sp}_4$. In particular, for all $n$ the Tate ring of $A^n$ is
generated by divisor classes.}
\begin{proof}
    The assertion for $\mathfrak{hg}$  is the same as \cite[\S 4.2. Type I(2)]{MZ1995}. The assertion for $\mathfrak{h}_l$ is slightly different from \textit{loc.cit.}, since the some of the cited results are only available in characteristic 0. In the following we write down the proof in full. First, we can assume that $l$ splits into $\lambda_1$,$\lambda_2$ in $E$. So $E\otimes\bQ_l$ splits as a product $E_{\lambda_1}\times E_{\lambda_2}=\bQ_l\times \bQ_l$. 
    Correspondingly, the representation of $\mathfrak{h}_l$ on $V_l$ splits as two irreducible 4 dimensional subreps $V_1\oplus V_2$. Similar to \cite[\S 4.2]{MZ1995}, there are symplectic pairings $\psi_i$ on $V_i$ so that $\mathfrak{h}_l\subseteq \mathfrak{sp}(V_1,\psi_1)\times \mathfrak{sp}(V_2,\psi_2)$. It suffices to show that this inclusion is an equality. 
    
    Let $\mathfrak{h}_i$ be the projection of $\mathfrak{h}_l$ to $\mathfrak{sp}(V_i,\psi_i)$. It suffices to show that $\mathfrak{h}_i=\mathfrak{sp}(V_i,\psi_i)$ for each $i$: by Goursat's lemma for Lie algebras and simpleness of $\mathfrak{sp}_4$, $\mathfrak{h}_l$ is either $\mathfrak{sp}(V_1,\psi_1)\times \mathfrak{sp}(V_2,\psi_2)$ or $\mathfrak{sp}(V_1,\psi_1)$ embedded diagonally. But the later is impossible since the endomorphism ring of $A$ would otherwise be much larger than $E$.  

Suppose that $\mathfrak{h}_i\subsetneq \mathfrak{sp}(V_i,\psi_i)$, then it must be a $\bQ_l$-form of $\mathfrak{sl}_2$. In \cite[\S 4.2]{MZ1995},  the authors use \cite[Theorem 0.9]{Zarhin_torsion} to show that $E$ must be of degree at least 4, hence arriving at a contradiction. However, the cited theorem is only for abelian varieties over number fields, whose proof breaks down in characteristic $p$. In the following, we give a modified argument based on \cite[\S 3]{Zarhin_torsion}  \footnote{Our notation used here is different from \textit{loc.cit.}: we use $A$ for the abelian variety whereas \textit{loc.cit.} uses $X$. The field $E_{\lambda_i}$ here is denoted by $E(\lambda_i)$ there, and the $\lambda_i$-adic Tate module $V_i$ here is written as $V(X)(\lambda_i)$ in there.}. Without loss of generality, we will write $i=1$ in what follows. 
Let $G$ be the image of the projection of $G_l(A)^\circ$ in $\GL(V_1)$.  Base change to $\mathbb{C}_l$, the irrep $\mathfrak{h}_1\acts V_1$ becomes an irrep $\mathfrak{sl}_2\acts \mathbb{C}_l^4$, and the purely group theoretic argument in \textit{loc.cit.} implies that there is an open subset $U\subseteq G$ such that each $u\in U$ is a semisimple automorphism of $ V_1$ and all eigenvalues of $u$ are of multiplicity 1 and they constitute a geometric progression in $\mathbb{C}_l$, whose ratio is not a  root of unity. Apply Chebotarev's density theorem for function fields to get a Frobenius element $\mathrm{Fr}(\mathrm{v}')\in U$, where $\mathrm{v}'$ is a prime in a field extension of $L$ over which $A$ has a good reduction $A_{\mathrm{v}'}$ with the same Newton polygon as $A$. The characteristic polynomial $P(t)$ of $\mathrm{Fr}(\mathrm{v}')$ lies in $E[t]$ and is of degree 4. We can write the roots as a geometric progression $\alpha \gamma^{j},j=0,1,2,3$, such that $\gamma$ is not a root of unity. Fix an embedding $\overline{\bQ}\hookrightarrow \bC$. Weil's conjecture implies that $\alpha \gamma^{j}$ has absolute value $\sqrt{q}$, and is a $\lambda$-adic unit when $\lambda$ is a finite place of $E$ not dividing $p$. This in particular implies that $\gamma$ has nontrivial $\mathfrak{p}$-adic valuation for some place $\mathfrak{p}$ in $E$ dividing $p$. In particular, the $\mathfrak{p}$-adic Newton polygon of $P(t)$ has 4 distinct slopes. It means that the Newton polygon of $A_{\mathrm{v'}}$, hence of $A$, has at least 4 distinct slopes. This contradicts the list of Newton polygons with Type I(2) endomorphisms in Theorem~\ref{eq:Newton--endo}. 
\end{proof}

\subsubsection{Type I(4)} \textit{The geometric endomorphism algebra of $A$ is a totally real quartic field $E$. Then 
$\mathfrak{hg}\simeq \mathfrak{sp}_{E}(V,\varphi)\simeq \mathfrak{sl}_{2,E}$ (resp.
$\mathfrak{h}_l\simeq \mathfrak{sp}_{E}(V_l,\varphi_l)\simeq \mathfrak{sl}_{2,E_l}$). In particular, for all $n$ the Tate ring of $A^n$ is
generated by divisor classes.}
\begin{proof}
 The proof of \cite[\S 4.3.Type I(4)]{MZ1995} is group theoretical, and it adapts verbatim to our setting.  
\end{proof}
\subsection{Type II}
\subsubsection{} \textit{The geometric endomorphism algebra of $A$ is  is an indefinite quaternion algebra
$D$ over a totally real field of degree $e\in \{1,2\}$ over $\bQ$. Then $\mathfrak{hg}$ is the centralizer of
$D$ in $\mathfrak{sp}(V,\varphi)$ (resp. $\mathfrak{h}_l$ is the centralizer of
$D$ in $\mathfrak{sp}(V_l,\varphi_l)$ ). In particular, for all $n$ the Tare ring of $A^n$ is generated by divisor classes.}
\begin{proof}
When $e=1$, one can choose a Rosati invariant real quadratic field $E\subseteq D$, and assume that $l$ splits in $E$. Similar to Type I(2) case, $(V_l,\varphi_l)$ splits into $(V_1,\psi_1)\oplus (V_2,\psi_2)$. In addition,
$\mathfrak{sp}_{D\otimes \bQ_l}(V_l,\varphi_l)\subseteq \mathrm{sp}(V_1,\psi_1)\times \mathrm{sp}(V_2,\psi_2)$ is a diagonally embedded subalgebra which is geometrically $\mathfrak{sp}_4$. To show that $\mathfrak{h}_l=\mathfrak{sp}_{D\otimes \bQ_l}(V_l,\varphi_l)$, it suffices to show that the projection of $\mathfrak{h}_l$ to $\mathrm{sp}(V_1,\psi_1)$ is $\mathrm{sp}(V_1,\psi_1)$. If not, then the image of the projection is geometrically $\mathfrak{sl}_2$ acting irreducibly on $V_1\otimes \mathbf{C}$. A similar argument as in Type I(2) case shows that $A$ has at least 4 distinct slopes, which contradicts Theorem~\ref{eq:Newton--endo}.

When $e=2$, write $E$ for the center of $D$ and assume that $l$ splits in $E$. Then again $V_l$ splits into $V_1\oplus V_2$. After base change to $\mathbb{C}_l$, 
$\mathfrak{sp}_{D\otimes \bQ_l}(V_l,\varphi_l)$ becomes $\mathfrak{sl}_2\times \mathfrak{sl}_2$, acting block diagonally on $V_1\oplus V_2$, and each $\mathfrak{sl}_2\acts V_i$ is a direct sum of two copies of standard reps. We must have $\mathfrak{h}_l=\mathfrak{sp}_{D\otimes \bQ_l}(V_l,\varphi_l)$ for otherwise $D$ would be much larger. 
\end{proof}
\subsection{Type III}
\subsubsection{Type III(1)}\textit{The geometric endomorphism algebra of $A$ is a definite quaternion aloebra $D$ over
$\bQ$. Then $\mathfrak{hg}$ is the centralizer of $D$ in $\mathfrak{sp}(V,\varphi)$ (resp. $\mathfrak{h}_l$ is the centralizer of
$D$ in $\mathfrak{sp}(V_l,\varphi_l)$), which is a $\bQ$-form (resp. $\bQ_l$-form) of $\mathfrak{so}_4$. The Tate conjecture holds for $A$. }
\begin{proof}
    The assertions on $\mathfrak{hg}$ and $\mathfrak{h}_l$ follow verbatim from the group theoretic treatment in \cite[\S 6.1]{MZ1995}. This in particular implies that $A$ satisfies the Mumford--Tate conjecture. Now the Hodge conjecture for abelian fourfolds together with Proposition~\ref{prop:relationbetweencycleconjectures} implies the Tate conjecture for $A$. 
\end{proof}

\subsection{Type IV}
\subsubsection{Type IV(1,1)}\label{subsub:TypeIV(1,1)} \textit{The geometric endomorphism algebra of $A$ is  an imaginary quadratic field $k$. Recall that we have fixed specific $\widetilde{A}\in \mathrm{S}(A)$ in Notation~\ref{not:compareMZ}. 
\begin{enumerate}
    \item If $k$ acts on $\widetilde{A}$ with signature $(1,3)$, then $\mathfrak{hg}\simeq \mathfrak{u}(V/k)$ (resp. $\mathfrak{h}_l\simeq \mathfrak{u}(V_l/k_l)$). In particular, the Tate ring of $A^n$ is generated by divisor classes. 
    \item If $k$ acts on $\widetilde{A}$ with signature $(2,2)$, then $\mathfrak{hg}\simeq \mathfrak{su}(V/k)$ (resp. $\mathfrak{h}_l\simeq \mathfrak{su}(V_l/k_l)$). In this case the Tate conjecture for $A$ is true. 
\end{enumerate}}\begin{remark}
    
 Here is a more intrinsic  classification of the $l$-adic monodromy without fixing $\widetilde{A}$ at the beginning: look at all lifts $\widetilde{A}\in \mathrm{S}(A)$ where $k$ also lifts: (1) if there exists a lift $\widetilde{A}$ such that $k$ acts with signature $(2,2)$, then $\mathfrak{h}_l\simeq \mathfrak{su}(V_l/k_l)$, (2) if for all lifts $k$ acts with signature $(1,3)$, then $\mathfrak{h}_l\simeq \mathfrak{u}(V_l/k_l)$. 
\end{remark}  
\begin{proof}
    The assertion for $\mathfrak{hg}$ follows from the same argument as in \cite[\S 7.4. Type IV(1,1)]{MZ1995}. In the following we prove the assertion for $\mathfrak{h}_l$: the proof is very different from (and is much harder than) that in \cite{MZ1995}. Once this is done, we immediately see that $A$ satisfies the Mumford--Tate conjecture, so we can use Hodge conjecture for abelian fourfolds together with Proposition~\ref{prop:relationbetweencycleconjectures} to conclude the Tate conjecture for $A$.

    By Theorem~\ref{eq:Newton--endo}, the Newton polygons that allow type IV(1,1) endomorphism are 
    $$\xi=\{0,(\frac{1}{3})^3, (\frac{2}{3})^3,1\},\;\{0^2,(\frac{1}{2})^4, 1^2\},\;\{(\frac{1}{4})^4, (\frac{3}{4})^4\}.$$ 
    We will treat them case by case: 
\begin{enumerate}[label=(\alph*)]
    \item Suppose that $\xi=\{0,(\frac{1}{3})^3, (\frac{2}{3})^3,1\}$. Then by (\ref{eq:d[]|h}), $p$ splits as $v\overline{v}$ in $k$. Accordingly, $\bH_{p,A}$ splits into two mutually dual pieces $\bH_{v,A}\oplus\bH_{\overline{v},A}$ and $V_p$ splits as $V_{v}\oplus V_{\overline{v}}$. Without loss of generality, suppose that $\bH_{v,A}$ contains slope 0. Our $G_p(A)^\circ$ is a subgroup of $\mathrm{GU}(V_p/k_p)$, which is  generated by the rank 1 torus of homothety $\mathbf{T}_p$ together with a unitary group $\mathrm{U}(V_p/k_p)\simeq\GL(V_{v})\simeq \GL(V_{\overline{v}})$. Now there are two possibilities: 
 \begin{itemize}
        \item \textbf{The (1,3) case}. The slopes of $\bH_{v,A}$ are $\{0,(\frac{1}{3})^3\}$. Then by Lemma~\ref{lm:signature_when_split} the signature of the $k$-action on $\widetilde{A}$ is $(1,3)$. We need to show that $G_p(A)^\circ=\mathrm{GU}(V_p/k_p)$.
        \item \textbf{The (2,2) case}. The slopes of $\bH_{v,A}$ are $\{0,(\frac{2}{3})^3\}$. Then by Lemma~\ref{lm:signature_when_split} the signature of the $k$-action on $\widetilde{A}$ is $(2,2)$. We need to show that $G_p(A)^\circ=\mathbf{T}_p\cdot \mathrm{SU}(V_p/k_p)$.
        \end{itemize}
 Let $G$ be the image of the projection $G_p(A)^\circ\rightarrow\GL(V_{v})$. Since $G_{p}(A)^\circ\supseteq\mathbf{T}_p$, we have $G\supseteq \mathbf{T}_{v}$, where $\mathbf{T}_{v}$ is the torus of homothety in $\GL(V_{v})$. In both cases, the Newton cocharacter $\nu:\bG_m\rightarrow G_p(A)^\circ$ gives rise to a  cocharacter $\nu':\mathbb{G}_m\rightarrow G$ that acts on $V_{v}$ with weights $\{0,1^3\}$. Since $G\supseteq \mathbf{T}_{v}$, we immediately see that, in both cases, $G$ contains the image of a cocharacter with weights $\{1,0^3\}$. We then apply \cite[Proposition A.3]{Pink1997} to conclude that  $G=\GL(V_{v})$. As a result, we see that $G_p(A)^\circ\supseteq \mathbf{T}_p\cdot \mathrm{SU}(V_p/k_p)$.
 
 If we are in the (2,2) case, then we are already done: we know that $\mathfrak{hg}=\mathfrak{su}(V/k)$, which forces that $G_p(A)^\circ= \mathbf{T}_p\cdot \mathrm{SU}(V_p/k_p)$. 
 
 If we are in the (1,3) case, consider the following map 
\begin{equation}\label{eq:prod_maps}G_p(A)^\circ\xrightarrow{\delta\times \delta'} \mathbb{G}_m\times \mathbb{G}_m
\end{equation}
where $\delta:G_p(A)^\circ\hookrightarrow \GL(V_p)\xrightarrow{\det} \mathbb{G}_m$ and $\delta':G_p(A)^\circ\rightarrow \GL(V_{v})\xrightarrow{\det} \mathbb{G}_m$. It suffices to show that $\delta\times \delta'$ is surjective. Let's identify the cocharacter lattice of the target as $\bZ\oplus \bZ$. The image torus $(\delta\times \delta')(\mathbf{T}_p)$ has cocharacter lattice $(2,1)\mathbb{Z}$. On the other hand, the composition of the Newton cocharacter with $\delta\times \delta'$ is a torus with cocharacter lattice $(4,1)\mathbb{Z}$. This immediately shows that $\delta\times \delta'$ is surjective. Therefore $G_p(A)^\circ=\mathrm{GU}(V_p/k_p)$.
 \item Suppose that $\xi=\{0^2,(\frac{1}{2})^4, 1^2\}$ and  $\widetilde{A}$ has signature $(2,2)$. Recall that this happens when $p$ does not split in $k$, or when $p$ splits in $k$ and $\bH_{p,A}$ decomposes into two pieces both with slopes $\{0,(\frac{1}{2})^2,1\}$. We have $\mathfrak{h}_p\subseteq \mathfrak{hg}_{\bQ_p}=\mathfrak{su}(V_p/k_p)$. We wish to show that the above inclusion is an equality. Let's base change everything to $\bC_p$. Denote by $\tau$ and $\overline{\tau}$ the two maps $k_p\rightarrow\bC_p$. Then $V_p\otimes \bC_p=V_{\tau}\oplus V_{\overline{\tau}}$ and $\mathfrak{su}(V_p/k_p)\otimes \bC_p\simeq \mathfrak{sl}(V_\tau)$. Moreover, the representation $\mathfrak{su}(V_p/k_p)\otimes \bC_p\acts V_p\otimes \bC_p$ identifies with $$\mathfrak{sl}(V_\tau)\stackrel{\mathrm{std}\oplus \mathrm{std}^\vee}{\acts} V_\tau\oplus V_{\overline{\tau}}.$$
 where $\mathrm{std}$ is the standard representation. Now $\mathfrak{h}_p\otimes \bC_p$ identifies with a reductive subalgebra of $\mathfrak{sl}(V_\tau)$ such that $\mathfrak{h}_p\otimes \bC_p\acts V_{\tau}\otimes \bC_p$ is irreducible and not self-dual (for otherwise the representation $\mathrm{std}\oplus \mathrm{std}^\vee$ restricted to $\mathfrak{h}_p\otimes \bC_p$ would possess more endomorphisms, producing non-existing endomorphisms of $A_{\overline{L}}$). By \cite[\S 3.1]{Hui2023}, the representation $\mathfrak{h}_p\otimes \bC_p\acts V_{\tau}\otimes \bC_p$ is irreducible only in the following four cases: 
\begin{table}[H]
\centering
\renewcommand{\arraystretch}{1.35}
\begin{tabular}{c|c|c}
\toprule
Case &$\mathfrak{h}_p\otimes \bC_p$ & 
representation\\
\midrule
1&$\mathfrak{sl}(V_\tau)$ & std\\ 
2&$\mathfrak{sp}_4$& std\\
3&$\mathfrak{sl}_2\times \mathfrak{sl}_2$ & std $\otimes$ std \\
4&$\mathfrak{sl}_2$ & $\mathrm{Sym}^3$\\
\hline
\end{tabular}
\caption{}\label{tb:irreps}
\end{table}

Among them only case 1 is not self-dual. Therefore we must have $\mathfrak{h}_p=\mathfrak{su}(V_p/k_p)$, as desired. 
\item Suppose that $\xi=\{0^2,(\frac{1}{2})^4, 1^2\}$ and  $\widetilde{A}$ has signature $(1,3)$. This happens when $p$  splits as $v\overline{v}$ in $k$, and $\bH_{p,A}=\bH_{v,A}\oplus \bH_{\overline{v},A}$ such that $\bH_{v,A}$ has Newton polygon $\{0^2,(\frac{1}{2})^2\}$. We have $\mathfrak{h}_p\subseteq\mathfrak{u}(V_p/k_p)\simeq \mathfrak{gl}(V_{v})$, and we want to show that this is an equality. Again, $\mathfrak{u}(V_p/k_p)\acts V_p$ identifies with a direct sum of the standard rep $\mathfrak{gl}(V_{v})\acts V_v$ with its dual. 

One may naïvely think that our strategy in (c) carries over, but this is wrong: while in (c), $\mathfrak{h}_p$ is contained in $\mathfrak{sl}(V_\tau)$, in this case $\mathfrak{h}_p$ is only contained in $\mathfrak{gl}(V_\tau)$. This difference on the one dimensional center ruins the argument based on self-duality: as long as  $\mathfrak{h}_p$ contains the 1 dimensional homothety  $\mathfrak{t}_v\subseteq \mathfrak{gl}(V_v)$, its rep on $V_v$ can never be self-dual.  

Of course, the rep $\mathfrak{h_p}\otimes\bC_p\acts V_{v}\otimes \bC_p$ is still irreducible, therefore the semisimple part $\mathfrak{h}_p^{\mathrm{ss}}$ and its representation still falls into the four cases in Table~\ref{tb:irreps}. In the following, we will use the structure of $p$-adic monodromy to rule out cases 2,3 and 4. Let $\nu$ be the Newton cocharacter, which factors through $V_v$ as $\nu_v$, and the weight of $\nu_v$ is $\{0^2,1^2\}$. Let $U_{\nu_v}$ be the unipotent of $G$ with respect to $\nu$, where $G$ is the projection of $G_p(A)^\circ$ to $\GL(V_v)$. Then 
by Lemma~\ref{lm:keylemma}, the $U_{\nu_v}$ needs to have positive even dimension. 

Now compute the dimension of $U_{\nu_v}$ for case 2,3,4: For case 2, the dimension  of $U_{\nu_v}$ is 3. For case 3, $\nu_v$ need to factor non-centrally through exactly one copy of $\mathfrak{sl}_2$ (note that if it factors through both copies, then there will be 4 different slopes upon taking tensor rep). Therefore the dimension of $U_{\nu_v}$ is 1. For case 4, the dimension is 1. 

As a result, cases 2,3,4 cannot happen, so $\mathfrak{h}_{p}^{\mathrm{ss}}=\mathfrak{sl}(V_v)$. Finally, a map similar to (\ref{eq:prod_maps}) shows that $\mathfrak{h}_{p}=\mathfrak{gl}(V_v)$, as desired. 
\item Suppose that $\xi=\{(\frac{1}{4})^4,(\frac{3}{4})^4\}$. Then $p$ does not split in $k$, and $\widetilde{A}$ has signature $(2,2)$. The proof is very similar to (b): denote by $\tau$ and $\overline{\tau}$ the two maps $k_p\rightarrow\bC_p$. Then $V_p\otimes \bC_p=V_{\tau}\oplus V_{\overline{\tau}}$ and $\mathfrak{su}(V_p/k_p)\otimes \bC_p\simeq \mathfrak{sl}(V_\tau)$. The representation $\mathfrak{su}(V_p/k_p)\otimes \bC_p\acts V_p\otimes \bC_p$ is the sum of standard rep of $\mathfrak{sl}(V_\tau)$ and its dual and  $\mathfrak{h}_p\otimes \bC_p$ identifies with a reductive subalgebra of $\mathfrak{sl}(V_\tau)$ such that $\mathfrak{h}_p\otimes \bC_p\acts V_{\tau}\otimes \bC_p$ is irreducible and not self-dual. Again use the classification Table~\ref{tb:irreps} to conclude. 
\end{enumerate}
    
\end{proof}

\subsubsection{Type IV(2,1)}\label{subsub:Type IV(2,1)} \textit{The geometric endomorphism algebra of $A$ is a quartic CM field $E$. Then $\mathfrak{hg}\simeq \mathfrak{u}_{E}(V,\psi)$. Moreover,\begin{enumerate}
    \item If $E$ contains an imaginary quadratic subfield $k$ where $p$ splits (i.e. $A$ is of exotic type IV), then $\mathfrak{h_l}\simeq \mathfrak{u}_k\oplus \mathfrak{su}_{E_l}(V_l,\psi_l)$. The Tate conjecture holds for $A$. 
    \item Otherwise,  $\mathfrak{h_l}\simeq \mathfrak{u}_{E_l}(V_l,\psi_l)$ and the Tate ring of $A^n$ is generated by divisor classes.
\end{enumerate}}
\begin{proof}
     The claim for $\mathfrak{hg}$ is the same as \cite[\S 7.5. Type IV(2,1)]{MZ1995}. The subtle part is $\mathfrak{h}_l$ (or $\mathfrak{h}_p$), whose center would be denoted by $\mathfrak{c}_l$. Before we do a case by case analysis, we make some group theoretical observations.   First, we have $\mathfrak{h}_l\subseteq \mathfrak{u}_{E_l}(V_l,\psi_l)$. The algebra $\mathfrak{u}_{E_l}(V_l,\psi_l)\otimes \bC_l$ acts on $V_l\otimes \bC_l$ in the following way: $V_l\otimes \bC_l$ splits into four two-dimensional subspaces: \begin{equation}\label{eq:splitin4}
         V_l\otimes \bC_l=V_{l,1}\oplus V_{l,2}\oplus V_{l,3}\oplus V_{l,4},
     \end{equation} and 
$\mathfrak{u}_{E_l}(V_l,\psi_l)\otimes \bC_l\simeq \mathfrak{d}_1\times \mathfrak{d}_2$ such that $\mathfrak{d}_1\simeq \mathfrak{d}_2\simeq\mathfrak{gl}_2$  and $\mathfrak{d}_1$ acts on $V_{l,1}$ by standard rep, acts on $V_{l,2}$ by its dual, and acts trivially on $V_{l,3}\oplus V_{l,4}$; $\mathfrak{d}_2$ acts on $V_{l,3}$ by standard rep, acts on $V_{l,4}$ by its dual, and acts trivially on $V_{l,1}\oplus V_{l,2}$. Similar holds for the $p$-adic monodromy.

 The center of $\mathfrak{u}_{E_l}(V_l,\psi_l)$ is $\mathfrak{u}_{E_l}=\{a\in E_l:a+\overline{a}=0\}$, which has rank 2. %
 It is clear that $\mathfrak{c}_l\subseteq \mathfrak{u}_{E_l}$.

By Theorem~\ref{eq:Newton--endo}, the Newton polygons that allow type IV(1,1) endomorphism are 
    $$\xi=\{0^2,(\frac{1}{2})^4, 1^2\},\;\{(\frac{1}{4})^4, (\frac{3}{4})^4\}.$$
    \begin{enumerate}
        \item 
      Suppose that $E$ contains an imaginary quadratic subfield $k$ where $p$ splits. Then $\xi=\{0^2,(\frac{1}{2})^4, 1^2\}$ by Lemma~\ref{lm:exo_criterion}. From the same lemma we find that $E$ is Galois with three quadratic subfields $k$ $E_0$ and another imaginary $k'$, and $p$ does not split in $E_0$ and $k'$.  We will denote by $\mathfrak{p},\overline{\mathfrak{p}}$ the two places of $E$ dividing $p$.   
      
 To show that $\mathfrak{h}_l$ is as claimed, it suffices to show that $\mathfrak{c}_l=\mathfrak{u}_{k_l}$ (this forces the semisimple part of $\mathfrak{h}_l$ to be $\mathfrak{su}_{E_l}(V_l,\psi_l)$, for otherwise $A$ will have more endomorphisms). The proof is a modification of the argument in \cite[\S 7.5. Type IV(2,1)]{MZ1995} ($l$-adic case). Spread out $A$ to an abelian scheme over $X_0$ and let $x\in X_0(\Fpbar)$. Let $x_0$ be a finite field model of $x$. By \cite[Theorem 1.2.2]{D2020a}, we can assume that $x$ and $x_0$ are chosen so that $G_{l}(\bH_{l,{x_0}},x)= G_{l}(\bH_{l,{x_0}},x)^\circ$ is a maximal torus of $G_{l}(\bH_{l,{X_0}},x)$ (in particular, the multiplicative group generated by the eigenvalues of $\mathrm{Fr}_{x_0}$ does not contain nontrivial roots of unity). Let $\mathrm{Fr}_{x_0}$ be the Frobenius element corresponding to $x_0$ and let $q$ be the cardinality of $x_0$. Viewing $V_l$ as a two dimensional $E_l$-space  and $\mathrm{Fr}_{x_0}$ as an $E_l$-linear operator on $V_l$, we define \begin{equation}\label{eq:central_unit}
     u= \det_{E_l}(\mathrm{Fr}_{x_0}^2/q)\in E_l^*.
 \end{equation} Then $u$ actually lies in $E^*$ and is not a root of unity. By Weil conjecture, $u$ has norm 1 under every complex embedding and is a $\lambda$-adic unit for every finite place $\lambda$ of $E$ not dividing $p$. Let $\tau\in\Gal(E/\bQ)$ be the element whose fixing field is $k$, then $\tau$ fixes $\mathfrak{p}$ and $\overline{\mathfrak{p}}$. Therefore $\tau(u)/u$ has norm 1 under every complex embedding, and is a $\lambda$-adic unit for every finite place of $E$ (including those dividing $p$), hence a root of unity. This means that $\tau$ fixes some power of $u$. By our choice of $\mathrm{Fr}_{x_0}$ at the beginning, we find that $u\in k$. Therefore we have proved that $\bQ(u)=k$, hence $\mathfrak{c_l}\supseteq \mathfrak{u}_{k_l}$.    Let $H_l=(G_l(A)\cap\SL(V_l))^\circ$. We can view any element $g\in \SL(V_l)$ as a linear operator on the two dimensional $E_l$-space $V_l$. Consider the group homomorphism $$\delta:H_l(\bQ_l)\rightarrow E^*_l: \,\,g\rightarrow \det_{E_l}(g).$$
Let $T\subseteq H_l(\bQ_l)$ be the maximal torus generated by $\mathrm{Fr}_{x_0}^2/q$. Then $\delta(T)=k_l^*$. On the other hand, if $E^*_l\subseteq T$, then $\delta$ would map $E^*_l$ isomorphically onto its image. This contradiction prohibits the center of $H_l$ being $E^*_l$, establishing that 
$\mathfrak{c_l}=\mathfrak{u}_{k_l}$.

We will leave the proof of the Tate conjecture in this case to \ref{itemTate(IV, 2)}.

    \item Suppose that $E$ does not contain an imaginary quadratic subfield where $p$ splits. \begin{enumerate}
        \item Case $\xi=\{0^2,(\frac{1}{2})^4, 1^2\}$. We claim that $\mathfrak{c}_l=\mathfrak{u}_{E_l}$. It is easy to see that $\mathfrak{c}_l\neq \mathbf{0}$. 
    Construct $u$ as in (\ref{eq:central_unit}). If $\dim\mathfrak{c}_l=1$, then arguing as in \cite[\S 7.5. Type IV(2,1)]{MZ1995} ($l$-adic case), we find that $k:=\bQ(u)$ is an imaginary  quadratic subfield of $E$ and $u\overline{u}=1$. By our assumption, $p$ does not split in $k$. So there is only one place of $k$, say $\mathfrak{P}$, dividing $p$. Then $0=v_{\mathfrak{P}}(u)+v_{\mathfrak{P}}(\overline{u})=2v_{\mathfrak{P}}(u)$. Since $u$ is a $\lambda$-unit for all places of $k$ no dividing $p$, and has norm 1 under any complex embedding, we see that $u$ is a root of unity. Contradiction.

It remains to show that the semisimple part $\mathfrak{h}_l^{\mathrm{ss}}=\mathfrak{su}_{E_l}(V_l,\psi_l)$. We will work in the $p$-adic setting: i.e., in the case where $l=p$. We spread out $A$ to an abelian scheme $A_{X_0}$ and let $x$ be a suitable base point. Let the image of the Newton cocharacter in $\mathfrak{h}_p$ be $\nu'$, which has weight $\{(-1)^2,1^2,0^4\}$ (up to scaling). Up to permutation, the distribution of $\nu'$-weights on the spaces $V_{p,i}$ in (\ref{eq:splitin4}) has three possibilities: 
        \begin{table}[H]
\centering
\renewcommand{\arraystretch}{1.35}
\begin{tabular}{c|cccc}
\toprule
&$V_{p,1}$&$V_{p,2}$&$V_{p,3}$&$V_{p,4}$\\
\midrule
A&$\{0,1\}$&$\{0,-1\}$&$\{0,1\}$&$\{0,-1\}$\\
B&$\{-1,1\}$&$\{-1,1\}$&$\{0,0\}$&$\{0,0\}$\\
C&$\{-1,-1\}$&$\{1,1\}$&$\{0,0\}$&$\{0,0\}$\\
\hline
\end{tabular}
\end{table}
Suppose that $\mathfrak{h}^{\mathrm{ss}}_p\otimes \bC_p\subsetneq \mathfrak{sl}_2\times \mathfrak{sl}_2$. Since the projection of $\mathfrak{h}^{\mathrm{ss}}_p\otimes \bC_p$ to each factor need to be $\mathfrak{sl}_2$ (otherwise you get non-existing endomorphisms of $A$),  $\mathfrak{h}^{\mathrm{ss}}_p\otimes \bC_p\simeq \mathfrak{sl}_2$ but embedded diagonally in $\mathfrak{sl}_2\times \mathfrak{sl}_2$. It is then easy to see that case B and C do not happen. So we are in case A. But then there is a nontrivial extension between the slope 0 part and the slope $\frac{1}{2}$ part of $\bH_{p,X_0}$. Using the computation in \S\ref{example:0-1/2-1} we find that the unipotent of $\mathfrak{h}_p$ with respect to $\nu'$ is at least of 2 dimensional. Contradiction.  

\item Case $\xi=\{(\frac{1}{4})^4, (\frac{3}{4})^4\}$. The same argument in (2) shows that $\mathfrak{c}_l=\mathfrak{u}_{E_l}$. By Lemma~\ref{lm:evendim two slope}, the Newton cocharacter $\nu$ should have positive even dimensional unipotent. This immediately implies that $\mathfrak{h}^{\mathrm{ss}}_p=\mathfrak{su}_{E_p}(V_p,\psi_p)$.
  \end{enumerate}
\item\label{itemTate(IV, 2)} We prove the Tate conjecture for $A$ in (1). This is done by performing a lifting to characteristic 0 and specializing Markman's result. Recall that $E$ contains a real quadratic subfield $E_0$ and two imaginary subfields $k$ and $k'$ such that $p$ splits in $k$ but not in $k'$. %
Since the action $\mathfrak{h}_l\acts V_{l}$ is made explicit at the beginning following (\ref{eq:splitin4}), it is not hard to compute that $$\mathcal{T}_l^2(A)=\mathscr{D}_l^2(A)\oplus \mathscr{W}_l(A),\qquad\mathscr{W}_l(A)=\mathscr{W}_{k',l}(A)=(\bigwedge_{k'\otimes\bQ_l}^4 V_l) (2),$$
where $\dim \mathscr{D}_l^2(A)=3 $ and $ \dim\mathscr{W}_l(A)=2$. More explicitly, $\mathscr{W}_l(A)$ is generated by the two $\mathfrak{h}_l$-invariant $k'$-lines: $$l_1:=\bigwedge^2V_{l,1}\wedge \bigwedge^2V_{l,4},\qquad l_2:=\bigwedge^2V_{l,2}\wedge \bigwedge^2V_{l,3}.$$
One checks that $\mathscr{W}_l(A):=l_1\oplus l_2$ descends to $\bQ_l$. On the other hand, the eigenspaces of  $k'\acts V_l$ are $k'$-spaces $V_{l,1}\oplus V_{l,4}$ and $V_{l,2}\oplus V_{l,3}$. Then $l_1 $ (resp. $l_2$) is the top wedge of $V_{l,1}\oplus V_{l,4}$ (resp. $V_{l,2}\oplus V_{l,3}$). This gives the desired formula for $\mathscr{W}_l(A)$
being the top wedge of $V_l$ as a $k'\otimes \bQ_l$-space. 

By Lemma~\ref{lm:exo_criterion}(3), there exists a lift $\widetilde{A}$ up to isogeny, together with $k'$-action of signature $(2,2)$. %
By the proof of \cite[Corollary 1.6.1]{Markman2025}, the two dimensional space of Hodge--Weil classes $$\mathscr{W}_{k',\mathrm{B}}(\widetilde{A}_{\bC})=(\bigwedge_{k'}^4 H^1_\mathrm{B}(\widetilde{A}_{\bC},\bQ))(2)$$
is spanned by algebraic cycles. By specializing, we see that $\mathscr{W}_l(A)$ is spanned by algebraic cycles. This proves the Tate conjecture for $A$. 
    \end{enumerate}
\end{proof}
\subsection{Alternative proof of the Tate conjecture for Type IV} We have seen that Type IV cases are the most involved. If our goal is only to prove the Tate conjecture in these cases, rather than the Mumford--Tate conjecture, then the proof can be simplified to remove any need for the signature-lifting statements in Theorem~\ref{thm:lift_endo_fourfold}.
The key is the following lemma adapted from \cite{Broe2026}. We will also later use the lemma to prove the Tate conjecture for certain non-simple abelian fourfolds (see step \ref{it:TT5} of the proof of Theorem \ref{thm:Tate_nonsimple}). Let $A$ be an abelian fourfold over a function field $L$ of positive characteristic $p$.

\begin{lemma}\label{lemWeilAlg13}
     Let $k \subseteq \End^0(A)$ be an imaginary quadratic field, and let $\mathscr{W}_l(A) := \mathscr{W}_{k, l}(A)$ denote the associated Weil classes. Lift $A$ up to isogeny, together with the $k$-action, to an abelian fourfold $\widetilde{A}$ in characteristic zero. Suppose that the $k$-action on $\widetilde{A}$ has signature $(1, 3)$. The following are equivalent: 
    \begin{enumerate}
        \item $\mathscr{W}_l(A)(2)$ consists of Tate classes;
        \item $\mathscr{W}_l(A)(2)$ consists of  algebraic classes;
        \item the prime $p$ does not split in $k$.
    \end{enumerate}
\end{lemma}
\begin{proof}
    We follow the argument of \cite[Section 7]{Broe2026}. Assume (1) holds, so that $\mathscr{W}_l(A)(2) \subseteq \mathscr{T}^2_l(A)$. One then sees by specializing $A$ to an abelian fourfold over a finite field, and arguing as in \cite[Lemma 7.1]{Broe2026}, that $p$ does not split in $k$. By Deuring's results \cite[Theorem 1.1]{Blake2014}, this implies that there exists an elliptic curve $\widetilde{E}$ over $\overline{\bQ}_p \cong \bC$ with CM by $k$ and supersingular reduction. 
    \par
    Fix an isomorphism $\tau: k \cong \End^0(\widetilde{E})$. Replacing $\tau$ with its conjugate if necessary, the diagonal action of $k$ on $\widetilde{Y} = \widetilde{A}_\bC \times_{\bC} \widetilde{E}^2$ has signature (3,3), so that $(\widetilde{Y}, k)$ is an abelian sixfold of Weil type. There further exists a polarization of $(\widetilde{Y}, k)$ of discriminant $-1$ \cite[Proposition 1.6]{Milne2022}. Then \cite[Theorem 1.5.1]{Markman2025} implies that the Hodge--Weil classes on $\widetilde{Y}$ are algebraic. We specialize to obtain algebraic classes on the reduction $Y \sim A \times_L E^2$. Using the supersingularity of $E^2$, we construct a suitable algebraic correspondence between $Y$ and $A$, and find that $\mathscr{W}_l(A)$ is spanned by algebraic cycles, exactly as in \cite{Broe2026}. 
    \par
    We have shown $(1) \implies (3) \implies (2)$, and $(2) \implies (1)$ is trivial, so the proof is complete.
\end{proof}

Now assume that $A$ is Type IV, with endomorphism algebra $E$. We explain how to use the above lemma to give an alternative proof of the Tate conjecture for $A$, which avoids the use of Theorem \ref{thm:lift_endo_fourfold}. By Proposition \ref{prop:lifting_field}, we may fix a lift $\widetilde{A}$ of $A$ to characteristic zero, such that $\End^0(\widetilde{A}) = E$.

\subsubsection{Type IV(1,1)}
If $p$ splits in $E$, then the earlier proof of the Tate conjecture for $A$ given in \ref{subsub:TypeIV(1,1)} has no dependence on Theorem \ref{thm:lift_endo_fourfold}: here the weaker Proposition \ref{prop:lifting_field} suffices for all uses of the theorem.
\par
We thus only treat the case where $p$ does not split in $E$. If the $E$-action on $\widetilde{A}$ has signature (2, 2), then the argument of part (b) of the proof of \ref{subsub:TypeIV(1,1)} yields the Tate conjecture for $A$. Otherwise the $E$-action on $\widetilde{A}$ has signature (1, 3). Then Lemma \ref{lemWeilAlg13} implies that $\mathscr{W}_{E, l}(A)(2)$ is spanned by algebraic classes, so $\mathscr{W}_{E, p}(A)(2)$ is also spanned by algebraic classes. This in turn implies that $\mathfrak{h}_p \subseteq \mathfrak{su}(V_p/k_p)$. We may again argue as in part (b) of the proof of \ref{subsub:TypeIV(1,1)}, and conclude.
    \par
    Contrast the above argument in the case where $p$ is nonsplit with that of \ref{subsub:TypeIV(1,1)}. The latter argument uses Theorem \ref{thm:lift_endo_fourfold} to reduce to the case where $\widetilde{A}$ has signature (2, 2).
    
\subsubsection{Type IV(2,1)}
If $E$ does not contain an imaginary quadratic subfield $k$ where $p$ splits, then we argue as in part (2) of the proof of \ref{subsub:Type IV(2,1)} and conclude. 
\par
We now assume that $E$ does contain such a $k$. The first part of the proof of Proposition \ref{lifting mom} (which does not depend on Theorem \ref{thm:lift_endo_fourfold}) then shows that $p$ is not totally split in $E$. Following part (1) of the proof of \ref{subsub:Type IV(2,1)}, we obtain an imaginary quadratic subfield $k'$ of $E$ in which $p$ is not split, and reduce to showing that $\mathscr{W}_{k', l}(A)$ is spanned by algebraic classes. If the $k'$-action on $\widetilde{A}$ has signature (2, 2), then we conclude as in part (3) of the proof of \ref{subsub:Type IV(2,1)}. Otherwise $\widetilde{A}$ has signature (1, 3), and we conclude by Lemma \ref{lemWeilAlg13}.

\section{The Tate conjecture for non-simple fourfolds}
 Let $A$ be a non-CM abelian variety of dimension $g\leq 4$ over a characteristic $p$ function field $L$ whose simple factors have dimension $\leq 3$. The main results of this section is 
\begin{theorem}\label{thm:Tate_nonsimple}
 The Tate conjecture holds for $A$. 
\end{theorem}
\begin{theorem}\label{thm:Mumford_Tate_nonsimple}
 The weak Mumford--Tate conjecture holds for $A$. If all simple factors of $A$ are non-CM and not of exotic type IV, then the Mumford--Tate conjecture holds for $A$.
\end{theorem}

\subsection{Simple abelian varieties of dimension $\leq 3$} 
Let $A$ be a simple non-CM abelian variety of dimension $g\leq 3$ over a characteristic $p$ algebraically closed field $k$. Its endomorphism algebra $D$ is well understood (cf. \cite[Proposition 4.1, Remark 6]{Hwang2022}): \begin{table}[H]
\centering
\renewcommand{\arraystretch}{1.35}
\begin{tabular}{|c|c|c|}
\toprule
$g$& {Type}   &$D$ \\
\midrule
1& I(1) & $\bQ$ \\\cline{1-3}
& I(1) & $\bQ$ \\\cline{2-2}
2& I(2) & Real quadratic field \\\cline{2-2}
 & II(1) & Indefinite quaternion over $\bQ$ \\
   \hline 
   & I(1) & $\bQ$ \\\cline{2-2}
3& I(3) & Totally real cubic field \\\cline{2-2}
 & IV(1,1) & Imaginary quadratic field \\
   \hline 
\end{tabular}
\end{table}
The following table lists possible types of $D$ associated to a given Newton polygon $\xi$. The proof is an easy version of Theorem~\ref{eq:Newton--endo}. We additionally use Theorem \ref{thm:nonvan} to exclude the case where $A$ is type IV(1,1) with Newton polygon $\{(\frac{1}{3})^3,(\frac{2}{3})^3\}$, and $p$ is split in $\End^0(A)$. We omit the proof.

\begin{table}[H]
\centering
\renewcommand{\arraystretch}{1.35}
\begin{tabular}{|c|c|c|}
\toprule
$g$&$\xi$& $D$  \\
\midrule
1&$\{0,1\}$ & I(1)\\\hline
2&$\{0^2,1^2\}$ & I(1), I(2), II(1)\\
&$\{0,(\frac{1}{2})^2,1\}$ & I(1), I(2)\\  \hline
&$\{0^3,1^3\}$& I(1), I(3), IV(1,1)\\ 
3&$\{0^2,(\frac{1}{2})^2,1^2\}$& I(1), I(3), IV(1,1)\\
&$\{0,(\frac{1}{2})^4,1\}$& I(1), I(3), IV(1,1)\\ 
&$\{(\frac{1}{3})^3,(\frac{2}{3})^3\}$& I(1), I(3)\\\hline
\end{tabular}
\caption{Type of $D$ vs. Newton polygon}\label{tb:Endo2}
\end{table}
The following is a lifting result in the the spirit of Theorem~\ref{thm:lift_endo_fourfold}:
\begin{theorem}\label{thm:lift3fold}
     Up to isogeny, there is a lift $\widetilde{A}$ over $\overline{W}$ such that  $\End^0(\widetilde{A})=D$. 
\end{theorem}
\begin{proof}
    All endomorphism algebras are fields except possibly ordinary abelian surfaces, so we just lift by Proposition~\ref{prop:lifting_field}.
For ordinary abelian surface we do the canonical lift.  \end{proof}

Mimicking Notation~\ref{not:compareMZ}, we define $\mathfrak{g},\mathfrak{hg}$ by $\widetilde{A}_{\bC}$, and 
$\mathfrak{g}_l,\mathfrak{h}_l$, $\mathfrak{g}_p,\mathfrak{h}_p$ in a similar manner. In the following, we compute the monodromy in each case. %
\begin{theorem}\label{thm:simple_leq 3}
Let $A$ be a simple non-CM $g\leq 3$ dimensional abelian variety over a finitely generated field over $\mathbb{F}_p$. If $A$ is not an abelian threefold of exotic type IV, then $A$ satisfies the Mumford--Tate conjecture, and for all $n\geq 1$ the Tate ring of $A^n$ is generated by divisor classes. If $A$ is an abelian threefold of exotic type IV, then $A$ satisfies the weak Mumford--Tate conjecture, and the Tate ring of $A$ is generated by divisor classes. The following chart summarizes the monodromy in each case:  
    \begin{table}[H]\label{tb:monLe3}
\centering
\renewcommand{\arraystretch}{1.35}
\begin{tabular}{|c|c|c|c|c|c|}
\toprule
$g$& $D$ & $\mathrm{End}_{D}(\mathfrak{sp}(V))$ &  $\mathrm{End}_{D}(\mathfrak{sp}(V))\otimes \bC$ &$\mathfrak{hg}$  & $\mathfrak{h}_l$ \\
\midrule
1&I(1) & $\mathfrak{sl}(V)$ &$\mathfrak{sl}_2$ & $\mathfrak{sl}(V)$ & $\mathfrak{sl}(V_l)$ \\\hline
&I(1) &$\mathfrak{sp}(V)$ & $\mathfrak{sp}_4$ & $\mathfrak{sp}(V)$ & $\mathfrak{sp}(V_l)$ \\
2&I(2) &$\mathfrak{sp}_D(V)$&$\mathfrak{sl}_2^2$ & $\mathfrak{sp}_D(V)$ & $\mathfrak{sp}_D(V_l)$ \\
&II(1) & $D^{\mathrm{Nm}=1}$ & $\mathfrak{sl}_2$&$D^{\mathrm{Nm}=1}$ & $D_{\bQ_l}^{\mathrm{Nm}=1}$ \\\hline
&I(1) & $\mathfrak{sp}(V)$ & $\mathfrak{sp}_6$& $\mathfrak{sp}(V)$ & $\mathfrak{sp}(V_l)$ \\
3&I(3) &$\mathfrak{sp}_D(V)$ &$\mathfrak{sl}_2^3$&  $\mathfrak{sp}_D(V)$ & $\mathfrak{sp}_D(V_l)$ \\
&IV(1,1)$^\mathrm{n}$& $\mathfrak{u}(V/D)$ & $\mathfrak{gl}_3$& $\mathfrak{u}(V/D)$ & $\mathfrak{u}(V_l/D_l)$ \\
&IV(1,1)$^\mathrm{e}$ & $\mathfrak{u}(V/D)$ & $\mathfrak{gl}_3$& $\mathfrak{u}(V/D)$ & $\mathfrak{su}(V_l/D_l)$ \\\hline
\end{tabular}
\end{table}
The symbol IV(1,1)$^\mathrm{e}$ (resp. IV(1,1)$^\mathrm{n}$) refers to exotic type IV (resp. non-exotic type IV). 
\end{theorem}
\begin{proof} The proof is quite involved. So we proceed in three steps: 
\begin{enumerate}
    \item 
The easier cases are: $g=1$, $g=2$ type I(2), II(1), and $g=3$ type I(3). In these cases, if $\mathfrak{hg}$ (resp. $\mathfrak{h}_l$) was smaller than the centralizer, then you would get much larger $D$.
\item 
Now for the type I(1) cases when $g=2,3$: First, the  fact that $\mathfrak{hg}$ is the full symplectic algebra is well-known. We then compute $\mathfrak{h}_l$ by classifying all reductive subalgebras of $\mathfrak{sp}_4$ or $\mathfrak{sp}_6$ whose actions on $V_l$ are still irreducible (these are classified in \cite[Proposition 4.5]{Gee2013} and its proof): For $g=2$ we may have (geometrically) $\mathfrak{h}_l=\mathfrak{sl}_2$ acting by $\mathrm{Sym}^{3}$. For $g=3$ we may have (geometrically) $\mathfrak{h}_l=\mathfrak{sl}_2$ acting by $\mathrm{Sym}^5$, or $\mathfrak{h}_l=\mathfrak{sl}_2\times \mathfrak{sl}_2$ acting by $\mathrm{std}\otimes \mathrm{Sym}^2$.

The $\mathrm{Sym}^\bullet$ reps of $\mathfrak{sl}_2$ for both $g=2$ and $g=3$ do not happen for the same reason we rule out 
$\mathrm{Sym}^7$ in Type I(1) (\S\ref{subsub:TypeI(1)}): for otherwise the Newton polygon has too many slopes.

For the $\mathrm{std}\otimes \mathrm{Sym}^2$ representation when $g=3$, we again look at $p$-adic monodromy $\mathfrak{h}_p$. Geometrically $\mathfrak{h}_p=\mathfrak{sl}_2\times \mathfrak{sl}_2$ acting by $\mathrm{std}\otimes \mathrm{Sym}^2$. Then the Newton cocharacter cannot factor nontrivially through both copies of $\mathfrak{sl}_2$, for otherwise there are too many different Newton slopes. So we are put into two subcases: \begin{enumerate}
    \item If the Newton cocharacter factors nontrivially through the first copy of $\mathfrak{sl}_2$, then the Newton polygon $\xi$ of $A$ has only two slopes. So $\xi=\{0^3,1^3\}$ or $\{(\frac{1}{3})^3,(\frac{2}{3})^3\}$. Note that $A$ can not be ordinary, since the Mumford--Tate conjecture holds for ordinary abelian threefolds by \cite{Jiang2025}. On the other hand, if $\xi=\{(\frac{1}{3})^3,(\frac{2}{3})^3\}$, then the unipotent of $\mathfrak{h}_p$ with respect to the Newton cocharacter need to be at least 3 dimensional by Lemma~\ref{lm:3slope}, which can not happen.  
    \item If the Newton cocharacter factors nontrivially through the second copy of $\mathfrak{sl}_2$. Then the Newton polygon $\xi$ of $A$ has three slopes and each appears with multiplicity 2. The only case that happens $\xi=\{0^2,(\frac{1}{2})^2,1^2\}$. Let $\Lie U_{\nu}\subseteq \mathfrak{h}_p$ be the unipotent of the Newton cocharacter, then it is of dimension 1. On the other hand, the very explicit $\mathrm{std}\otimes \mathrm{Sym}^2$-rep lets us explicitly compute the action $\Lie U_{\nu}\acts V_p$, and we find that there is non-trivial extension between the slope 0 and the slope $\frac{1}{2}$ part of the $F$-isocrystal associated to (a suitable spreadout) of $A$. The Frobenius structure on $\Lie U_{\nu}$ then forces $\Lie U_{\nu}$ to contain a $\sigma$-isocrystal of slope $1/2$, hence of dimension at least 2 (similar to the computation we made in \S\ref{example:0-1/2-1}). This contradicts the fact that $\dim\Lie U_{\nu}=1$.
\end{enumerate}
\item 
What remains is abelian threefolds of type IV(1,1). For $\mathfrak{hg}$, this follows from \cite[Theorem 3]{Ribet1983} (note that the signature of the complex multiplication needs to be $(2,1)$). Now we determine the $l$-adic (or $p$-adic) monodromy.   

In the IV(1,1)$^{\mathrm{e}}$ case, $p$ cannot split in $D$: otherwise $V_p$ would split into two mutually dual direct summands, of respective slopes $\{0^2, 1/2\}$ and $\{1/2, 1^2\}$, which is impossible. Define $u=\det_{D_l}(\mathrm{Fr}_{x_0}^2/q)\in E^*$ as in (\ref{eq:central_unit}), which is of norm 1 under all complex embeddings and is a $\lambda$-unit for every finite place $\lambda\nmid p$. Let $\tau\in \mathrm{Gal}(D/\bQ)$ be the complex conjugation. Since $p$ is not split in $D$,  $\tau$ fixes the unique prime above $p$. So $\tau(u)/u$ is a root of unity. This means that a power of $u$ lies in $ \bQ^*$ hence must be $\pm 1$. Therefore $u$ is a root of unity. This means that the maximal torus of $H_l:=(G_l(A)\cap \SL(V_l))^\circ$ generated by (a suitable power of)  $\mathrm{Fr}_{x_0}^2/q$ lies in $\mathrm{SU}(V_l/D_l)$, so $\mathfrak{h}_l\subseteq \mathfrak{su}(V_l/D_l)$. But then we must have $\mathfrak{h}_l= \mathfrak{su}(V_l/D_l)$ since otherwise the monodromy is too small and you get non-existing endomorphism.

For the IV(1,1)$^{\mathrm{n}}$ case, the Newton polygons that can appear are $\{0^3,1^3\}$ and $\{0,(\frac{1}{2})^4,1\}$.  We work with $\mathfrak{h}_p$. Base change everything to $\bC_p$. Denote by $\tau$ and $\overline{\tau}$ the two maps $D_p\rightarrow\bC_p$. Then $V_p\otimes \bC_p=V_{\tau}\oplus V_{\overline{\tau}}$ and 
$\mathfrak{u}(V_p/D_p)\otimes \bC_p\simeq \mathfrak{gl}(V_\tau)\simeq \mathfrak{gl}(V_{\overline{\tau}})$. Moreover, the representation $\mathfrak{u}(V_p/D_p)\otimes \bC_p\acts V_p\otimes \bC_p$ identifies with 
$$\mathfrak{gl}(V_\tau)\stackrel{\mathrm{std}\oplus \mathrm{std}^\vee}{\acts} V_\tau\oplus V_{\overline{\tau}}.$$
Let $\nu$ be the Newton cocharacter, which factors as $\nu'$ through $\mathfrak{u}(V_p/D_p)$. Then the above duality implies that $\nu'$ acts on $V_\tau$ and $V_{\overline{\tau}}$ with opposite weights. On the other hand, by Theorem~\ref{thm:nonvan}, $\nu'$ cannot act on 
$V_\tau$ via pure weights. It is then easy to check that for both Newton polygons, $\nu'$ acts on the two parts by weights $(-1)^2,1$ and $-1,1^2$ (up to scaling).

If we denote the projection of the $p$-adic monodromy group $G_p(A)^\circ$ to $ \GL(V_\tau)$ by $G$, then $G$ contains the torus of homothety in $ \GL(V_\tau)$, and the above weight analysis implies that $G$ contains the image of a cocharacter of weights $\{0^2,1\}$. We then apply \cite[Proposition A.3]{Pink1997} to conclude that  $G=\GL(V_{\tau})$. As a result, $\mathfrak{h}_p\supseteq \mathfrak{su}(V_p/D_p)$. If $\mathfrak{h}_p=\mathfrak{su}(V_p/D_p)$, then the weights of $\nu'$ on $V_{\tau}$ need to sum up to 0. This is clearly not true. Therefore $\mathfrak{h}_p=\mathfrak{u}(V_p/D_p)$. 
\end{enumerate}
\end{proof}

\subsection{The Tate conjecture}
Let $A$ be a non-CM non-simple abelian fourfold over finitely generated field $L/\mathbb{F}_p$. We will always assume that $L$ is sufficiently large so that $G_l(A)$ is connected.%
To prove the Tate conjecture, we need to analyze $H^4_{\text{ét}}(A_{\overline{L}},\bQ_l)$. If $A\sim A_1\times A_2$, then the above cohomology decomposes by Kunneth's formula. Let $$V_l(A)_{i,j}:=H^i_{\text{ét}}(A_{1,\overline{L}},\bQ_l)\otimes H^j_{\text{ét}}(A_{2,\overline{L}},\bQ_l)=\bigwedge^iV_l(A_1)\otimes \bigwedge^jV_l(A_2).$$ We will use the group $H_l(A)$ (or its Lie algebra $\mathfrak{h}_l(A)$) instead of $G_l(A)$ for the representation in order to avoid dealing with Tate twists. Two cases can happen:\begin{enumerate}
    \item $\dim A_1=1$. Then for all $i+j=4$ the fixed subspace $(V_l(A)_{i,j})^{\mathfrak{h}_l(A)}$ is spanned by divisor classes except possibly when $i=1$. 
    \item $\dim A_1=2$. Then for all $i+j=4$ the fixed subspace $(V_l(A)_{i,j})^{\mathfrak{h}_l(A)}$ is spanned by divisor classes except possibly when $i=2$ (note that when $i=1$, we have duality $H^3_{\text{ét}}(A_{2,\overline{L}},\bQ_l)\simeq H^1_{\text{ét}}(A_{2,\overline{L}},\bQ_l)^\vee$ and it reduces to endomorphism between $A_1,A_2$. Similar for $i=3$). In this case, the $\mathfrak{h}_l(A_1)$-module 
$V_l(A_1)_{2,2}$ splits as $V_l(A_1)_{2,2}^{\mathfrak{h}_l(A_1)}\oplus U_1$, and likewise the $\mathfrak{h}_l(A_2)$-module 
$V_l(A_2)_{2,2}$ splits as $V_l(A_2)_{2,2}^{\mathfrak{h}_l(A_2)}\oplus U_2$. The possible exceptional classes in $(V_l(A)_{2,2})^{\mathfrak{h}_l(A)}$ all lie in 
$(U_1\otimes U_2)^{\mathfrak{h}_l(A)}$.

\end{enumerate}

\begin{proof}[Proof of Theorem~\ref{thm:Tate_nonsimple}] We will use the notation and  observation as above. We will only deal with the case of abelian fourfolds. The $g\leq 3$ cases are much easier and are omitted.
\begin{enumerate}
\item\label{it:TT1} Suppose that $A$ is isogenous to the product of four elliptic curves. The case is easy and we leave it to the reader. 
 \item\label{it:TT2} 
 Suppose that $A\sim A_1\times A_2$ and $A_1$ is a simple non-CM abelian surface. If $A_1$ is isogenous to $A_2$, then the Tate conjecture holds for $A$ by Theorem~\ref{thm:simple_leq 3}. Otherwise, we show that $(V_l(A)_{2,2})^{\mathfrak{h}_l(A)}$ is spanned by divisor classes, i.e., showing that $(U_1\otimes U_2)^{\mathfrak{h}_l(A)}=0$. In fact, by the classification of monodromy in the Theorem~\ref{thm:simple_leq 3} we can explicitly compute the representation 
 $\mathfrak{h}_l(A_1)\acts U_1$ in the three cases I(1), II(2) and II(1), and show that it is absolutely irreducible. We also note that in all three cases, $\mathfrak{h}_l(A_1)$ is semi-simple (and simple if we choose a nice $l$). 

If $(U_1\otimes U_2)^{\mathfrak{h}_l(A)}\neq 0$ then by Goursat's lemma $\mathfrak{h}_l(A)=\mathfrak{h}_l(A_1)\times_\mathfrak{h} \mathfrak{h}_l(A_2)$, where $\mathfrak{h}$ is a nontrivial common quotient. The semi-simplicity of $\mathfrak{h}_l(A_1)$ would then imply the existence of a nonzero homomorphism between $A_1$ and $A_2$, and we run into a contradiction. Therefore $(U_1\otimes U_2)^{\mathfrak{h}_l(A)}=0$ as desired. 
 \item\label{it:TT3}  Suppose that $A\sim A_1\times A_2$ and $A_1\sim E_1\times E_2$ is isogenous to a product of two non-CM elliptic curves and $A_2$ is a simple CM abelian surface. Again, we want to show that $(V_l(A)_{2,2})^{\mathfrak{h}_l(A)}$ is spanned by divisor classes. Note that $H^2_{\text{ét}}(A_{1,\overline{L}},\bQ_l)$ further decomposes by Kunneth, and $U_1=V_l(E_1)\otimes V_l(E_2)$ is 4-dimensional when $E_1\not\sim E_2$; and $U_1$ is a 3-dimensional direct summand of $\End(E_1,E_2)\otimes \bQ_l\simeq (V_l(E_1)\otimes V_l(E_2))^{\mathfrak{h}_l(A_1)}$ when $E_1\sim E_2$ (in fact, it is the $\mathrm{Sym}^2$-rep of $\mathfrak{sl}_2$). 

 If $(U_1\otimes U_2)^{\mathfrak{h}_l(A)}\neq 0$, we again use a Goursat's lemma argument as above to derive a contradiction (it is even easier: Goursat implies that $\mathfrak{h}_l(A_1)$ has a non-trivial commutative quotient, which cannot happen since $\mathfrak{h}_l(A_1)$ is either $\mathfrak{sl}_2$ or $\mathfrak{sl}_2^2$). 

 \item\label{it:TT4}  Suppose that $A\sim A_1\times A_2$,  $A_1$ is a non-CM elliptic curve and $A_2$ is CM.  As we have seen above, the possible exceptional classes lie in  $(V_l(A)_{1,3})^{\mathfrak{h}_l(A)}$. We claim that $(V_l(A)_{1,3})^{\mathfrak{h}_l(A)}=0$. Again, if this is not 0, then Goursat's lemma implies that $\mathfrak{h}_l(A_1)$ and $\mathfrak{h}_l(A_2)$ have a nontrivial common quotient, which does not happen. 
 \item\label{it:TT5}  Suppose that $A\sim A_1\times A_2$, where $A_1$ is a non-CM simple threefold. The possible exceptional classes lie in  $(V_l(A)_{3,1})^{\mathfrak{h}_l(A)}$. 
 \par
 If $A_2$ is not CM, then we claim that $(V_l(A)_{3,1})^{\mathfrak{h}_l(A)}=0$. This follows, again, from the classification of monodromy in Theorem~\ref{thm:simple_leq 3} together with Goursat's lemma: when $A_2$ is not CM, then $(V_l(A)_{3,1})^{\mathfrak{h}_l(A)}\neq 0$ implies that $\mathfrak{h}_l(A)$ is a graph of $\mathfrak{h}_l(A_1)\twoheadrightarrow \mathfrak{h}_l(A_2)=\mathfrak{sl}_2$. This cannot happen, except possibly in the case where $A_1$ is type I(3) and the natural map of $\bQ_l$-algebras $\bQ_l \to \End^0(A_1) \otimes_{\bQ} \bQ_l$ admits a section. But the ``$l$-independence'' of $G_l(A)^\circ$ (cf. Definition \ref{def:mon}) implies that if the invariants $(V_l(A)_{3,1})^{\mathfrak{h}_l(A)}$ vanish for a single $l \neq p$, then the same is true for all $l \neq p$. We thus reduce to the case where no such section exists.
 \par
 If $A_2$ is supersingular, then $\mathfrak{h}_l(A_2)=0$, and $(V_l(A)_{3,1})^{\mathfrak{h}_l(A)}\neq 0$ implies that $\mathfrak{h}_l(A_1)\acts\bigwedge^3_{\bQ_l} V_l(A_1)$ has a nonzero fixed vector. If $A_1$ is not exotic type IV, then such a possibility is excluded by the classification of monodromy. If $A_1$ is exotic type IV with $k := \End^0(A_1)$, then
 \begin{gather*}
     \mathfrak{h}_l(A) = \mathfrak{h}_l(A_1) = \mathfrak{su}(V_l(A_1) / k_l),
 \end{gather*}
 and
 \begin{align*}
     (V_l(A)_{3,1})^{\mathfrak{h}_l(A)} &= \mathscr{W}_{k, l}(A_1) \otimes_{\bQ_l} V_l(A_2).
 \end{align*}
 We may write
 \begin{gather*}
     \mathscr{W}_{k, l}(A_1) \otimes_{\bQ_l} V_l(A_2) = (\mathscr{W}_{k, l}(A_1) \otimes_{k_l, \iota_l} V_l(A_2)) \oplus (\mathscr{W}_{k, l}(A_1) \otimes_{k_l, \overline{\iota}_l} V_l(A_2))
 \end{gather*}
 where $\iota: k \hookrightarrow \End^0(A_2)$ is a fixed embedding, and $\iota_l$ denotes the induced $k_l$-module structure on $V_l$. We note that
 \begin{gather*}
     \mathscr{W}_{k, l}(A_1) \otimes_{k_l, \iota_l} V_l(A_2) = \mathscr{W}_{k, l}(A)
 \end{gather*}
 is the space of Weil classes with respect to the embedding $\Id \oplus \iota: k \hookrightarrow \End^0(A)$. It now suffices to show that $\mathscr{W}_{k, l}(A)$ is spanned by algebraic classes. For this, use Proposition \ref{prop:lifting_field} to lift $A$ together with the $k$-action to an abelian fourfold $\widetilde{A}$ in characteristic zero. If $k$ acts on $\widetilde{A}$ with signature (2, 2), then we conclude by the Hodge conjecture for $\widetilde{A}$. Otherwise $k$ acts with signature $(1, 3)$. Since $p$ is nonsplit in $k$, we conclude by Lemma \ref{lemWeilAlg13}. The argument with $\overline{\iota}$ in place of $\iota$ is identical.
 \par
 Now assume that $A_2$ is CM and ordinary, with $k_2 := \End^0(A_2)$. By Goursat's lemma for algebraic groups, we either have $H_l(A) = H_l(A_1) \times H_l(A_2)$, in which case $(V_l(A)_{3,1})^{\mathfrak{h}_l(A)} = 0$, or we have that $H_l(A) \to H_l(A_1)$ is an isogeny. The latter forces $A_1$ to be non-exotic type IV, by the classification of monodromy. We then let $k_1 := \End^0(A_1)$. The isogeny $H_l(A) \to H_l(A_1)$ and the projection $H_l(A) \to H_l(A_2)$ induce an isogeny between the respective central tori $U_{k_1} \otimes_{\bQ} \bQ_l$ of $H_l(A_1)$ and $U_{k_2} \otimes_{\bQ} \bQ_l$ of $H_l(A_2)$. This implies that $k_{1, l} \cong k_{2, l}$, and via the ``$l$-independence'' of $G_l(A)^\circ$, we reduce to the case where $k_1 \cong k_2$. Then we can lift $A$ up to isogeny to $\widetilde{A}:=\widetilde{A}_1\times \widetilde{A}_2\in \mathrm{S}(A)$, where $\widetilde{A}_1$ is as per Theorem~\ref{thm:lift3fold} and $\widetilde{A}_2$ is the canonical lift. It is easy to show that $\Lie G_{\mathrm{B}}({\widetilde{A}}_{\bC})\otimes \bQ_l$ is isomorphic to $\mathfrak{g}_l$, so the Mumford--Tate conjecture holds for $A$.
 \par
    \item Now we finish the theorem. Suppose the largest dimension of a simple factor of $A$ is $N$.  If $N=1$ then we are in case (\ref{it:TT1}). If $N=2$, write the simple abelian surface factor as $A_1$, and write $A\sim A_1\times A_2$. We have following cases:\begin{enumerate}[label=(\alph*)]
        \item $A_1$ is non-CM, then we are in case (\ref{it:TT2}), 
        \item $A_1$ is CM but $A_2$ is simple non-CM, then we are again in case (\ref{it:TT2}),
        \item $A_1$ is CM but $A_2$ splits as product of elliptic curves. If two elliptic curves are all non-CM, then we are in case (\ref{it:TT3}). If one elliptic curve is CM while the other one is non-CM, then we are in case (\ref{it:TT4}). 
    \end{enumerate}
   Finally suppose that $N=3$. If the simple threefold factor is CM, then we are in case (\ref{it:TT4}). If the simple threefold factor is non-CM, then we are in (\ref{it:TT5}). 

   Therefore, all possible cases are exhausted, and the proof of the Tate conjecture for non-CM non-simple abelian fourfolds is complete. 
\end{enumerate}

\end{proof}

\subsection{The Mumford--Tate conjecture}
Let $A$ be an abelian variety over a finitely generated field $L$. As before, let $G_l(A)$ be the $l$-adic monodromy group of $A$ and let $H_l = (G_l(A) \cap \SL(V_l(A)))^\circ$. We give some general conditions under which the $H_l$ of a non-simple $A$ decomposes as a product over the $H_l$s of the factors.
\par
\begin{proposition}\label{prop:prod}
    Let $A_1',...,A_m', A_1'',...,A_n''$ be abelian varieties over $L$. Let
    \begin{align*}
        A' &= \prod\limits_{i=1}^m A_i',\\
        A'' &= \prod\limits_{j=1}^n A_j''.
    \end{align*}
    Assume that for all simple isogeny factors $T$ of $A = A' \times_L A''$, the following hold:
    \begin{enumerate}
        \item $T$ is geometrically simple;
        \item $T$ has odd relative dimension \cite[Definition 1.1]{Lombardo2016}, or $\dim T \le 2$;
        \item $T$ satisfies the minuscule weights conjecture \ref{conj:MinusculeWt}.
    \end{enumerate}
    Also assume that all simple isogeny factors of $A''$ are type IV, that none of the simple isogeny factors of $A'$ are type IV, and that $\Hom_{\ol{L}}(A_i', A_j') = 0$ for $i \neq j$. Then the projection maps induce an isomorphism of algebraic groups over $\mathbb{Q}_l$
    \begin{gather*}
        H_l(A' \times_L A'') \cong \prod\limits_{i=1}^m H_l(A_i') \times_{\mathbb{Q}_l} H_l(A'').
    \end{gather*}
\end{proposition}
\begin{proof}
    When $L$ has characteristic zero, the claim follows from the argument of \cite[Theorem 4.7]{Lombardo2016}, with minor modifications: namely, we may relax the hypothesis in \textit{loc.cit.} that all $A_i'$, $A_j''$ be absolutely simple, using the argument of \cite[Proposition 2.8]{Lombardo2016}. We may further relax the hypothesis on odd relative dimension in \cite[Theorem 4.7]{Lombardo2016} to accommodate surface isogeny factors of even relative dimension, using \cite[Remark 2.25]{Lombardo2016}. 
    \par
    When $L$ has positive characteristic $p$, we similarly follow the argument of \cite[Theorem 5.9]{Lombardo2016}, with minor modifications. First, note that \cite[Section 5]{Lombardo2016} mistakenly states that $G_l$ of a nonzero abelian variety over $L$ may fail to contain the torus of homotheties. In fact, such a $G_l$ always contains the homotheties, by Lemma \ref{lemHomotheties}. Once this error is corrected, the proof of \cite[Theorem 5.9]{Lombardo2016} applies with the above definition of $H_l$ in place of \cite[Definition 5.5]{Lombardo2016}. 
    \par
    As in the case where $L$ has characteristic zero, we may relax the hypothesis that the $A_i'$ and $A_j''$ be absolutely simple in \cite[Theorem 5.9]{Lombardo2016}, and the hypothesis on odd relative dimension. Here we use \cite[Theorem 0.2.1]{Zarhin1980} to address the case of a surface isogeny factor $S$ with geometric endomorphism algebra $\mathbb{Q}$: such an $S$ has even relative dimension. We may also relax the hypothesis on ordinary reduction in \cite[Theorem 5.9]{Lombardo2016}, using our assumption on the minuscule weights conjecture. This completes the proof of the claim.
\end{proof}

\begin{remark}
    Recall that the minuscule weights conjecture was previously known to hold for abelian varieties of dimension $\le 2$ \cite[Remark 0.4.2.1]{Zarhin1985}.
\end{remark}

\begin{proof}[Proof of Theorem~\ref{thm:Mumford_Tate_nonsimple}]
We have seen from Theorem~\ref{thm:simple_leq 3} that a simple non-CM abelian variety $A$ of dimension $g\leq 3$ over $L$ satisfies the Mumford--Tate conjecture, except when $A$ is exotic type IV. We first make the following: 
\begin{claim}
   A non-simple abelian variety $A$ of dimension $g\leq 4$ such that each simple factor is non-CM and not of exotic type IV satisfies the Mumford--Tate conjecture. If a factor is of exotic type IV, then $A$ satisfies the weak Mumford--Tate conjecture. 
\end{claim}
Let $A_1,A_2,\dots, A_n$ ($n\leq 4$) be the simple factors of $A$ which are mutually non-isogenous (of course we allow $A$ to have repeated simple factors, but a repeated factor does not contribute to the monodromy). If $\dim A_i\leq 2$ for all $i$, then no type IV factor appears, so Proposition~\ref{prop:prod} applies to show that \begin{equation}\label{eq:prodddd}
    H_l(A)=\prod_{i=1}^n H_l(A_i).
\end{equation}
Recall that for each $i=1,\dots,n$, $A_i$ satisfies the Mumford--Tate conjecture. Let $\widetilde{A}_i\in\mathrm{S}(A_i)$ be such that $H_{\mathrm{B}}(\widetilde{A}_{i,\bC})\otimes \bQ_l= H_l(A_i)$ (here $H_{\mathrm{B}}(-)$ is the Hodge group). Then we can construct an $\widetilde{A}\in \mathrm{S}(A)$ simply by taking the product of all $\widetilde{A}_i$ with repetition (if a simple factor $A_i$ is repeated $r$ times, then $\widetilde{A}_i$ is also repeated $r$ times). Then $H_{\mathrm{B}}(\widetilde{A}_\bC)\subseteq \prod_{i=1}^4 H_{\mathrm{B}}(\widetilde{A}_{i,\bC})$. On the other hand $H_l(A)\subseteq H_{\mathrm{B}}(\widetilde{A}_\bC)\otimes \bQ_l$. So  
(\ref{eq:prodddd}) forces the two inclusions in the last sentence to be equality. Therefore $A$ satisfies the Mumford--Tate conjecture. 

The case remains is that $A\sim A_1\times A_2$ and $A_1$ is a simple non-CM threefold and $A_2$ is a non-CM elliptic curve. In this case, $H_l(A_2)$ is $\mathrm{SL}_{2,\bQ_l}$ and by the classification in Theorem~\ref{thm:simple_leq 3} together with Goursat's lemma, we find that $H_l(A)= H_l(A_1)\times H_l(A_2)$ (here we use the ``l-independence'' of $G_l(A)^\circ$ to address the case where $A_1$ is type I(3), as in step \ref{it:TT5} of the proof of Theorem \ref{thm:Tate_nonsimple}). We then treat the Betti side in the same way as the last paragraph, and it follows that the Mumford--Tate conjecture holds for $A$ unless $A_1$ is of exotic type IV.

This proves the \textit{Claim}.

Finally, suppose that $A$ has CM factors. Then the factor contributes nothing to the derived group of the monodromy. More precisely, let $A=A_1\times A_2$ where each simple factor of $A_1$ is non-CM and each simple factor of $A_2$ is CM. Then $H_l(A)^{\mathrm{der}}=H_l(A_1)^{\mathrm{der}}$.
Since the weak Mumford--Tate conjecture holds for $A_1$ by the \textit{Claim}, we can construct an $\widetilde{A}\in \mathrm{S}(A)$ by taking $\widetilde{A}_1\in \mathrm{S}(A_1)$ with
$H_l(A_1)^{\der}= H_{\mathrm{B}}(\widetilde{A}_{1,\bC})^{\der}\otimes \bQ_l$, and taking $\widetilde{A}_2\in \mathrm{S}(A_2)$ to be CM. Let $\widetilde{A}=\widetilde{A}_1\times \widetilde{A}_2$. Then $H_l(A)^{\mathrm{der}}= H_l(A_1)^{\mathrm{der}} = H_\mathrm{B}(\widetilde{A}_{\bC})^{\mathrm{der}} \otimes \bQ_l$, by the fact that $\widetilde{A}_2$ and $A_2$ are CM. So the weak Mumford--Tate conjecture holds for $A$.
\end{proof}

\bibliographystyle{alpha}
\bibliography{references}

\end{document}